\documentclass[12pt]{article}
\usepackage{amssymb}
\usepackage{times}

\usepackage{epsfig}

\usepackage{amssymb,amsmath,amsthm}
\usepackage{amsfonts}

\newtheorem{defn}{Definition}[section]

\newtheorem{cor}[defn]{Corollary}

\newtheorem{lem}[defn]{Lemma}

\newtheorem{prop}[defn]{Proposition}

\newtheorem{rk}[defn]{Remark}

\newtheorem{thm}[defn]{Theorem}

\newcommand{\beeqn}{\begin{equation}}
\newcommand{\be}{\begin{equation}}
\newcommand{\ee}{\end{equation}}
\newcommand{\bea}{\begin{eqnarray}}
\newcommand{\eea}{\end{eqnarray}}
\newcommand{\beas}{\begin{eqnarray*}}
\newcommand{\eeas}{\end{eqnarray*}}

\newcommand{\goto}{\rightarrow}

\newcommand{\ink}{\rule{.5\baselineskip}{.55\baselineskip}}

\newcommand{\ds}{\displaystyle}

\newcommand{\noi}{\noindent}

\newcommand{\lan}{\langle}
\newcommand{\ran}{\rangle}

\newcommand{\ve}{\varepsilon}

\newcommand{\skp}{\vspace{\baselineskip}}

\newcommand{\R}{{\mathbb R}}

\newcommand{\N}{\mathbb N}
\newcommand{\Nb}{{\mathbb N}_b}
\newcommand{\nb}{{\mathbb N}_b}
\newcommand{\nbjstar}{{\mathbb N}_b \setminus \{\jstarplusone\}}
\newcommand{\nzero}{{\mathbb N}_0}

\newcommand{\Z}{\mathbb Z}

\newcommand{\alphabc}{\alpha_b(c)}
\newcommand{\alphaonec}{\alpha_1(c)}
\newcommand{\ankm}{A_{N,K,m}}
\newcommand{\anbm}{A_{N,b,m}}

\newcommand{\barthetanb}{\overline{\Theta}_{N,b}}
\newcommand{\bnkm}{B_{N,K,m}}
\newcommand{\bnbm}{B_{N,b,m}}
\newcommand{\betam}{\beta_m}
\newcommand{\gammam}{\gamma_m}

\newcommand{\djay}{D_j}
\newcommand{\dnu}{\mathcal{D}_\nu}

\newcommand{\deltanbmnu}{\Delta_{N,b,m;\nu}}

\newcommand{\fraconen}{\frac{1}{N}}

\newcommand{\gammankmell}{\gamma_{N,K,m;\ell}}
\newcommand{\gbalpha}{g_b(\alpha)}
\newcommand{\hatb}{\widehat{B}}
\newcommand{\hatf}{\overline{F}}

\newcommand{\gammabalpha}{\gamma_b(\alpha)}

\newcommand{\jstar}{j^\star}
\newcommand{\jstarminusone}{j^\star-1}
\newcommand{\jstarplusone}{j^\star+1}
\newcommand{\jstarplustwo}{j^\star+2}
\newcommand{\jstarmminusone}{j^\star+m-1}
\newcommand{\jstarm}{j^\star+m}

\newcommand{\kappam}{\kappa_m}

\newcommand{\kell}{K_\ell}

\newcommand{\Kbarell}{\overline{K}_\ell}
\newcommand{\kbarell}{\overline{k}_\ell}
\newcommand{\lamn}{\Lambda_N}

\newcommand{\Nj}{N_j}
\newcommand{\nj}{N_j}
\newcommand{\nuj}{\nu_j}
\newcommand{\nun}{\nu^{(N)}}
\newcommand{\nunj}{\nu^{(N)}_j}
\newcommand{\nujstar}{\nu_{j^\star}}
\newcommand{\nujstarplusone}{\nu_{j^\star+1}}

\newcommand{\omegankm}{\Omega_{N,K,m}}
\newcommand{\omegan}{\Omega_N}
\newcommand{\omeganbm}{\Omega_{N,b,m}}
\newcommand{\omeganb}{\Omega_{N,b}}

\newcommand{\omeganzerom}{\Omega_{N,0,m}}
\newcommand{\omeganone}{\Omega_{N,1}}
\newcommand{\omeganonem}{\Omega_{N,1,m}}
\newcommand{\Phijstarm}{\Phi(j^\star,m)}
\newcommand{\pn}{\mathcal{P}_\N}
\newcommand{\Pnbm}{P_{N,b,m}}
\newcommand{\pnbm}{P_{N,b,m}}
\newcommand{\Pnb}{P_{N,b}}
\newcommand{\pnb}{\mathcal{P}_{\N_b}}
\newcommand{\pr}{\mathcal{P}_{r}}
\newcommand{\Pnc}{\mathcal{P}_{N,c}}
\newcommand{\pnc}{\mathcal{P}_{\N,c}}
\newcommand{\pnbc}{\mathcal{P}_{\N_b,c}}

\newcommand{\pnbbc}{\mathcal{P}_{\N_b,[b,c]}}
\newcommand{\pnkm}{P_{N,K,m}}
\newcommand{\Pm}{\mathcal{P}_m}
\newcommand{\Rm}{\R^m}

\newcommand{\rhoalphab}{\rho_{b,\alpha}}
\newcommand{\rhoalphac}{\rho_{\alpha(c)}}
\newcommand{\rhoalphabc}{\rho_{b,\alpha_b(c)}}

\newcommand{\rhoalphabj}{\rho_{b,\alpha;j}}
\newcommand{\rhoalphabjstarplusone}{\rho_{b,\alpha;j^\star + 1}}

\newcommand{\rhoalphabcj}{\rho_{b,\alpha_b(c);j}}
\newcommand{\rhoalphaonec}{\rho_{1,\alpha_1(c)}}

\newcommand{\rhobc}{\rho_{b,c}}
\newcommand{\rhobcj}{\rho_{b,c;j}}

\newcommand{\thetan}{\theta^{(N)}}
\newcommand{\thetanlower}{\theta^{(n)}}

\newcommand{\thetanlowerprime}{\theta^{(n^\prime)}}
\newcommand{\thetanlowerprimej}{\theta^{(n^\prime)}_j}
\newcommand{\thetanj}{\theta^{(N)}_j}
\newcommand{\thetaj}{\theta_j}
\newcommand{\Thetankm}{\Theta_{N,K,m}}
\newcommand{\Thetanb}{\Theta_{N,b}}

\newcommand{\Thetanbj}{\Theta_{N,b;j}}

\newcommand{\thetanbnu}{\theta_{N,b,\nu}}

\newcommand{\thetanbnuj}{\theta_{N,b,\nu;j}}
\newcommand{\thetankm}{\theta_{N,K,m}}

\newcommand{\thetajstar}{\theta_{j^\star}}
\newcommand{\thetajstarplusone}{\theta_{j^\star+1}}

\newcommand{\varphiplus}{\varphi^+}
\newcommand{\varphiminus}{\varphi^-}
\newcommand{\D}{\mathcal{D}}
\newcommand{\W}{\mathcal{W}}
\newcommand{\X}{\mathcal{X}}
\newcommand{\Y}{\mathcal{Y}}
\newcommand{\xelln}{X_{\ell,N}}
\newcommand{\mathZ}{\mathcal{Z}}

\newcommand{\zalphab}{Z_b(\alpha)}

\newcommand{\zbalpha}{Z_b(\alpha)}

\newcommand{\zalphabminusone}{Z_{b-1}(\alpha)}
\newcommand{\zbminusonealpha}{Z_{b-1}(\alpha)}

\newcommand{\zbalphac}{Z_b(\alpha_b(c))}

\newcommand{\zalphabc}{Z_b(\alpha_b(c))}

\newcommand{\zbc}{Z_b(c)}

\newcounter{bean}
\newcommand{\benuma}{\setlength{\labelwidth}{.25in}
\begin{list}%
{(\alph{bean})}{\usecounter{bean}}}
\newcommand{\eenuma}{\end{list}}

\def\theequation{\thesection.\arabic{equation}}
\def\theequation{\arabic{section}.\arabic{equation}}
\def\thedefn{\arabic{section}.\arabic{defn}}
\newcommand{\beginsec}{\setcounter{equation}{0}}

\begin{document}


\title{Detailed Large Deviation Analysis \\ 
of a Droplet Model Having a \\
Poisson Equilibrium Distribution}
\author{Richard S.\ Ellis\normalsize{$\,^1$} \vspace{-.1in} \\
\small{rsellis@math.umass.edu} \vspace{-.125in} \\ \\
Shlomo Ta'asan\normalsize{$\,^2$}  \normalsize \vspace{-.1in}\\  
\small{shlomo@andrew.cmu.edu} \vspace{-.125in} \\ \\
\normalsize{$^1$ Department of Mathematics and Statistics} \vspace{-.05in} \\ 
\normalsize{University of Massachusetts} \vspace{-.05in} \\ 
\normalsize{Amherst, MA 01003} \vspace{-.1in}\\ \\
\normalsize{$^2$ Department of Mathematical Sciences} \vspace{-.05in} \\ 
\normalsize{Carnegie Mellon University} \vspace{-.05in} \\ 
\normalsize{Pittsburgh PA 15213}}
\maketitle

\begin{abstract}
\noi 

One of the main contributions of this paper is to illustrate how large deviation theory can be used to determine the equilibrium distribution of a basic droplet model that underlies a number of important models in material science and statistical mechanics. The model is simply defined.
Given $b \in \N$ and $c > b$,
$K$ distinguishable particles are placed, each with equal probability $1/N$,
onto the $N$ sites of a lattice, where the ratio $K/N$, the average number of particles per site, equals $c$. We focus on configurations for which each site is occupied by a minimum of $b$ particles.
The main result is the large deviation principle (LDP), in the limit where 
$K \goto \infty$ and $N \goto \infty$ with $K/N = c$, for a sequence of random, number-density measures, which are the empirical measures of dependent
random variables that count the droplet sizes. The rate function in the LDP is the relative entropy 
$R(\theta | \rho^\star)$, where $\theta$ is a possible asymptotic configuration of the number-density measures and $\rho^\star$ is a Poisson distribution restricted to the set of positive integers
$n$ satisfying $n \geq b$. This LDP reveals that $\rho^*$ is the equilibrium distribution of the number-density measures, 
which in turn implies that $\rho^*$ is the equilibrium distribution of the random variables that count the droplet sizes. 
We derive the LDP via a local large deviation estimate of the probability that the 
number-density measures equal $\theta$ for any probability measure $\theta$ in the range of these random measures. 
\end{abstract}

\noi
{\it American Mathematical Society 2010 Subject Classifications{\em :}}  60F10 (primary), 82B05 (secondary)
\skp

\noi
{\it Key words and phrases{\em :}} large deviation principle, microcanonical ensemble, number-density measures,
relative entropy

\section{Introduction}
\label{section:intro}
\beginsec

This paper contains the material in the companion paper \cite{EllisTaasan1} together with the following:
full details of several routine proofs omitted from \cite{EllisTaasan1}, additional appendices, and extra background information.

These two papers are motivated by a natural and simply stated question. Given $b \in \N$ and $c > b$, $K$ distinguishable particles are placed, each with 
equal probability $1/N$, onto the $N$ sites of a lattice. Under the assumption
that $K/N = c$ and that each site is occupied by a minimum of $b$ particles, what is the equilibrium distribution, as $N \goto \infty$, of the number of particles per site? We prove in 
Corollary \ref{cor:equilibrium} that this equilibrium distribution is a Poisson distribution $\rhoalphabc$ restricted to the set of positive
integers $n$ satisfying $n \geq b$;
the parameter $\alphabc$ is chosen so that the mean of $\rhoalphabc$ equals $c$.  As we explain at the end
of the introduction, this equilibrium distribution has important applications to technologies using sprays and powders.

We answer this question about the equilibrium distribution by first proving a large deviation principle (LDP) for a sequence of random, number-density
measures, which are the empirical measures of a sequence of dependent random variables that count the droplet sizes. This LDP is stated in
Theorem \ref{thm:ldpthetankm}. 
The space for which we prove the LDP is a natural choice, being
the smallest convex subset of probability measures containing the range of the number-density 
measures. Our proof of the LDP avoids general results 
in the theory of large deviations, many of which do not apply because the space for which we prove the LDP is not a complete, separable metric space. 
Our proof is completely self-contained and starts from first principles, using techniques that are familiar in statistical mechanics. 
For example, the proof of the local large deviation estimate in Theorem \ref{thm:mainestimate}, a key step in the proof
of the LDP for the number-density measures, is based on combinatorics, Stirling's formula, and Laplace asymptotics. Our self-contained proof of the LDP perfectly matches the simplicity and elegance of our main result on the equilibrium distribution stated in the preceding paragraph.

In order to define the droplet model and to formulate the LDP for the number-density measures, a standard probabilistic model is introduced. 
We begin as in the first paragraph. Given $b \in \N$ and $c > b$, 
$K$ distinguishable particles are placed, each with equal probability $1/N$,
onto the $N$ sites of 
the lattice $\Lambda_N = \{1,2,\ldots,N\}$. In section \ref{section:model} we also consider the case $b=0$.
The large deviation limit --- or in statistical mechanical terminology, the
thermodynamic limit --- is defined by taking $K \goto \infty$ and 
$N \goto \infty$ with $K/N$ equal to $c$. The ratio $K/N$ equals the average number of particles per site or the average size of a droplet.
The configuration space for the droplet model is the set $\Omega_{N} = \Lambda_N^K$ consisting of all $\omega = (\omega_1, \omega_2, \ldots, \omega_K)$,
where $\omega_i$ denotes the site in $\Lambda_N$ occupied by the $i$'th particle.  The cardinality of $\Omega_{N}$ equals $N^K$. 
Denote by $P_N$ the uniform probability
measure that assigns equal probability $1/N^K$ to each of the $N^K$ configurations $\omega \in \Omega_{N}$. For subsets $A$ of $\Omega_{N}$, 
$P_N(A) = \mbox{card}(A)/N^K$, where card denotes cardinality.

The asymptotic analysis of the droplet model involves the following two 
random variables, which are functions of the configuration $\omega \in \Omega_{N}$: for $\ell \in \Lambda_N$,
$K_\ell(\omega)$ denotes the number of particles occupying the site $\ell$ in the 
configuration $\omega$; for $j \in \N \cup \{0\}$, $\Nj(\omega)$ denotes the number of sites $\ell \in \Lambda_N$ for which 
$K_\ell(\omega) = j$. 

We focus on the subset of $\Omega_{N}$ consisting of all configurations $\omega$ for which every site of $\Lambda_N$
is occupied by at least $b$ particles. Because of this restriction $\nj(\omega)$ is indexed by $j \in \nb = \{n \in \Z : n \geq b\}$.
It is useful to think of each particle as having one unit of mass and of the set of particles
at each site $\ell$ as defining a droplet. With this interpretation, for each configuration $\omega$,
$K_\ell(\omega)$ denotes the mass or size of the droplet at site $\ell$.
The $j$'th droplet class has $\Nj(\omega)$ droplets and mass $j \Nj(\omega)$.  
Because the number of sites in $\Lambda_N$ equals $N$ and the sum of the masses
of all the droplet classes equals $K$, the following conservation laws hold for such configurations:
\be 
\label{eqn:conserveintro}
\sum_{j \in \Nb} \Nj(\omega) = N \ \mbox{ and } \ \sum_{j \in \Nb} j \Nj(\omega) = K .
\ee
In addition, since the total number of particles is $K$, it follows that $\sum_{\ell \in \Lambda_N} K_\ell = K$. These equality
constraints show that the random variables $N_j$ and the random variables $K_\ell$ are not independent.

In order to carry out the asymptotic analysis of the droplet model, we introduce
a quantity $m = m(N)$ that converges to $\infty$ sufficiently slowly with respect to $N$;
specifically, we require that $m(N)^2/N \goto 0$ as $N \goto \infty$. 
In terms of $b$ and $m$ we define the subset $\omeganbm$ of $\Omega_{N}$ consisting of all configurations $\omega$ for which every site of $\Lambda_N$
is occupied by at least $b$ particles and at most $m$ of the quantities $N_j(\omega)$ 
are positive. This second condition is a key technical device that allows us to control the errors in several estimates. 

The random quantities in the droplet model for which we formulate an LDP 
are the number-density measures $\Thetanb$. For $\omega \in \omeganbm$ 
these random probability measures assign to $j \in \Nb$ the probability $N_j(\omega)/N$, which is the number density of 
the $j$'th droplet class. Thus for any subset $A$ of $\Nb$
\[
\Thetanb(\omega,A) = \sum_{j \in \Nb} \Thetanbj(\omega) \delta_j(A) = \sum_{j \in A} \Thetanbj(\omega), \mbox{ where } \Thetanbj(\omega) = \frac{N_j(\omega)}{N}.
\]
Because of the two conservation laws in (\ref{eqn:conserveintro}) and because $K/N = c$, for $\omega \in \omeganbm$, $\Thetanb(\omega)$ 
is a probability measure on $\Nb = \{n \in \Z : n \geq b\}$ having mean 
\[
\sum_{j \in \Nb} j \Thetanbj(\omega) = \frac{1}{N} \sum_{j \in \Nb} j \Nj(\omega) = \frac{K}{N} = c.
\] 
Thus $\Thetanb$ takes values in $\pnbc$, which is defined to be the set of probability measures on $\Nb$ having mean $c$. $\pnbc$ is topologized
by the topology of weak convergence.

The probability measure $\Pnbm$ defining the droplet model is obtained by restricting the uniform measure $P_{N}$ to the set of configurations 
$\omeganbm$. Thus $\Pnbm$ 
equals the conditional probability $P_{N}(\cdot | \omeganbm)$. For subsets $A$ of $\omeganbm$,
$\Pnbm(A)$ takes the form
\[
\Pnbm(A) = \frac{1}{\mbox{card}(\omeganbm)} \cdot \mbox{card}(A).
\]
In the language of statistical mechanics $\Pnbm$ defines a microcanonical ensemble that incorporates the conservation laws for number
and mass expressed in (\ref{eqn:conserveintro}).

A natural question is to determine two equilibrium distributions: the equilibrium distribution
$\rho^\star$ of the number-density measures and 
the equilibrium distribution $\rho^{**} = \sum_{j \in \Nb} \rho^{**}_j \delta_j$ of the droplet-size random variables $K_\ell$. These distributions
are defined by the following two limits: for any $\ve > 0$, any $\ell \in \Lambda_N$, and all $j \in \Nb$ 
\[
\lim_{N \goto \infty} \Pnbm(\Thetanb \in B(\rho^{*},\ve)) \goto 1 \ \mbox{ and } \ \lim_{N \goto \infty} \Pnbm(K_\ell = j) = \rho^{**}_j,
\]
where $B(\rho^{*},\ve)$ denotes the open ball with center $\rho^{*}$ and radius $\ve$ defined with respect to an appropriate metric on $\pnbc$.
We make the following observations concerning these equilibrium distributions.
\begin{enumerate}
  \item The equilibrium distributions $\rho^*$ for $\Thetanb$ and $\rho^{**}$ for $K_\ell$ coincide.
  \item We first determine the equilibrium distribution $\rho^*$ of $\Thetanb$
and then prove that $\rho^*$ is also the equilibrium distribution of $K_\ell$.
  \item As in many models in statistical mechanics, an efficient way to determine the equilibrium distribution $\rho^*$ of $\Thetanb$
is to prove an LDP for $\Thetanb$, which we carry out in Theorem \ref{thm:ldpthetankm}.
\end{enumerate}

The content of Theorem \ref{thm:ldpthetankm} is the following: as $N \goto \infty$ 
the sequence of number-density measures $\Thetanb$
satisfies the LDP on $\pnbc$ with respect to the measures $\Pnbm$. The rate function is the relative entropy
$R(\theta|\rhoalphab)$ of $\theta \in \pnbc$ with respect to the Poisson distribution $\rhoalphab$ on $\Nb$ having components
\[
\rhoalphabj = \frac{1}{\zalphab} \cdot \frac{\alpha^j}{j!} \mbox{ for } j \in \Nb.
\]
In this formula $\zalphab$ is the normalization that makes $\rhoalphab$ a probability measure, and $\alpha$ equals the unique value $\alphabc$ 
for which $\rhoalphabc$ has mean $c$ [Thm.\ \ref{thm:alphabc}(a)].
Using the fact that $R(\theta|\rhoalphabc)$ equals 0 at the unique measure $\theta = \rhoalphabc$, 
we apply the LDP for $\Thetanb$ to conclude in Theorem \ref{thm:equilibrium}
that $\rhoalphabc$ is the equilibrium distribution of $\Thetanb$. Corollary \ref{cor:equilibrium} then implies 
that $\rhoalphabc$ is also the equilibrium distribution of $K_\ell$.

The space $\pnbc$ is the most natural space on which to formulate the LDP for $\Thetanb$ in Theorem \ref{thm:ldpthetankm}. Not only
is $\pnbc$ the smallest convex set of probability measures containing the range of $\Thetanb$ for all $N \in \N$, but also the union
over $N \in \N$ of the range of $\Thetanb$ is dense in $\pnbc$. As we explain in part (a) of Theorem \ref{thm:pnc}, $\pnbc$
is not a complete, separable metric space, a situation that prevents us from applying the many general results in the theory of large deviations
that require the setting of a complete, separable metric space. In our opinion the fact that we avoid using such general results makes our 
self-contained proof of the LDP even more attractive.

The droplet model is defined in section \ref{section:model}. 
Our proof of the LDP for $\Thetanb$ consists of the following three steps, the first of which 
is the topic of section \ref{section:proof1} and the second and third of which are the topics of section \ref{section:proof2}.

\begin{enumerate}
  \item Step 1 is to derive 
the local large deviation estimate in part (b) of Theorem \ref{thm:mainestimate}. This local estimate, one of the centerpieces of the paper,
gives information not available in the LDP for $\Thetanb$, which involves global estimates.
It states that as $N \goto \infty$, for any probability measure $\theta$ in the range of the number-density measure $\Thetanb$
\be 
\label{eqn:localintro}
\frac{1}{N} \log \Pnbm(\Thetanb = \theta) = -R(\theta|\rhoalphabc) + \mbox{o}(1),
\ee
where $\mbox{o}(1)$ is an error term converging to 0 uniformly for all measures $\theta$ in the range of $\Thetanb$. 
Showing that the parameter of the Poisson distribution $\rhoalphabc$ in the local large deviation estimate equals $\alphabc$
is one of the crucial elements of the proof. The proof of the local large deviation estimate involves combinatorics, Stirling's formula, 
and Laplace asymptotics. 
  \item Step 2 is to lift
this local large deviation estimate to the large deviation limit for $\Thetanb$ lying in open balls and certain other subsets of $\pnbc$. 
This is done in Theorem \ref{thm:ldlimitballs} as a consequence of the general formulation given in Theorem \ref{thm:balls} and the approximation
procedure proved in appendix B. 
  \item Step 3 is to lift the large deviation limit for open balls and certain other subsets to the LDP for $\Thetanb$ stated in Theorem 
\ref{thm:ldpthetankm}, thus proving this LDP. This is done by applying the general formulation given in Theorem \ref{thm:ballstoldp}. 
\end{enumerate}

The paper has four appendices. In appendix A we derive properties of the relative entropy needed in a number of our results.
Appendix B is devoted to the proof of the approximation procedure to which we just referred in item 2 above. In appendix
C we prove the existence of the quantity $\alphabc$ that defines the Poisson distribution $\rhoalphabc$ 
and derive a number of properties of this quantity. Our proof of the existence of $\alphabc$ for general $b$ is subtle. This proof should be
contrasted with the straightforward proof of the existence of $\alphabc$ for $b=1$, which is given in Theorem \ref{thm:alphac}.
We now explain the contents of appendix D.
In order to control several errors in our self-contained proof of the LDP, we must introduce the restriction involving the quantity $m = m(N)$
that, as mentioned earlier, 
requires no more than $m$ of the quantities $N_j$ to be positive. This restriction is explained in detail in section \ref{section:model}; 
it is incorporated in the definition (\ref{eqn:omegankm}) of the set of configurations $\omeganbm$
and the definition (\ref{eqn:condprob}) of the microcanonical ensemble $\pnbm$. In appendix D 
we present evidence supporting the conjecture that this restriction can be eliminated.
Eliminating this restriction would enable us to present our results in a more natural form. 

The paper \cite{EllisTaasan2} explores how our work on the droplet model was inspired
by the work of Ludwig Boltzmann on a simple model of a random ideal gas, for which the Maxwell-Boltzmann is the equilibrium distribution.
The form of the Maxwell-Boltzmann distribution can be proved using Sanov's theorem, 
which proves the LDP for the empirical measures of i.i.d.\ random variables \cite[\S 4]{EllisTaasan2}.
As we show just before Corollary \ref{cor:equilibrium},
$\Thetanb$ is the empirical measure of the random variables $K_\ell$. However, Sanov's theorem for empirical measures of i.i.d.\
random variables cannot be applied because the $K_\ell$ are dependent and, since their distributions depend on $N$, they form a triangular array. 
In section 7 of \cite{EllisTaasan2} we explore how 
Sanov's theorem, although not applicable as stated, can be used to give a heuristic motivation of the LDP for $\Thetanb$.

The main application of the results in this paper is to technologies using sprays and powders,
which are ubiquitous in many fields, including agriculture, the chemical and pharmaceutical industries, consumer products, electronics, manufacturing, material science, medicine, mining, paper making, the steel industry, and waste treatment. In this paper we focus on sprays; our theory also applies to powders with only changes in terminology.
The behavior of sprays might be complex depending on various parameters including evaporation, temperature, and viscosity. Our goal here is to consider the simplest model where the only assumption is made on the average size of droplets in the spray. In many situations
it is important to have good control over the sizes of the droplets, which can be translated into properties of probability distributions. The size distributions are important because they determine reliability and safety in each particular application. 

Interestingly, there does not seem to be a rigorous theory that predicts the equilibrium distribution of droplet sizes, analogous to the Maxwell--Boltzmann distribution of energy levels in a random ideal gas 
\cite{MugEva,SelBrz}. Our goal in the present paper is to provide such a theory. We do so by focusing
on one aspect of the problem related to the relative entropy, an approach that characterizes
the equilibrium distribution of droplet sizes as being a Poisson distribution restricted to 
$\Nb$. We expect that this distribution will dominate experimental observations. A full understanding of droplet behavior under dynamic conditions requires treating many other aspects and is beyond the scope of this paper. A comparison of our results with experimental data will appear elsewhere. 
In addition we plan to apply the ideas in this paper to understand the entropy of dislocation networks. 

Because of the length of this paper and its many technicalities, we would like to help the reader by summarizing the main results and explaining how
one proceeds from the local large deviation estimate stated in (\ref{eqn:localintro}) and proved in part (b) of Theorem \ref{thm:mainestimate} to the LDP 
for the number-density measures $\Thetanb$ stated in Theorem \ref{thm:ldpthetankm}. We also summarize the theorems proved in appendices A, B, C, and D.
\begin{itemize}
  \item {\bf Theorem \ref{thm:ldpthetankm}.} This theorem states that the sequence of $\pnkm$-distributions of the
number-density measures $\Thetanb$ on $\pnc$ satisfies the LDP
on $\pnc$ with rate function $R(\theta|\rhoalphabc)$. 
  \item {\bf Theorem \ref{thm:equilibrium}.} In this theorem we identify the Poisson distribution $\rhoalphabc$
as the equilibrium distribution of $\Thetanb$ with respect to $\pnbm$. It is a consequence of Theorem \ref{thm:ldpthetankm}. 
    \item {\bf Corollary \ref{cor:equilibrium}.} The Poisson distribution $\rhoalphabc$ is shown in this corollary to be also
the equilibrium distribution of the droplet-size random variables $K_\ell$ with respect to $\pnbm$. It is a consequence of Theorem \ref{thm:equilibrium}.
    \item {\bf Theorem \ref{thm:pnc}.} This theorem proves a number of properties of two spaces of probability measures that arise in
the large deviation analysis of $\Thetanb$. 
  \item {\bf Theorem \ref{thm:mainestimate}.} In part (a) of this theorem we show that there exists a unique value $\alpha = \alphabc \in (0,\infty)$
for which the measure $\rhoalphabc$ has mean $c$; the components of $\rhoalphabc$ are defined in (\ref{eqn:rhoj}).
In part (b) we prove the local large deviation estimate (\ref{eqn:localintro}). 
  \item {\bf Theorems \ref{thm:ldlimitballs} and \ref{thm:balls}.} Theorem \ref{thm:ldlimitballs} shows how 
to lift the local large deviation estimate in part (b) of Theorem \ref{thm:mainestimate} to the large deviation limit for $\Thetanb$ lying in open
balls and certain other subsets of $\pnbc$. Theorem \ref{thm:ldlimitballs} is derived as a 
consequence of the general formulation stated in Theorem \ref{thm:balls}. 
  \item {\bf Theorem \ref{thm:ballstoldp}.} This theorem is a general formulation that allows us to lift the 
large deviation limit for open balls and certain other subsets in Theorem \ref{thm:ldlimitballs} to the LDP stated in Theorem \ref{thm:ldpthetankm},
thus proving this LDP.
  \item {\bf Theorem \ref{thm:relentropy}}. In this theorem we collect a number of properties of the relative entropy used throughout the paper.
  \item {\bf Theorem \ref{thm:approximate}.} This result is an approximation theorem that
allows us to approximate an arbitrary probability measure $\theta \in \pnbc$
by a sequence of probability measures $\thetan$ in the range of $\Thetanb$
having the following property: the sequence of relative entropies $R(\thetan | \rhoalphabc)$ converges to $R(\theta|\rhoalphabc)$ as $N \goto \infty$.  
This approximation theorem is applied in two key places. First, it allows us to prove
the asymptotic estimate in Lemma \ref{lem:omegankm}, which is a basic ingredient in the proof of the local large deviation
estimate in part (b) of Theorem \ref{thm:mainestimate}. Second, it allows us to lift this local large deviation estimate to the large
deviation limit for open balls and certain other subsets as formulated in Theorem \ref{thm:ldlimitballs}. 
  \item {\bf Theorem \ref{thm:alphabc}.} This theorem studies a number of properties
of the quantity $\alphabc$ that defines the Poisson-type equilibrium distribution $\rhoalphac$.
  \item{\bf Theorem \ref{thm:alphac}.} This theorem studies a number of properties of the quantity $\alphabc$ for $b = 1$.
  \item {\bf Theorems \ref{thm:msqrtn}, \ref{thm:transfer}, and \ref{thm:bender} and Proposition \ref{prop:omegankmagain}.} 
These results address issues related to the constraint involving
the quantity $m = m(N)$ in the definition (\ref{eqn:omegankm}) of the set of configurations $\omeganbm$ and
the definition (\ref{eqn:condprob}) of the microcanonical ensemble $\pnbm$. We discuss how, if we could eliminate this constraint, our results would have a more natural form. Theorem \ref{thm:bender} is based on a deep, classical result on the asymptotic
behavior of Stirling numbers of the second kind.
\end{itemize}

\skp
\noi 
{\bf Acknowledgments.}
The research of Shlomo Ta'asan is supported in part by a grant from the National Science Foundation
(NSF-DMS-1216433). Richard S.\ Ellis thanks Jonathan Machta for sharing his insights into statistical
mechanics and for useful comments on this introduction, 
Luc Rey-Bellet for valuable conversations concerning large deviation theory, and Michael Sullivan for his generous
help with a number of topological issues arising in this paper. We are also grateful to Jonathan Machta 
for suggesting the generalization, explained in section \ref{section:model}, from a minimum of 1 particle
at each site to a minimum of $b$ particles at each site, where $b$ is any positive integer, and for helping us with the proof
of part (a) of Theorem \ref{thm:alphabc}.

\section{Definition of Droplet Model and Main Theorem}
\label{section:model}
\beginsec

After defining the droplet model, we state the main theorem in the paper, Theorem \ref{thm:ldpthetankm}. The content
of this theorem is the LDP for the sequence of random, number-density measures, which are the empirical measures of a sequence of dependent random variables that count the droplet sizes in the model. As we show in Theorem 
\ref{thm:equilibrium} and in Corollary \ref{cor:equilibrium}, the LDP enables us to identify a Poisson distribution as the equilibrium distribution both of the number-density measures and of the droplet-size random variables. Finally, in Theorem \ref{thm:pnc} we prove a number of properties of two spaces of probability measures in terms of which the LDP  for the number-density measures is formulated.

We start by fixing parameters $b \in \N \cup \{0\}$ and $c \in (b,\infty)$. 
The droplet model is defined by a probability measure $P_{N,b}$ parametrized by $N \in \N$ and the nonnegative integer $b$. The measure depends on 
two other positive integers, $K$ and $m$, where $2 \leq m \leq N < K$. Both $K$ and $m$ are functions of $N$ in the large deviation
limit $N \goto \infty$. In this limit --- which is the same as
the thermodynamic limit in statistical mechanics ---
we take $K \goto \infty$ and $N \goto \infty$, where $K/N$, the average number of particles per site, stays equal to $c$. Thus $K = Nc$.
In addition, 
we take $m \goto \infty$ sufficiently slowly by choosing $m$ to be a function $m(N)$ satisfying $m(N) \goto \infty$ and 
$m(N)^2/N \goto 0$ as $N \goto \infty$; e.g., $m(N) = N^\delta$ for some $\delta \in (0,1/2)$. 
Throughout this paper we fix such a function $m(N)$. The parameter $b$ and the function $m = m(N)$ first appear
in the definition of the set of configurations $\omeganbm$ in (\ref{eqn:omegankm}), where these quantities will be explained. 

Because $K$ and $N$ are integers, 
$c$ must be a rational number. This in turn imposes a restriction on the values of $N$ and $K$. If $c$ is a positive
integer, then $N \goto \infty$ along the positive integers
and $K \goto \infty$ along the subsequence $K = cN$. If $c = x/y$, where $x$ and $y$ are positive integers with
$y \geq 2$ and $x$ and $y$ relatively prime, then $N \goto \infty$ along the subsequence $N = yn$ for $n \in \N$ and
$K \goto \infty$ along the subsequence $K = cN = xn$. Throughout this paper, when we write $N \in \N$ or $N \goto \infty$,
it is understood that $N$ and $K$ satisfy the restrictions discussed here.

In the droplet model $K$ distinguishable particles are placed, each with equal probability $1/N$, onto the sites of the lattice $\Lambda_N
= \{1, 2, \ldots, N\}$. This simple description corresponds to a simple probabilistic model.
The configuration space is the set $\Omega_{N} = \lamn^K$ consisting of all sequences
$\omega = (\omega_1, \omega_2, \ldots, \omega_K)$, where $\omega_i \in \lamn$ denotes the site in $\Lambda_N$ occupied by the $i$'th particle.
Let $\rho^{(N)}$ be the measure on $\lamn$ that assigns equal probability $1/N$ to each site in $\lamn$, and let
$P_{N} = (\rho^{(N)})^K$ be the product measure on $\Omega_{N}$ with equal one-dimensional marginals $\rho^{(N)}$. Thus $P_N$ is the uniform
probability measure that assigns equal probability $1/N^K$ to each of the $N^K$ configurations 
$\omega \in \Omega_{N}$; for subsets $A$ of $\Omega_{N}$ we have $P_{N}(A) = \mbox{card}(A)/N^K$,
where card denotes cardinality. 

The asymptotic analysis of the droplet model involves two random variables that we now introduce. Our goal is to prove a large deviation principle (LDP)
for a sequence of random probability measures defined in terms of these random variables. The LDP is
stated in Theorem \ref{thm:ldpthetankm}. 
\begin{itemize}
  \item For $\ell \in \lamn$ and $\omega \in \Omega_{N}$, $K_\ell(\omega)$ denotes the 
  number of particles occupying site $\ell$ in the configuration $\omega$. In other words, $K_\ell(\omega) =
  \mbox{card}\{i \in \{1,2,\ldots,K\} : \omega_i = \ell\}$.
   \item For $j \in \N \cup \{0\}$ and $\omega \in \Omega_{N}$, $N_j(\omega)$ denotes the 
   number of sites $\ell \in \lamn$ for which $K_\ell(\omega) = j$. 
\end{itemize}
The dependence of $K_\ell(\omega)$ and $N_j(\omega)$ on $N$ is not indicated in the notation. Because the distributions of both random variables
depend on $N$, both $K_\ell$ and $N_j$ form triangular arrays.

We now specify the role played by the nonnegative integer $b$, first focusing on the case where $b$ is a positive integer.
The case where $b = 0$ is discussed later. 
For $\omega \in \omegan$, in general there exist sites $\ell \in \Lambda_N$ for which $K_\ell(\omega) = 0$; i.e., sites that are occupied by 0 particles. For this reason the quantity $N_j(\omega)$ just defined is indexed by $j \in \N \cup \{0\}$.
The next step in the definition of the droplet model is to specify a subset $\omeganbm$ 
of configurations $\omega \in \Omega_{N}$ for which every site is occupied by at least $b$ particles and another constraint holds. In the following definition
of $\omeganbm$, $\nb$ denotes the set $\{n \in \Z : n \geq b\}$. Thus $\N_0$ is the set of nonnegative integers.

\begin{enumerate}
  \item Given $b \in \N$, for any configuration $\omega \in \omeganbm$ every site of $\lamn$ is occupied by at least $b$ particles. In other words, for each 
$\ell \in \lamn$ there exists at least $b$ values of $i \in \{1,2,\ldots,K\} $ such that $\omega_i = \ell$. Equivalently, in the configuration $\omega$
and for each $\ell \in \lamn$ we have $K_\ell(\omega) \geq b$. It follows that for $\omega \in \omeganbm$, $N_j(\omega)$ is indexed by $j \in \nb$.
  \item For any configuration $\omega \in \omeganbm$ at most $m$ of the components $N_j(\omega)$ for $j \in \N_b$ are positive.
As specified at the start of this section, $m = m(N) \goto \infty$ and $m(N)^2/N \goto 0$ as $N \goto \infty$.
\end{enumerate}

We denote by $N(\omega)$ the sequence $\{N_j(\omega), j \in \N_b\}$ and define
\[
|N(\omega)|_+ = \mbox{card}\{j \in \N_b : N_j(\omega) \geq 1 \} .
\]
In terms of this notation
\be 
\label{eqn:omegankm}
\omeganbm = \{\omega \in \Omega_{N} : K_\ell(\omega) \geq b \ \forall \ell \in \lamn \ \mbox{ and } \ |N(\omega)|_+ \leq m = m(N)\} .
\ee

Constraint 2, which restricts the number of positive components of $N(\omega)$, is a useful technical device that
allows us to control the errors in several estimates. In appendix D we explain why we impose this constraint
and give evidence supporting the conjecture that this restriction can be eliminated.
Because of the two constraints, the maximum number of particles that can occupy any site is $K-b(N-1) = N(c-b) + b$. It follows
that $N_j(\omega) = 0$ for all $j \geq N(c-b)+b$. 

When $b$ is a positive integer, for each $\omega \in \omeganbm$ each site in $\Lambda_N$ is occupied by at least $b$ particles. 
In this case it is useful to think of each particle as having one unit of mass and of the set of particles
at each site $\ell$ as defining a droplet. With this interpretation, for each configuration $\omega$,
$K_\ell(\omega)$ denotes the mass or the size of the droplet at site $\ell$.
The $j$'th droplet class has $N_j(\omega)$ droplets and mass $j N_j(\omega)$. 
Because the number of sites in $\Lambda_N$ equals $N$ 
and the sum of the masses of all the droplet classes equals $K$, it follows that the quantities $\Nj(\omega)$ 
satisfy the following conservation laws for all $\omega \in \omeganbm$:
\be 
\label{eqn:conserve}
\sum_{j \in \N_b} \Nj(\omega) = N \ \mbox{ and } \ \sum_{j \in \N_b} j \Nj(\omega) = K .
\ee

We now consider the modifications that must be made in these definitions when $b = 0$. 
In this case constraint 1 in the definition of $\omeganbm$ disappears because we allow sites to be occupied by 0 particles, 
and therefore $N_j(\omega)$ is indexed by $j \in \nzero = \N \cup \{0\}$.  On the other hand, we retain constraint 2 in the definition of $\Omega_{N,0,m}$, which 
requires that for any configuration $\omega \in \Omega_{N,0,m}$ at most $m$ of the components $N_j(\omega)$ for $j \in \N_0$ are positive.
In terms of $|N(\omega)|_+$ the definition of $\Omega_{N,0,m}$ becomes
\[
\Omega_{N,0,m} = \{\omega \in \Omega_{N} : |N(\omega)|_+ \leq m = m(N)\}.
\]
Because the choice $b=0$ allows sites to be empty, we lose the interpretation of the set of particles at each site as being a droplet.
However, for $\omega \in \Omega_{N,0,m}$ the two conservation laws (\ref{eqn:conserve}) continue to hold. 

For the remainder of this paper we work with any fixed nonnegative integer $b$.
The probability measure $\Pnbm$ defining the droplet model is obtained by restricting the uniform measure $P_{N}$ to
the set $\omeganbm$. Thus $\Pnbm$ equals the conditional probability $P_N(\cdot | \omeganbm)$. For subsets $A$ of $\omeganbm$, $\Pnbm(A)$ takes the form
\bea
\label{eqn:condprob}
\Pnbm(A) & = & P_{N}(A \, | \, \omeganbm) = \frac{1}{P_{N}(\omeganbm)} \cdot P_{N}(A) \\
\nonumber & = & \frac{1}{\mbox{card}(\omeganbm)} \cdot \mbox{card}(A).
\eea
The second line of this formula follows from the fact that $P_{N}$ assigns equal probability $1/N^K$ to every $\omega \in \omeganbm$.
In the language of statistical mechanics $\Pnbm$ defines
a microcanonical ensemble that incorporates the conservation laws for number and mass expressed in (\ref{eqn:conserve}).

Having defined the droplet model, we introduce the random probability measures 
whose large deviations we will study. For $\omega \in \omeganbm$ these measures are the number-density measures $\Thetanb$
that assign to $j \in \Nb$ the probability $N_j(\omega)/N$. This ratio represents the number density of droplet class $j$.
Thus for any subset $A$ of $\Nb$
\be 
\label{eqn:thetankm}
\Thetanb(\omega,A) = \sum_{j \in \Nb} \Thetanbj(\omega) \delta_j(A) = \sum_{j \in A} \Thetanbj(\omega), \mbox{ where } \Thetanbj(\omega) = \frac{N_j(\omega)}{N}.
\ee
By the two formulas in (\ref{eqn:conserve}) 
\be 
\label{eqn:sunday}
\sum_{j \in \Nb} \Thetanbj(\omega) = 1 \ \mbox{ and } \ \sum_{j \in \Nb} j \Thetanbj(\omega) = \frac{K}{N} = c .
\ee
Thus $\Thetanb(\omega)$ is a probability measure on $\Nb$ having mean $c$. 

We next introduce several spaces of probability measures that arise in the large deviation analysis of the droplet model.
$\pnb$ denotes the set of probability measures on $\Nb = \{n \in \Z : n \geq b\}$. Thus $\theta \in \pnb$ has the form
$\sum_{j \in \Nb} \theta_j \delta_j$, where the components $\theta_j$ satisfy $\theta_j \geq 0$ and 
$\theta(\Nb) = \sum_{j \in \Nb} \theta_j =  1$.
We say that a sequence of measures $\{\theta^{(n)}, n \in \N\}$ in $\pnb$ converges weakly to $\theta \in \pnb$, and write 
$\thetan \Rightarrow \theta$, if for any bounded function $f$ mapping $\Nb$ into $\R$
\[
\lim_{n \goto \infty} \int_{\Nb} f d\thetanlower = \int_{\Nb} f d\theta .
\]
$\pnb$ is topologized by the topology of weak convergence. There is a standard technique for introducing a metric structure on $\pnb$  
for which we quote the main facts. Because $\N$ is a complete, separable metric space with metric $d(x,y) = |x-y|$, 
there exists a metric $\pi$ on $\pnb$ called the Prohorov metric with the following properties:
\begin{itemize}
  \item Convergence with respect to the Prohorov metric is equivalent to weak convergence \cite[Thm.\ 3.3.1]{EthierKurtz}; i.e., $\thetanlower \Rightarrow \theta$ 
if and only if $\pi(\thetanlower,\theta) \goto 0$ as $N \goto \infty$. 
  \item With respect to the Prohorov metric, $\pnb$ is a complete, separable metric space \cite[Thm.\ 3.1.7]{EthierKurtz}. 
\end{itemize}  

We denote by $\pnbc$ the set of measures in $\pnb$ having mean $c$. Thus $\theta \in \pnbc$ has the form
$\sum_{j \in \Nb} \theta_j \delta_j$, where the components $\theta_j$ satisfy $\theta_j \geq 0$, $\sum_{j \in \Nb} \theta_j = 1$,
and $\int_{\N} x \theta(dx) = \sum_{j \in \Nb} j \theta_j = c$. By (\ref{eqn:sunday})
the number-density measures $\Thetanb$ defined in (\ref{eqn:thetankm}) take values in $\pnbc$.

In part (a) of Theorem \ref{thm:pnc} we prove two properties of $\pnbc$: with respect to the Prohorov metric,
$\pnbc$ is a relatively compact, separable subset of $\pnb$; however, $\pnbc$ is not a closed subset 
of $\pnb$ and thus is not a compact subset or a complete metric space. The fact that $\pnbc$ is not a closed subset of $\pnb$
is easily motivated. If $\thetanlower$ is a sequence in $\pnbc$ such that $\thetanlower \Rightarrow \theta$ for some $\theta \in \pnb$,
then some of the mass of $\thetanlower$ could escape to $\infty$, causing $\theta$ to have a mean strictly less than $c$;
an example is given in (\ref{eqn:nonclosed}).
Although $\pnbc$ is the natural space in which to formulate the LDP for $\Thetanb$ in Theorem \ref{thm:ldpthetankm},
the fact that $\pnbc$ is not a closed subset of $\pnb$ gives rise to a number of unique features in the LDP. 

Because $\pnbc$ is not a closed subset of $\pnb$, it is natural to introduce the closure of $\pnbc$ in $\pnb$.
As we prove in part (b) of Theorem \ref{thm:pnc},
the closure of $\pnbc$ in $\pnb$ equals $\pnbbc$, which is the set of measures in $\pnb$ having mean lying in the closed interval $[b,c]$. 
For any $\theta \in \pnb$ the minimum value of the mean of $\theta$ is $b$, which occurs if and only if $\theta = \delta_b$. Being the closure of
the relatively compact, separable metric space $\pnbc$, $\pnbbc$ is a compact, separable metric space with respect to the Prohorov metric.
This space appears in the formulation of the large deviation upper bound in part (c) of Theorem 
\ref{thm:ldpthetankm}. 
 
We next state Theorem \ref{thm:ldpthetankm}, which is the LDP for the sequence of distributions
$\Pnbm(\Thetanb \in d\theta)$ on $\pnbc$ as $N \goto \infty$.
The rate function in the LDP 
is the relative entropy of $\theta$ with respect to a certain measure $\rhoalphabc = \sum_{j \in \Nb} \rhoalphabcj \delta_j$ defined 
in (\ref{eqn:rhoj}), where each $\rhoalphabcj > 0$. 
Thus any $\theta \in \pnbc$ is absolutely continuous with respect to $\rhoalphabc$. For $\theta \in \pnbc$ the relative entropy of $\theta$
with respect to $\rhoalphabc$ is defined by
\be 
\label{eqn:definerelentropy}
R(\theta | \rhoalphabc) = 
\sum_{j \in \Nb} \theta_j \log (\theta_j/\rhoalphabcj) .
\ee
If $\theta_j = 0$, then $\theta_j \log (\theta_j/\rhoalphabcj) = 0$. For $A$ a subset of $\pnbc$
or $\pnbbc$, $R(A | \rhoalphabc)$ denotes the infimum
of $R(\theta | \rhoalphabc)$ over $\theta \in A$.

For $j \in \Nb$ the components of the measure $\rhoalphabc$ appearing in the LDP have the form 
\be 
\label{eqn:rhoj}
\rhoalphabcj = \frac{1}{\zalphabc} \cdot \frac{[\alphabc]^j}{j!} ,
\ee
where $\alphabc \in (0,\infty)$ is chosen so that $\rhoalphabc$ has mean $c$
and $\zalphabc$ is the normalization making $\rhoalphabc$ a probability measure; 
thus $Z_0(\alpha_0(c)) = e^{\alpha_0(c)}$, and for $b \in \N$, $\zalphabc = e^{\alphabc} - \sum_{j=0}^{b-1} [\alphabc]^j/j!$. 
As we show in part (a) of Theorem \ref{thm:alphabc}, there exists a unique value of 
$\alphabc$. For $b \in \N$
the Poisson-type distribution $\rhoalphabc$ differs
from a standard Poisson distribution because the former has 0 mass at $0, 1, \ldots, b-1$ while the latter has positive mass at these points.
In fact, $\rhoalphabc$ can be identified as the distribution of a Poisson random variable $\Xi_{\alphabc}$ with parameter $\alphabc$ 
conditioned on $\Xi_{\alphabc} \in \Nb$ [Thm.\ \ref{thm:alphabc}(d)]. Despite this difference we shall also refer to $\rhoalphabc$ as a Poisson distribution.

According to part (a) of Theorem \ref{thm:ldpthetankm} $R(\cdot|\rhoalphabc)$ has compact
level sets in $\pnbc$. It is well known that the relative entropy has compact level sets in the complete space $\pnb$. The level
sets are also compact in $\pnbbc$ because the latter is a compact subset of $\pnb$. However,
because $\pnbc$ is not closed in $\pnb$, the compactness of the level sets in $\pnbc$ is not obvious.

As a consequence of the fact that $\pnbc$ is not closed in $\pnb$, 
the large deviation upper bound takes two forms depending on whether the subset $F$ of $\pnbc$ is compact or 
whether $F$ is closed. When $F$ is compact, in part (b)
we obtain the standard large deviation upper bound for $F$ with $-R(F | \rhoalphabc)$ on the
right hand side. When $F$ is closed, in part (c) we obtain a variation of the standard large deviation upper bound; $-R(F | \rhoalphabc)$ on the right hand side is
replaced by $-R(\overline{F} | \rhoalphabc)$, where $\overline{F}$ is the closure of $F$ in the compact space $\pnbbc$ and is therefore compact.
When $F$ is compact, its closure in $\pnbbc$ is $F$ itself. In this case the large deviation upper bounds in parts (b) 
and (c) coincide.

The refinement in part (c) is important. It is applied in the proof of Theorem \ref{thm:equilibrium} 
to show that $\rhoalphabc$ is the equilibrium distribution of the number-density measures $\Thetanb$. In turn, Theorem \ref{thm:equilibrium}
is applied in the proof of Corollary \ref{cor:equilibrium} to show that $\rhoalphabc$ is the equilibrium distribution of the droplet-size 
random variables $K_\ell$.

In the next theorem we assume that $m$ is the function $m(N)$ appearing in the definition of $\omeganbm$ in (\ref{eqn:omegankm})
and satisfying $m(N) \goto \infty$ and $m(N)^2/N \goto 0$ as $N \goto \infty$.
The assumption that $m(N)^2/N \goto 0$ is used to control error terms in Lemmas \ref{lem:deltankmnu}, \ref{lem:omegankm}, and \ref{lem:nujstar}.
This assumption on $m(N)$ is optimal in the sense that it is a minimal assumption guaranteeing
that an error term in the lower bound in part (a) of Lemma \ref{lem:nujstar} and in the upper bound in part (b) of the lemma converge to 0.

\begin{thm}
\label{thm:ldpthetankm}
Fix a nonnegative integer $b$ and
a rational number $c \in (b,\infty)$. Let $m$ be the function $m(N)$ appearing in the definition of $\omeganbm$ in {\em (\ref{eqn:omegankm})}
and satisfying $m(N) \goto \infty$ and $m(N)^2/N \goto 0$ as $N \goto \infty$.
Let $\rhoalphabc \in \pnbc$ be the distribution having the components defined in {\em (\ref{eqn:rhoj})}.
Then as $N \goto \infty$, with respect to the measures $\Pnbm$,
the sequence $\Thetanb$ satisfies the large deviation principle on $\pnbc$ with rate function $R(\theta|\rhoalphabc)$ in the following sense.

{\em (a)} $R(\theta|\rhoalphabc)$ maps $\pnbc$ into $[0,\infty]$, and for any $M< \infty$ the level set $\{\theta \in \pnbc : R(\theta | \rhoalphabc) \leq M\}$ is compact.

{\em (b)} For any compact subset $F$ of $\pnbc$ we have the large deviation upper bound
\[
\limsup_{N \goto \infty} \frac{1}{N} \log \Pnbm(\Thetanb \in F) \leq -R(F | \rhoalphabc).
\]

{\em (c)} For any closed subset $F$ of $\pnbc$, let $\overline{F}$ denote the closure of $F$ in $\pnbbc$. We have the large deviation upper bound
\[
\limsup_{N \goto \infty} \frac{1}{N} \log \Pnbm(\Thetanb \in F) \leq -R(\overline{F} | \rhoalphabc).
\]

{\em (d)} For any open subset $G$ of $\pnbc$ we have the large deviation lower bound
\[
\liminf_{N \goto \infty} \frac{1}{N} \log \Pnbm(\Thetanb \in G) 
\geq -R(G | \rhoalphabc) .
\]
\end{thm}

As noted in the comments after the statement of Theorem \ref{thm:ballstoldp},  
Theorem \ref{thm:ldpthetankm} is a consequence of that theorem and several other results 
proved in the paper. Part (b) of Theorem \ref{thm:mainestimate} proves a local large
deviation estimate for probabilities of the form $\Pnbm(\Thetanb = \theta)$, where $\theta$ is a probability measure
in the range of $\Thetanb$. This local estimate is one of the centerpieces of this paper,
giving information not available in the LDP for $\Thetanb$, which involves global estimates.
In Theorem \ref{thm:ldlimitballs} we show how to lift this local estimate to the large deviation
limit for $\Thetanb$ lying in open balls and certain other subsets of $\pnbc$ defined in terms of open balls. 
Theorem \ref{thm:ldlimitballs} is proved as an application of the
general formulation given in Theorem \ref{thm:balls}. Finally we show how to lift the large deviation limit for open balls 
and certain other subsets defined in terms of open balls to the LDP stated in Theorem
\ref{thm:ldpthetankm}. We do so by applying the general formulation given in Theorem \ref{thm:ballstoldp}.
In part (d) of Theorem \ref{thm:relentropy} we prove
that the level sets of $R(\theta | \rhoalphabc)$ in $\pnbc$ are compact. 

The rate function in Theorem \ref{thm:ldpthetankm}
has the property that for $\theta \in \pnbbc$, $R(\theta | \rhoalphabc) \geq 0$ with equality if and only
if $\theta = \rhoalphabc$ [Thm.\ \ref{thm:relentropy}(a)]. As we explain in the next theorem, the large deviation upper
bound and this property of the relative entropy allow us to interpret the Poisson distribution
$\rhoalphabc$ as the equilibrium distribution of the number-density measures $\Thetanb$. In this theorem 
$[B_\pi(\rhoalphabc,\ve)]^c$  denotes the complement in $\pnbc$ of 
the open ball in $\pnbc$ with center $\rhoalphabc$ and radius $\ve > 0$ with respect to the Prohorov metric $\pi$. This open ball is defined by
\[
B_\pi(\rhoalphabc,\ve) = \{\nu \in \pnbc : \pi(\rhoalphabc,\nu) < \ve\}.
\]
$[\widehat{B}_\pi(\rhoalphabc,\ve)]^c$ denotes the complement in $\pnbbc$ of the open ball defined by 
\[
\widehat{B}_\pi(\rhoalphabc,\ve) = \{\nu \in \pnbbc : \pi(\rhoalphabc,\nu) < \ve\} .
\]

There is a subtlety in the proof in the next theorem that $\rhoalphabc$ is the equilibrium distribution of $\Thetanb$.
To prove this, we need an exponentially decaying estimate on the probability that $\Thetanb \in [B_\pi(\rhoalphabc,\ve)]^c$.
Since $[B_\pi(\rhoalphabc,\ve)]^c$ is closed in $\pnbc$ but is not compact, we obtain this estimate by applying the large deviation upper bound in part (c)
of Theorem \ref{thm:ldpthetankm} to $[B_\pi(\rhoalphabc,\ve)]^c$ and using the fact that the closure of this set in $\pnbbc$ is a subset
of $[\widehat{B}_\pi(\rhoalphabc,\ve)]^c$. 

\begin{thm}
\label{thm:equilibrium}
We assume the hypotheses of Theorem {\em \ref{thm:ldpthetankm}}. The following results hold for any $\ve > 0$. 

{\em (a)} The quantity
$x^\star = \inf\{R(\theta | \rhoalphabc) : \theta \in [\widehat{B}_\pi(\rhoalphabc,\ve)]^c\}$ is strictly positive.

{\em (b)} For any number $y$ in the interval $(0,x^\star)$ and all sufficiently large $N$
\[
\Pnbm(\Thetanb \in [B_\pi(\rhoalphabc,\ve)]^c) \leq \exp[-N y] \ \mbox{ as } N \goto \infty .
\]
This upper bound implies that as $N \goto \infty$
\[
\lim_{N \goto \infty} \Pnbm(\Thetanb \in B_\pi(\rhoalphabc,\ve)) = 1 \ \mbox{ and } \ 
\lim_{\ve\goto 0}\lim_{N \goto \infty} \Pnbm(\Thetanb \in B_\pi(\rhoalphabc,\ve)) = 1 .
\]
These limits allow us to interpret the Poisson distribution $\rhoalphabc$ having the components defined in {\em (\ref{eqn:rhoj})}
as the equilibrium distribution of the number-density measures $\Thetanb$ with respect to $\Pnbm$.
\end{thm}

\noi 
{\bf Proof.}  The starting point is the large deviation upper bound in part (c) of Theorem \ref{thm:ldpthetankm} applied
to the closed set $[B_\pi(\rhoalphabc,\ve)]^c$, which is a subset of $[\widehat{B}_\pi(\rhoalphabc,\ve)]^c$. 
We denote the closure of $[B_\pi(\rhoalphabc,\ve)]^c$ in $\pnbbc$ by $\overline{[B_\pi(\rhoalphabc,\ve]^c}$. 
We claim that $\overline{[B_\pi(\rhoalphabc,\ve)]^c} \subset [\widehat{B}_\pi(\rhoalphabc,\ve)]^c$. Indeed, any $\nu \in \overline{[B_\pi(\rhoalphabc,\ve]^c}$
is the weak limit of a sequence $\nu^{(n)} \in [B_\pi(\rhoalphabc,\ve]^c \subset \pnbc$. Since the closure of $\pnbc$ in $\pnb$ equals $\pnbbc$, in
general we have $\nu \in \pnbbc$. In addition, since $\nu^{(n)} \in [\widehat{B}_\pi(\rhoalphabc,\ve)]^c$, it follows that 
$\nu \in [\widehat{B}_\pi(\rhoalphabc,\ve)]^c$. 
This proves the claim that 
$\overline{[B_\pi(\rhoalphabc,\ve)]^c} \subset [\widehat{B}_\pi(\rhoalphabc,\ve)]^c$.
Because of this relationship, the large deviation upper bound in part (c) of Theorem \ref{thm:ldpthetankm} takes the form
\bea
\label{eqn:lefteqn}
\lefteqn{
\limsup_{N \goto \infty} \frac{1}{N} \log \Pnbm(\Thetanb \in [B_\pi(\rhoalphabc,\ve)]^c\}} \\
\nonumber
&& \hspace{.5in} \leq -R(\overline{[B_\pi(\rhoalphabc,\ve)]^c} | \rhoalphabc) \leq -R([\widehat{B}_\pi(\rhoalphabc,\ve)]^c | \rhoalphabc).
\eea

We now prove part (a) of Theorem \ref{thm:equilibrium}. 
Since $R(\theta | \rhoalphabc)$ has compact level sets in $\pnbbc$, it attains its infimum $x^\star$ on the closed set $[\widehat{B}_\pi(\rhoalphabc,\ve)]^c$. 
If $x^\star = 0$, then there would exist $\theta \in [\widehat{B}_\pi(\rhoalphabc,\ve)]^c$ such that $R(\theta | \rhoalphabc) = 0$. But
on $\pnbbc$, $R(\theta | \rhoalphabc)$ attains its infimum of 0
at the unique measure $\theta = \rhoalphabc$. Hence we obtain a contradiction because $\rhoalphabc 
\not \in [\widehat{B}_\pi(\rhoalphabc,\ve)]^c$. This completes the proof of 
part (a). The inequality in part (b) is an immediate consequence of part (a) and the large deviation upper bound (\ref{eqn:lefteqn}). 
This inequality yields the two limits in the next display. The proof of Theorem \ref{thm:equilibrium} is complete. \ \ink

\skp
We now apply Theorem \ref{thm:equilibrium} to prove that $\rhoalphabc$ is also the equilibrium distribution of the random variables $K_\ell$,
which count the droplet sizes at the sites of $\Lambda_N$. Although these random variables are identically distributed, they are dependent
because for each $\omega \in \omeganbm$ they satisfy the equality constraint $\sum_{\ell \in \Lambda_N} K_\ell(\omega) = K$. 
Except for one step the proof that $\rhoalphabc$ is also the equilibrium distribution of $K_\ell$ is completely algebraic and requires only the condition that the $K_\ell$ are identically distributed. Their dependence does not affect the proof.
A key observation needed in the proof is that $\Thetanb$ is the empirical measure of these random variables; i.e.,
for $\omega \in \omeganbm$, $\Thetanb(\omega)$ assigns to subsets $A$ of $\Nb$ the probability
\[
\Thetanb(\omega,A) = \frac{1}{N} \sum_{\ell=1}^N \delta_{K_\ell(\omega)}(A).
\]
This characterization of $\Thetanb$ follows from the fact that the empirical measure of $K_\ell$ assigns to $j \in \Nb$ the probability
\be 
\label{eqn:empirical}
\frac{1}{N} \sum_{\ell=1}^N \delta_{K_\ell(\omega)}(\{j\}) = \frac{N_j(\omega)}{N} = \Thetanbj(\omega).
\ee

\begin{cor}
\label{cor:equilibrium}
We assume the hyotheses of Theorem {\em \ref{thm:ldpthetankm}}. Then for any site $\ell \in \Lambda_N$ and any $j \in \Nb$
\[ 
\lim_{N \goto \infty} \Pnbm(K_\ell = j) = \rhoalphabcj = \frac{1}{\zbalphac} \cdot \frac{[\alphabc]^j}{j!}.
\]
\end{cor}

\noi 
{\bf Proof.} Since the random variables $K_\ell$ are identically distributed, it suffices to prove the corollary for $\ell = 1$. 
Theorem \ref{thm:equilibrium} implies that if $g$ is any bounded continuous function mapping $\pnbc$ into $\R$, then 
\be 
\label{eqn:nblim}
\lim_{N \goto \infty} \int_{\omeganbm} g(\Thetanb) d\Pnbm = g(\rhoalphabc).
\ee
Given $\varphi$ any bounded function mapping $\Nb$ into $\R$ we define for $\theta \in \pnb$ the bounded function
\[
g(\theta) = \sum_{j \in \nb} \varphi(j) \thetaj.
\]
By the definition of weak convergence, $g$ is continuous on $\pnbc$. 
Equation (\ref{eqn:empirical}) now yields
\beas
g(\Thetanb(\omega)) & = & \sum_{j \in \nb} \varphi(j) \Thetanbj (\omega) \\
& = & \frac{1}{N} \sum_{\ell \in \Lambda_N} \sum_{j \in \Nb} \varphi(j) \delta_{K_\ell(\omega)}(\{j\})
= \frac{1}{N} \sum_{\ell \in \Lambda_N} \varphi(K_\ell(\omega)).
\eeas
Since the $K_\ell$ are identically distributed, it follows from (\ref{eqn:nblim}) that 
\beas
\lefteqn{
\lim_{N \goto \infty} \int_{\omeganbm} \varphi(K_1) d\Pnbm} \\
& & = \lim_{N \goto \infty} \frac{1}{N} \sum_{\ell=1}^N \int_{\omeganbm} \varphi(K_\ell) d\Pnbm \\
& & = \lim_{N \goto \infty} \int_{\omeganbm} g(\Thetanb) d\Pnbm = g(\rhoalphabc) = \sum_{j \in \nb} \varphi(j) \rhoalphabcj.
\eeas
Setting $\varphi = 1_{j'}$ for any $j' \in \Nb$ yields 
\[
\lim_{N \goto \infty} \Pnbm(K_1 = j') = \rho_{b,\alphabc;j'}.
\]
This completes the proof of the corollary. \ink

\skp
The last theorem in this section proves several properties of $\pnbc$ and $\pnbbc$ with respect to the Prohorov metric
that are needed in the paper.

\begin{thm}
\label{thm:pnc} Fix a nonnegative integer $b$ and a real number $c \in (b,\infty)$. The metric spaces $\pnbc$ and $\pnbbc$ have the following properties.

{\em (a)} $\pnbc$, the set of probability measures on $\N_b$ having mean $c$, is a relatively compact, separable subset of $\pnb$. However, $\pnbc$ is not a closed
subset of $\pnb$ and thus is not a compact subset or a complete metric space.

{\em (b)} $\pnbbc$, the set of probability measures on $\N_b$ having mean lying in the closed interval $[b,c]$, is the closure of $\pnbc$ in $\pnb$. $\pnbbc$
is a compact, separable subset of $\pnb$. 
\end{thm}

\noi 
{\bf Proof.} 
(a) For $\xi \in \N$ satisfying $\xi \geq b$ let $\Psi_\xi$ denote the compact subset $\{b,b+1,\ldots,\xi\}$ of $\Nb$, and let
$[\Psi_\xi]^c$ denote its complement. For any $\theta \in \pnbc$
\[ 
c = \sum_{j \in \Nb} j \thetaj \geq \sum_{j \geq \xi + 1} j \thetaj \geq \xi \sum_{j \geq \xi + 1} \thetaj = \xi \theta([\Psi_\xi]^c) .
\]
It follows that $\pnbc$ is tight; i.e., for any $\ve > 0$ there exists $\xi \in \N$ such that
\[
\sup_{\theta \in \pnbc} \theta([\Psi_\xi]^c) < \ve .
\]
Prohorov's Theorem implies that $\pnbc$ is relatively compact \cite[Thm.\ 3.2.2]{EthierKurtz}. The separability
of $\pnbc$ is proved in Corollary \ref{cor:dense}. 

In the present setting the relative compactness of $\pnbc$ is easy to prove from the tightness of $\pnbc$ without
appealing to the general formulation of Prohorov's Theorem. Given any sequence $\thetanlower \in \pnbc$, a diagonal
argument yields a subsequence $\thetanlowerprime$ such that $\thetaj = \lim_{n \goto \infty} \thetanlowerprimej$ exists for all $j \in \Nb$.
Define $\theta = \sum_{j \in \Nb}\thetaj \delta_j$. We claim that $\thetanlowerprime \Rightarrow \theta$. To see this let $f$ be any nonzero bounded
function mapping $\Nb$ into $\R$. Given $\ve > 0$ choose $\xi \in \Nb$ so large that 
\[
\sup_{n'} \thetanlowerprime([\Psi_\xi]^c) < \ve/[2\|f\|_\infty] \ \mbox{ and } \ 
\theta([\Psi_\xi]^c) < \ve/[2\|f\|_\infty].
\]
The latter bound is possible since by Fatou's Lemma $c = \liminf_{n' \goto \infty} \sum_{j \in \Nb}j \thetanlowerprimej \geq \sum_{j \in \Nb} j \thetaj$.
It follows that
\bea
\label{eqn:nprime}
\left|\int_{\Nb} f d\thetanlowerprime - \int_{\Nb} f d\theta\right| & \leq & \sum_{j=b}^\xi |f(j)| |\thetanlowerprimej - \theta_j| 
+ \sum_{j \geq \xi+1} |f(j)|(\thetanlowerprimej + \theta_j) \\
\nonumber & \leq & \sum_{j=b}^\xi |f(j)| |\thetanlowerprimej - \theta_j| + \ve.
\eea
Since $\thetanlowerprimej \goto \thetaj$ for $j \in \{b,b+1,\ldots,\xi\}$ and $\ve > 0$ is arbitrary, 
the weak convergence of $\thetanlowerprime$ to $\theta$ is proved.
Taking $f$ to be identically 1 verifies that $\theta \in \pnb$, which must be the case since $\pnb$ is complete.

We now prove that $\pnbc$ is not a closed subset of $\pnb$ by exhibiting a sequence $\thetanlower \in \pnbc$
having a weak limit that does not lie in $\pnbc$.  To simplify the notation, 
we denote the mean of $\sigma \in \pnb$ by $\langle\sigma\rangle$. Let $\theta$ be any measure in $\pnb$ with mean 
$\langle \theta \rangle = \beta \in [b,c)$; thus $\theta \not \in \pnbc$. The sequence
\be 
\label{eqn:nonclosed}
\thetanlower = \frac{n-c}{n-\beta} \theta + \frac{c-\beta}{n-\beta} \delta_n \mbox{ for } n \in \N, n > c
\ee
has the property that $\thetanlower \in \pnbc$ and that $\thetanlower \Rightarrow \theta \not \in \pnbc$. We conclude
that $\pnbc$ is not a closed subset of $\pnb$. This completes the proof of part (a). 

(b) Since $\pnbc$ is a separable subset of $\pnb$ and $\pnbc$ is dense in $\pnbbc$, it follows that $\pnbbc$ is separable.
We prove that $\pnbbc$ is the closure of $\pnbc$ in $\pnb$. 
Let $\thetanlower$ be a sequence in $\pnbc$ converging weakly to $\theta \in \pnb$. Since $\thetanlower \Rightarrow \theta$ implies that 
$\thetanlower_j \goto \theta_j$
for each $j \in \Nb$, Fatou's Lemma implies that
\[
c = \liminf_{n \goto \infty} \langle\thetanlower\rangle \geq \langle\theta\rangle.
\]
Since for any $\theta \in \pnb$ we have $\langle\theta\rangle \geq b$, 
it follows that $c \geq \langle\theta\rangle \geq b$. This shows that the 
closure of $\pnbc$ in $\pnb$ is a subset of $\pnbbc$. 

We next prove that $\pnbbc$ is a subset of the closure of $\pnbc$ in $\pnb$ by showing that for
any $\theta \in \pnbbc$ there exists a sequence $\thetanlower \in \pnbc$ such that $\thetanlower \Rightarrow \theta$. 
If $\langle\theta\rangle = c$, then 
we choose $\thetanlower = \theta$ for all $n \in \N$. If $\langle\theta\rangle = \beta \in [b,c)$, then we use the sequence
$\thetanlower$ in (\ref{eqn:nonclosed}), which converges weakly to $\theta$. 
We conclude
that $\theta$ lies in the closure of $\pnbc$ and thus that $\pnbbc$ is a subset of the closure of $\pnbc$ in $\pnb$. This completes the proof of part (b). 
The proof of Theorem \ref{thm:pnc} is done. \ink

\skp
We end this section by giving examples of closed, noncompact subsets of $\pnbc$ and compact subsets of $\pnbc$. We do this to emphasize 
the care that must be taken in dealing with the non-closed metric space $\pnbc$ and the necessity of having
separate large deviation upper bounds for compact sets in part (b) of Theorem \ref{thm:ldpthetankm} and for closed sets in 
part (c) of Theorem \ref{thm:ldpthetankm}. We construct these examples as level sets of lower semicontinuous functions $I$ mapping $\pnbc$ into $[0,\infty]$
and having the form 
\[
I(\theta) = \int_{\Nb} g d\theta = \sum_{j \in \Nb} g(j) \theta_j, \mbox{ where } g(j) \geq 0 \mbox{ for all } j \in \Nb.
\]
Since $\thetanlower \Rightarrow \theta \in \pnbc$ implies that $\thetanlower_j \goto \theta_j$ for each $j \in \Nb$, Fatou's Lemma implies that
$I$ is lower semicontinuous on $\pnbc$. Thus for any $M < \infty$ the level set
\[
U_M = \{\theta \in \pnbc : I(\theta) \leq M\}
\]
is closed in $\pnbc$. 

For the next set of examples, we assume that $g$ is a nondecreasing function mapping $\Nb$ into $[0,\infty)$ 
and satisfying $g(j) \goto \infty$
and $g(j)/j \goto 0$ as $j \goto \infty$. In this case, as in the proof of part (a) of Theorem
\ref{thm:pnc} that $\pnbc$ is relatively compact, Prohorov's Theorem implies that the level set $U_M$ is relatively compact. However, in general $U_M$
is not compact because it is not closed in $\pnb$. A sequence showing that $U_M$ is not closed in $\pnb$
is given by $\thetanlower \in \pnbc$ defined in (\ref{eqn:nonclosed}), where $\theta$ has mean $\beta \in [b,c)$. 
For all sufficiently large $n$, $\thetanlower$ lies in the level set $U_{\beta+1}$, but $\thetanlower \Rightarrow \theta$, which is not in $\pnbc$.  

For the final set of examples, we assume that $g$ is a nondecreasing function mapping $\Nb$ into $[0,\infty)$
and satisfying 
$g(j)/j \goto \infty$ as $j \goto \infty$. Again Prohorov's Theorem implies 
that $U_M$ is relatively compact. In addition, because of the assumption on $g$, $U_M$ is uniformly integrable; i.e.,
\[
\lim_{D \goto \infty} \sup_{\theta \in U_M} \int_{\{x \in \Nb : x \geq D\}} x \theta(dx) = 0.
\]
This implies that if $\thetanlower \in U_M$ converges weakly to $\theta \in \pnb$, then $c = \langle \thetanlower \rangle \goto \langle \theta \rangle$.
This standard consequence of uniform integrability, proved in Proposition 2.3 in the appendix of \cite{EthierKurtz}, can be proved in the present
setting as in (\ref{eqn:nprime}) if $\thetanlowerprime$ is replaced by $\thetanlower$ and $f(j)$ is replaced by $j$ for $j \in \Nb$.
It follows that $\theta$ has mean $c$ and so lies in $\pnbc$ and therefore in $U_M$ because $U_M$ is closed in $\pnbc$.
We conclude that $U_M$ is both relatively compact and closed in $\pnbc$, implying that $U_M$ is compact.

The rate function in Theorem \ref{thm:ldpthetankm} is the relative entropy $R(\theta|\rhoalphac)$, a lower semicontinuous
function mapping $\pnbc$ into $[0,\infty]$ that does not have the simple form of $I$. 
The proof that $R(\cdot|\rhoalphac)$ has compact level sets in $\pnbc$ relies on 
Lemma 5.1 in \cite{DonVar3} and
the fact that $\rhoalphac$ has a finite moment generating function $\int_{\Nb} \exp(wx) \rhoalphac(dx)$ for all $w \in (0,\infty)$
[Thm.\ \ref{thm:relentropy}(d)].

In the next section we present the local large deviation estimate that will be used in section \ref{section:proof2} to prove the LDP
for $\Thetanb$ in Theorem \ref{thm:ldpthetankm}.

\section{Local Large Deviation Estimate Yielding Theorem \ref{thm:ldpthetankm}}
\label{section:proof1}
\beginsec

The main result needed to prove the LDP in Theorem \ref{thm:ldpthetankm} is the local large deviation estimate
stated in part (b) of Theorem \ref{thm:mainestimate}. The first step is to introduce a set $\anbm$ that plays a central role in 
this paper. Fix a nonnegative integer $b$ and a rational number $c \in (b,\infty)$. Given $N \in \N$ define $K = Nc$ and let
$m$ be the function appearing in the definition of $\omeganbm$ in (\ref{eqn:omegankm})
and satisfying $m(N) \goto \infty$ and $m(N)^2/N \goto 0$ as $N \goto \infty$. Define $\nb = \{n \in \Z : n \geq b\}$; thus 
$\N_0$ is the set of nonnegative integers.
Let $\nu$ be a sequence $\{\nuj, j \in \Nb\}$ for which each $\nuj \in \N_0$; thus $\nu \in \N_0^{\Nb}$.
We define $\anbm$ to be the set of $\nu \in \N_0^{\nb}$ satisfying 
\be 
\label{eqn:nuj3}
\sum_{j\in\Nb} \nu_j = N, \ \sum_{j\in\Nb} j \nu_j = K, \ \mbox{and} \ |\nu|_+ \leq m = m(N),
\ee
where $|\nu|_+ = \mbox{card}\{j \in \Nb : \nuj \geq 1\}$. Because $\nuj \in \N_0$, the two sums 
involve only finitely many terms. 

For $\omega \in \omeganbm$ the 
components $\Thetanbj(\omega)$ of the number-density measure defined in (\ref{eqn:thetankm}) are $\Nj(\omega)/N$ for $j \in \Nb$,
where $N_j(\omega)$ denotes the number of sites in $\Lambda_N$ containing $j$ particles in the configuration $\omega$. 
We denote by $N(\omega)$ the sequence $\{N_j(\omega), j \in \Nb\}$. 
By definition, for every $\omega \in \omeganbm$ each site $\ell \in \Lambda_N$ is occupied by at least $b$ particles,
and $|N(\omega)|_+ \leq m = m(N)$. It follows that $\anbm$ is the range of $N(\omega)$ for $\omega \in \omeganbm$; the two sums involving $\nuj$ in (\ref{eqn:nuj3}) 
correspond to the two sums involving $N_j(\omega)$ in (\ref{eqn:conserve}).

Since the range of $N(\omega)$ is $\anbm$, for $\omega \in \omeganbm$ the range of $\Thetanb(\omega)$ is the set of probability measures $\thetanbnu$
whose components for $j \in \Nb$ have the form
\be 
\label{eqn:thetankmj}
\thetanbnuj = \frac{\nu_j}{N} \ \mbox{ for } \nu \in \anbm .
\ee
By (\ref{eqn:nuj3}) $\thetanbnu$ takes values in $\pnbc$, the set of probability measures on $\Nb$ having mean $c$.
It follows that the set  
\be 
\label{eqn:bnkm}
\bnbm = \{\theta \in \pnbc : \theta_j = \nu_j/N \mbox{ for } j \in \Nb \mbox{ for some } \nu \in \anbm\} 
\ee
is the range of $\Thetanb(\omega)$ for $\omega \in \omeganbm$. 

In part (b) of the next theorem we state the local large deviation estimate for the event $\{\Thetanb = \thetanbnu\}$.
In part (a) we introduce the Poisson distribution $\rhoalphabc$ that appears in the local estimate. This Poisson distribution is the 
restriction to $\Nb$ of a standard Poisson distribution on $\N \cup \{0\}$; $\rhoalphabc$ is defined
in terms of a parameter $\alphabc$ guaranteeing that it has mean $c$. If $b=0$, then $\alpha_0(c) = c$, while if $b \in \N$, then $\alphabc < c$
[Thm.\ \ref{thm:alphabc}(b)]. 

In Theorem \ref{thm:alphac} we give the straightforward proof of the existence of $\alphabc$ for $b=1$. 
The proof of the existence
of $\alphabc$ for general $b \in \N$ is much more subtle than the proof for $b=1$. The proof for general
$b \in \N$ is 
given in appendix C in the present paper, where it is the content of part (a) of Theorem \ref{thm:alphabc}. 
Parts (b)--(d) of that theorem explore other properties of $\alphabc$. In particular, in part (b)
we prove that $\alphabc$ is asymptotic to $c$ as $c \goto \infty$. 

We comment on the proof of part (a) of the next theorem for $b \in \N$
because the existence of $\alphabc$ is crucial to the paper. 
Define $\gamma_b(\alpha) = \alpha \zbminusonealpha/\zalphab$, where
$\zbalpha = e^\alpha - \sum_{j=0}^{b-1} \alpha^j/j!$. According to part (a),
if for a given $c \in (b,\infty)$ 
there exists a unique solution $\alpha = \alphabc \in (0,\infty)$ of $\gamma_b(\alpha) = c$,
then it follows that $\rhoalphabc \in \pnbc$. 
The existence of such a solution is a consequence of the following three steps, which are carried out in appendix C: 
$\lim_{\alpha \goto 0^+} \gamma(\alpha) = b$; $\lim_{\alpha \goto \infty} \gamma(\alpha) = \infty$;
$\gamma_b'(\alpha) > 0$ for $\alpha \in (0,\infty)$. To carry out step 3, we note that because $Z_b'(\alpha) =
\zbminusonealpha$, we can write $\gamma_b(\alpha) = (\alpha \log \zalphab)'$ and 
$\gamma_b'(\alpha) = (\alpha \log \zalphab)''$. To prove that $\gamma_b'(\alpha) > 0$, we express
$\zbalpha$ first in terms of an incomplete gamma function and then in terms of a moment generating function. The log-convexity
of the moment generating function and a short calculation involving power series completes the proof.

\begin{thm}
\label{thm:mainestimate}
{\em (a)} Fix a nonnegative integer $b$ and a real number $c \in (b,\infty)$. 
For $\alpha \in (0,\infty)$ let $\rhoalphab$ be the measure on $\Nb$ having components 
\[
\rhoalphabj = \frac{1}{\zalphab} \cdot \frac{\alpha^j}{j!} \ \mbox{ for } j \in \Nb,
\]
where $Z_{0,\alpha} = e^\alpha$, and for $b \in \N$, $\zalphab = e^\alpha - \sum_{j=0}^{b-1}\alpha^j/j!$.
Then there exists a unique value $\alphabc \in (0,\infty)$ such that 
$\rhoalphabc$ lies in the set $\pnbc$ of probability measures on $\Nb$ having mean $c$. If $b=0$, then 
$\alpha_0(c) = c$. If $b \in \N$, then $\alphabc$ is the unique solution in $(0,\infty)$ of $\alpha \zbminusonealpha/\zalphab = c$.

{\em (b)} Fix a nonnegative integer $b$ and a rational number $c \in (b,\infty)$. Let $m$ be the function $m(N)$
appearing in the definitions of $\omeganbm$ in {\em (\ref{eqn:omegankm})} 
and satisfying $m(N) \goto \infty$ and $m(N)^2/N \goto 0$ as $N \goto \infty$.
For any $\nu \in \anbm$ we define $\thetanbnu \in \pnbc$ to have the components $\thetanbnuj = \nuj/N$ for $j \in \Nb$. Then
the relative entropy $R(\thetanbnu | \rhoalphabc)$ is finite, and we have the local large deviation estimate
\[
\frac{1}{N} \log \Pnbm(\Thetanb = \thetanbnu) = -R(\thetanbnu | \rhoalphabc) + \ve_N(\nu).
\]
The quantity $\ve_N(\nu) \goto 0$ uniformly for $\nu \in \anbm$ as $N \goto \infty$.
\end{thm}

\skp

We now prove the local large deviation estimate in part (b) of Theorem \ref{thm:mainestimate}. This proof
is based on a combinatorial argument that is reminiscent of, and as natural as, the combinatorial argument used to prove
Sanov's theorem for empirical measures defined in terms of i.i.d.\ random variables having a finite state space \cite[\S 3]{EllisTaasan2}. 
Part (b) of Theorem \ref{thm:mainestimate} is proved
by analyzing the asymptotic behavior of the product of two multinomial coefficients that we now introduce.

Given $\nu \in \anbm$, our goal is to estimate the probability $\Pnbm(\Thetanb = \thetanbnu)$, where $\thetanbnu$
has the components $\thetanbnuj = \nu_j/N$ for $j \in \Nb$. 
A basic observation is that the set 
$\{\omega \in \omeganbm : \Thetanb(\omega) = \thetanbnu\}$
coincides with the set
\be 
\label{eqn:deltankmnu}
\deltanbmnu = \{\omega \in \omeganbm : \Nj(\omega) = \nu_j \mbox{ for }
j \in \Nb\} .
\ee
It follows that
\bea
\label{eqn:pnkmthetadelta}
\Pnbm(\Thetanb = \thetanbnu) & = & \Pnbm(\deltanbmnu) \\
\nonumber
& = & \frac{1}{\mbox{card}(\omeganbm)} \cdot \mbox{card}(\deltanbmnu).
\eea
Our first task is to determine the asymptotic behavior of $\mbox{card}(\deltanbmnu)$.
In determining the asymptotic behavior of $\mbox{card}(\omeganbm)$, we will use the fact that
$\omeganbm$ can be written as the disjoint union
\be 
\label{eqn:omegankmcup}
\omeganbm = \bigcup_{\nu \in \anbm} \deltanbmnu .
\ee

Let $\nu \in \anbm$ be given. We start by expressing the cardinality of $\mbox{card}(\deltanbmnu)$ as a product of two multinomial coefficients.
For each configuration $\omega \in \deltanbmnu$, $K$ particles are distributed onto the $N$ sites of the lattice $\Lambda_N$ with $j$ particles going
onto $\nu_j$ sites for $j \in \Nb$. We carry this out in two stages. In stage one $K$ particles are placed into $N$ bins, 
$\nu_j$ of which have $j$ particles for $j \in \Nb$. The number of ways of making this placement equals the multinomial coefficient
\be 
\label{eqn:multikminteger}
\frac{K!}{\ds\prod_{j \in \Nb} (j!)^{\nuj}} .
\ee
This multinomial coefficient is well-defined
since $\sum_{j \in \Nb} j\nu_j  = K$. 
Given this placement of $K$ particles into $N$ bins, the number of ways of moving the particles from the bins onto the sites $1,2,\ldots,N$ of the lattice 
$\Lambda_N$ equals the multinomial coefficient
\be 
\frac{N!}{\ds\prod_{j \in \Nb} \nuj!} .
\label{eqn:multinnu}
\ee
This second multinomial coefficient is well-defined since $\sum_{j \in \Nb} \nu_j = N$. 
We conclude that
the cardinality of $\deltanbmnu$ is given by the product of these two multinomial coefficients: 
\be
\label{eqn:cardkminteger}
\mbox{card}(\deltanbmnu) = \ \frac{N!}{\ds\prod_{j \in \Nb}\nuj!} \cdot \frac{K!}{\ds\prod_{j \in \Nb} (j!)^{\nu_j}} .
\ee
Since $|\nu|_+ \leq m$, at most $m$ of the components $\nuj$ are positive. A related version of
this formula, well known in combinatorial analysis, is derived in Example III.23 of \cite{FlaSed}. 

The next two steps in the proof of the local estimate given in part (b) of Theorem \ref{thm:mainestimate} is to prove
the asymptotic formula for $\mbox{card}(\deltanbmnu)$ in Lemma \ref{lem:deltankmnu} and the asymptotic formula for 
$\mbox{card}(\omeganbm)$ in part (b) of Lemma \ref{lem:omegankm}.
The proof of Lemma \ref{lem:deltankmnu} is greatly simplified by a substitution in line 3 of (\ref{eqn:twocoeffagain}).
This substitution involves a parameter $\alpha \in (0,\infty)$, which, we emphasize, 
is arbitrary in this lemma. 
The substitution in line 3 of 
(\ref{eqn:twocoeffagain}) allows us to express the asymptotic behavior of both $\mbox{card}(\deltanbmnu)$ 
in Lemma \ref{lem:deltankmnu} and $\mbox{card}(\omeganbm)$ in Lemma \ref{lem:omegankm} directly in terms of the
relative entropy $R(\thetanbnu|\rhoalphab)$, where $\rhoalphab$ is the probability measure on $\Nb$ having the components defined
in part (a) of Theorem \ref{thm:mainestimate}.
One of the major issues in the proof of part (b) of Theorem \ref{thm:mainestimate} is to show that the arbitrary parameter $\alpha$
appearing in Lemmas \ref{lem:deltankmnu} and \ref{lem:omegankm} must take the value $\alpha_b(c)$, which is the unique
value of $\alpha$ guaranteeing that $\rhoalphab \in \pnbc$ [Thm.\ \ref{thm:mainestimate}(a)]. We show that $\alpha$ must equal
$\alphabc$ after the statement of Lemma \ref{lem:omegankm}. 

\begin{lem}
\label{lem:deltankmnu}
Fix a nonnegative integer $b$ and a rational number $c \in (b,\infty)$. Let $\alpha$ be any real number in $(0,\infty)$, and let $m$ be the function $m(N)$ 
appearing in the definition of $\omeganbm$ in {\em (\ref{eqn:omegankm})} 
and satisfying $m(N) \goto \infty$ and $m(N)^2/N \goto 0$ as $N \goto \infty$.
We define
\[
f(\alpha,b,c,K) = \log\zalphab - c \log \alpha + c\log K - c  .
\]
For any $\nu \in \anbm$, we define $\thetanbnu \in \pnbc$ to have the components $\thetanbnuj = \nuj/N$ for $j \in \Nb$. Then
\begin{eqnarray*}
\lefteqn{
\frac{1}{N} \log \mbox{\em card}(\deltanbmnu)} \\
& & = -R(\thetanbnu | \rhoalphab) + f(\alpha,b,c,K) + \zeta_N(\nu) .
\end{eqnarray*}
The quantity $\zeta_N(\nu) \goto 0$ uniformly for $\nu \in \anbm$ as $N \goto \infty$.
\end{lem}

\noi 
{\bf Proof.} The proof is based on a weak form of Stirling's approximation, which states that 
for all $N \in \N$ satisfying $N \geq 2$ and for all $n \in \N$ satisfying $1 \leq n \leq N$ 
\be 
\label{eqn:stirling}
1 \leq \log (n!) - (n \log n - n) \leq 2 \log N .
\ee
We summarize (\ref{eqn:stirling}) by writing 
\be 
\label{eqn:weakform}
\log (n!) = n \log n - n + \mbox{O}(\log N) \ \forall N \in \N, N \geq 2 \mbox{ and } \forall
n \in \{1,2,\ldots,N\} .
\ee
By (\ref{eqn:stirling}) the term denoted by $\mbox{O}(\log N)$ satisfies $1 \leq \mbox{O}(\log N) \leq 2 \log N$.
We will also use (\ref{eqn:stirling}) with $N$ replaced by $K$ and by other quantities in the model.

To simplify the notation, we rewrite (\ref{eqn:cardkminteger}) in the form
\[
\mbox{card}(\deltanbmnu) = M_1(N,\nu) \cdot M_2(K,\nu) ,
\]
where $M_1(N,\nu)$ denotes the first multinomial coefficient on the right side of (\ref{eqn:cardkminteger}), and $M_2(K,\nu)$ denotes
the second multinomial coefficient on the right side of (\ref{eqn:cardkminteger}). We have
\be 
\label{eqn:twocoeff}
\frac{1}{N} \log \mbox{card}(\deltanbmnu) = \frac{1}{N} \log \mbox{card}(M_1(N,\nu)) + \frac{1}{N} \log \mbox{card}(M_2(K,\nu)) .
\ee

The asymptotic behavior of the first term on the right side of 
the last display is easily calculated. Since $\nu \in \anbm$, there are $|\nu|_+ \in \{1,2,\ldots,m\}$ positive components $\nuj$. 
Because of this restriction on the number $|\nu|_+$ of positive components of $\nu$, we are able to control the error
in line 3 of (\ref{eqn:firststirling}).
We define $\Psi_N(\nu) = \{j \in \Nb : \nuj \geq 1\}$. For each $j \in \Psi_N(\nu)$, since the components $\nuj$ satisfy 
$1 \leq \nuj \leq N$, 
we have 
\[
\log (\nuj!) = \nuj \log \nuj - \nuj + \mbox{O}(\log N) \mbox{ for all } N \geq 2.
\]
Using the fact that $\sum_{j \in \Psi_N(\nu)} \nu_j = N$, we obtain
\begin{eqnarray}
\label{eqn:firststirling}
\lefteqn{
\frac{1}{N} \log \mbox{card}(M_1(N,\nu)) } \\
\nonumber
& & = \frac{1}{N} \log(N!) - \frac{1}{N}\sum_{j \in \Psi_N(\nu)} \log (\nuj!) \\
\nonumber
& & = \frac{1}{N}(N \log N - N + \mbox{O}(\log N)) - \frac{1}{N} \sum_{j \in \Psi_N(\nu)}(\nuj \log \nuj - \nuj + \mbox{O}(\log N)) \\
\nonumber
& & = - \sum_{j \in \Nb}(\nuj/N) \log (\nuj/N) + \frac{\mbox{O}(\log N)}{N} - \frac{1}{N} \sum_{j \in \Psi_N(\nu)}\mbox{O}(\log N) \\
\nonumber 
& & = - \sum_{j \in \Nb}\thetanbnuj \log \thetanbnuj + \zeta^{(1)}_N - \zeta^{(2)}_N(\nu),
\eea
where $\zeta_N^{(1)} = [\mbox{O}(\log N)]/N \goto 0 \mbox{ as } N \goto \infty$ and
\[
\zeta_N^{(2)}(\nu) = \frac{1}{N} \sum_{j \in \Psi_N(\nu)} \mbox{O}(\log N) .
\]
By the inequality noted after (\ref{eqn:weakform}) and the fact that $|\nu|_+ \leq m$
\[
0 \leq \max_{\nu \in \anbm} \zeta_N^{(2)}(\nu) \leq \max_{\nu \in \anbm} \frac{2}{N} \sum_{j \in \Psi_N(\nu)} \log N \leq \frac{2 m \log N}{N} .
\]
Since $(m \log N)/N \goto 0$ as $N \goto \infty$, we conclude that $\zeta_N^{(2)}(\nu) \goto 0$ uniformly for $\nu \in \anbm$ as $N \goto \infty$.

We now study the asymptotic behavior of the second term on the right side of (\ref{eqn:twocoeff}). 
Since $K = Nc$, we obtain for all $K \geq 2$
\begin{eqnarray}
\label{eqn:needstirling}
\lefteqn{
\frac{1}{N} \log \mbox{card}(M_2(K,\nu)) } \\
\nonumber
\nonumber
& & = \frac{1}{N} \log(K!) - \frac{1}{N}\sum_{j \in \Nb}\nuj \log (j!) \\
\nonumber
& & = \frac{1}{N}(K \log K - K + \mbox{O}(\log K)) - \sum_{j \in \Nb}\thetanbnuj \log (j!) \\
\nonumber
& & = c \log K - c - \sum_{j \in \Nb}\thetanbnuj \log (j!) + \zeta_N^{(3)} .
\end{eqnarray}
where 
\[
0 \leq \zeta^{(3)}_N = \frac{\mbox{O}(\log K)}{N} = \frac{\mbox{O}(\log N)}{N} \goto 0 \mbox{ as } N \goto \infty .
\]
The weak form of Stirling's formula is used to rewrite the term $\log(K!)$ in the last display, but
not to rewrite the terms $\log(j!)$, which we leave untouched. 

Substituting (\ref{eqn:firststirling}) and (\ref{eqn:needstirling}) into (\ref{eqn:twocoeff}), we obtain
\bea
\label{eqn:twocoeffmore}
\lefteqn{ \frac{1}{N} \log \mbox{card}(\deltanbmnu)} \\
\nonumber 
&& = \frac{1}{N} \log \mbox{card}(M_1(N,\nu)) + \frac{1}{N} \log \mbox{card}(M_2(K,\nu)) \\
\nonumber
&& = - \sum_{j \in \Nb}\thetanbnuj \log \thetanbnuj - \sum_{j \in \Nb}\thetanbnuj \log (j!) + c \log K - c + \zeta_N(\nu) \\
\nonumber 
&& = - \sum_{j \in \Nb}\thetanbnuj \log(\thetanbnuj j!) + c \log K - c + \zeta_N(\nu) .
\eea
In this formula $\zeta_N(\nu) = \zeta_N^{(1)} - \zeta_N^{(2)}(\nu) + \zeta_N^{(3)}$. As $N \goto \infty$
\[
\max_{\nu \anbm} |\zeta_N(\nu)| \leq \zeta_N^{(1)} + \max_{\nu \in \anbm} \zeta_N^{(2)}(\nu) + \zeta_N^{(3)} \goto 0.
\]
We conclude that $\zeta_N(\nu) \goto 0$ uniformly for $\nu \in \anbm$ as $N \goto \infty$.

Now comes the key step, the purpose of which is to express the sum 
in the last line of (\ref{eqn:twocoeffmore}) as the relative entropy $R(\thetanbnuj | \rhoalphab)$, where $\alpha \in (0,\infty)$ 
is arbitrary. 
To express the sum in the last line of (\ref{eqn:twocoeffmore}) as $R(\thetanbnu | \rhoalphab)$, we rewrite the sum as shown in line 3 of the next display:
\bea
\label{eqn:twocoeffagain} 
\lefteqn{ \frac{1}{N} \log \mbox{card}(\deltanbmnu)} \\
\nonumber 
&& = - \sum_{j \in \Nb}\thetanbnuj \log(\thetanbnuj j!) + c \log K - c + \zeta_N(\nu) \\
\nonumber 
&& = - \sum_{j \in \Nb}\thetanbnuj \log\!\left(\frac{\thetanbnuj}{\alpha^j/(\zalphab\cdot j!)} 
\cdot \frac{\alpha^j}{\zalphab} \right) + c \log K - c + \zeta_N(\nu) \\
\nonumber
&& = - \sum_{j \in \Nb}\thetanbnuj \log(\thetanbnuj / \rhoalphabj) + (\log \zalphab) \sum_{j \in \Nb}\thetanbnuj \\
\nonumber 
&& \hspace{.5in} - (\log \alpha) \sum_{j \in \Nb}j \thetanbnuj  + c \log K - c + \zeta_N(\nu) \\
\nonumber
&& = - R(\thetanbnu | \rhoalphab) + \log \zalphab - c \log \alpha + c \log K - c + \zeta_N(\nu) \\
\nonumber 
&& = - R(\thetanbnu | \rhoalphab) + f(\alpha,b,c,K) + \zeta_N(\nu) .
\eea
We obtain the next-to-last equality by using the fact that since $\thetanbnu \in \pnbc$, 
\[
\sum_{j \in \Nb}\thetanbnuj = 1 \ \mbox{ and } \ \sum_{j \in \Nb}j \thetanbnuj = c.
\] 
The proof of Lemma \ref{lem:deltankmnu} is complete. \ \ink

\skp
The local large deviation estimate in Lemma \ref{lem:deltankmnu} suggests a beautiful connection with Boltzmann's calculation of the
Maxwell--Boltzmann distribution for the random ideal gas.  This connection and
Boltzmann's calculation are described in \cite{EllisTaasan2}. \

The next step in the proof of the local large deviation estimate in part (b) of Theorem \ref{thm:mainestimate} is to prove
the asymptotic formula for $\mbox{card}(\omeganbm)$ stated in part (b) of the next lemma. 
The proof of this lemma uses Lemma \ref{lem:deltankmnu} in a fundamental way.
After the statement of this lemma we show how to apply it and Lemma \ref{lem:deltankmnu} to prove part (b) of Theorem \ref{thm:mainestimate}. 
An important component of this proof
is to calculate the quantity $\min_{\theta \in \pnbc} R(\theta | \rhoalphab)$, which appears in part (b) of the next lemma. 
The proof of part (b) of the lemma depends on part (a), which is also used to verify hypothesis (i) of Theorem \ref{thm:balls} in 
the setting of Theorem \ref{thm:ldlimitballs}.

\begin{lem}
\label{lem:omegankm}
Fix a nonnegative integer $b$ and a rational number $c \in (b,\infty)$. The following conclusions hold. 

{\em (a)} The set $\anbm$ defined at the beginning of section {\em \ref{section:proof1}} has the property that
\[
\lim_{N \goto \infty} \frac{1}{N} \log{\mbox{\em card}}(\anbm) = 0.
\]

{\em (b)} Let $\alpha$ be the positive real number in Lemma {\em \ref{lem:deltankmnu}}, 
and let $m$ be the function $m(N)$ appearing in the definition of $\omeganbm$ in {\em (\ref{eqn:omegankm})}
and satisfying $m(N) \goto \infty$ and $m(N)^2/N \goto 0$ as $N \goto \infty$.
We define
\[
f(\alpha,b,c,K) = \log\zalphab - c \log \alpha + c\log K - c  .
\]
Then $R(\theta | \rhoalphab)$ attains its infimum over $\theta \in \pnbc$, and 
\be 
\label{eqn:frac1n}
\frac{1}{N} \log \mbox{\em card}(\omeganbm) = f(\alpha,b,c,K) - \min_{\theta \in \pnbc}R(\theta | \rhoalphab) + \eta_N .
\ee
The quantity $\eta_N \goto 0$ as $N \goto \infty$.
\end{lem}

Before proving Lemma \ref{lem:omegankm}, we derive the local large deviation estimate in part (b) of
Theorem \ref{thm:mainestimate} by applying Lemmas \ref{lem:deltankmnu} and \ref{lem:omegankm}.  An integral part 
of the proof is to show how the arbitrary value of $\alpha \in (0,\infty)$ appearing in these lemmas is replaced
by the specific value $\alphabc$ appearing in Theorem \ref{thm:mainestimate}. 
As in the statement of part (b) of Theorem \ref{thm:mainestimate}, let $\nu$ be any vector in $\anbm$ and define $\thetanbnu \in \pnbc$ to have the components 
$\thetanbnuj = \nuj/N$ for $j \in \Nb$. 
By (\ref{eqn:pnkmthetadelta})
\bea
\label{eqn:choref}
\lefteqn{
\frac{1}{N} \log \Pnbm(\Thetanb = \thetanbnu) } \\ 
\nonumber & & = \frac{1}{N} \log \Pnbm(\deltanbmnu) \\
\nonumber & & = \frac{1}{N} \log \mbox{card}(\deltanbmnu) - \frac{1}{N} \log \mbox{card}(\omeganbm) .
\eea
Substituting the asymptotic formula for $\log \mbox{card}(\deltanbmnu)$ derived in Lemma
\ref{lem:deltankmnu} and the asymptotic formula for $\log \mbox{card}(\omeganbm)$ given in part (b) of Lemma \ref{lem:omegankm} yields
\bea
\label{eqn:lastone}
\lefteqn{
\frac{1}{N} \log \Pnbm(\Thetanb = \thetanbnu)} \\ 
\nonumber && = -R(\thetanbnu | \rhoalphab) + f(\alpha,b,c,K) + \zeta_N(\nu)  \\
\nonumber && \hspace{.5in} - \left(f(\alpha,b,c,K)- \min_{\theta \in \pnbc} R(\theta | \rhoalphab) + \eta_N\right) \\
\nonumber && = -R(\thetanbnu | \rhoalphab) + \min_{\theta \in \pnbc} R(\theta | \rhoalphab) + \ve_N(\nu) .
\eea
The error term $\ve_N(\nu)$ equals $\zeta_N(\nu) - \eta_N$; $\zeta_N(\nu)$ is the error
term in Lemma \ref{lem:deltankmnu}, and $\eta_N$ is the error term in Lemma \ref{lem:omegankm}.
As $N \goto \infty$, $\zeta_N(\nu) \goto 0$ uniformly for $\nu \in \anbm$, and $\eta_N \goto 0$. It follows that $\ve_N(\nu) \goto 0$
uniformly for $\nu \in \anbm$ as $N \goto \infty$. 

We now consider the first two terms on the right side of the last line of (\ref{eqn:lastone}).  By 
assertion (ii) in part (f) of Theorem \ref{thm:relentropy} applied to $\theta = \thetanbnu \in \pnbc$,
for any $\alpha \in (0,\infty)$
\[
R(\thetanbnu | \rhoalphab) - \min_{\theta \in \pnbc} R(\theta | \rhoalphab) = R(\thetanbnu | \rhoalphabc).
\]
With this step we have succeeded in replacing the relative entropy $R(\thetanbnu | \rhoalphab)$ with respect to $\rhoalphab$, which
appears in Lemma \ref{lem:deltankmnu}, by
the relative entropy $R(\thetanbnu | \rhoalphabc)$ with respect to $\rhoalphabc$, which appears
in Theorem \ref{thm:mainestimate}. Substituting the last equation into (\ref{eqn:lastone})
gives 
\[
\frac{1}{N} \log \Pnbm(\Thetanb = \thetanbnu) = - R(\thetanbnu | \rhoalphabc) + \ve_N(\nu),
\]
where $\ve_N(\nu) \goto 0$ uniformly for $\nu \in \anbm$ as $N \goto \infty$. 
This is the conclusion of part (b) of Theorem \ref{thm:mainestimate}.

We now complete the proof of part (b) of Theorem \ref{thm:mainestimate} by proving Lemma \ref{lem:omegankm}.

\skp
\noi 
{\bf Proof of Lemma \ref{lem:omegankm}.}  (a) To estimate the cardinality of $\anbm$ we write
\[
\anbm \subset \left\{\nu \in \N_0^N : \sum_{j \in \Nb}\nu_j = N, |\nu|_+ \leq m\right\}
= \bigcup_{k=1}^m \left\{\nu \in \N_0^N : \sum_{j \in \Nb}\nu_j = N, |\nu|_+ = k\right\} .
\]
Thus we can bound the cardinality of $\anbm$ by bounding separately the cardinality of each of the disjoint sets in the union. 
By \cite[Cor.\ 2.5]{Charal} the number of elements in the set indexed by $k$ equals the binomial
coefficient $C(N-1,k-1)$. Since by assumption $m/N \goto 0$
as $N \goto \infty$, for all sufficiently large $N$ the quantities $C(N-1,k-1)$ are increasing and are maximal when $k = m$. 
Since $C(N-1,k-1) \leq C(N,k)$, it follows that
\[
\mbox{card}(\anbm) \leq \sum_{k=1}^m C(N,k) \leq m C(N,m) = m \frac{N!}{m! (N-m)!} .
\]
An application of the weak form of Stirling's formula yields for all $m \geq 2$ and all $N \geq m+2$
\beas
0 & \leq & \frac{1}{N} \log \mbox{card}(\anbm) \\
& \leq & \frac{1}{N} (\log m + \log(N!) - \log(m!) - \log((N-m)!))) \\
& = & \frac{\log m}{N} - \frac{m}{N}\log \frac{m}{N} - \left(1-\frac{m}{N}\right) \log \left(1-\frac{m}{N}\right) + 
\frac{\mbox{O}(\log N)}{N} .
\eeas
Since $m/N \goto 0$ as $N \goto \infty$, we conclude that as $N \goto \infty$
\beas
0 & \leq & \frac{1}{N} \log \mbox{card}(\anbm) \\
\nonumber
& \leq & \frac{\log m}{N} - \frac{m}{N}\log \frac{m}{N} - 
\left(1-\frac{m}{N}\right) \log \left(1-\frac{m}{N}\right) + \frac{\mbox{O}(\log N)}{N}
\goto 0 .
\eeas
This completes the proof of part (a). 

\skp
(b) The starting point is 
(\ref{eqn:omegankmcup}), which states that
\[
\omeganbm = \bigcup_{\nu \in \anbm} \deltanbmnu .
\]
For distinct $\nu \in \anbm$ the sets $\deltanbmnu$ are disjoint. Hence
\bea
\label{eqn:lineone}
\lefteqn{
\frac{1}{N} \log \mbox{card}(\omeganbm)} \\ 
\nonumber & & = \frac{1}{N} \log \sum_{\nu\in\anbm} \mbox{card}(\deltanbmnu) \\
\nonumber & & = \frac{1}{N} \log \left(\max_{\nu\in\anbm} \mbox{card}(\deltanbmnu) \cdot \sum_{\nu\in\anbm} 
\frac{\mbox{card}(\deltanbmnu)}{\max_{\nu\in\anbm} \mbox{card}(\deltanbmnu)}\right) \\
\nonumber 
& & = \frac{1}{N} \log \!\left(\max_{\nu\in\anbm} \mbox{card}(\deltanbmnu)\right) + \delta_N ,
\eea
where 
\[
0 < \delta_N = \frac{1}{N} \log \left(\sum_{\nu \in \anbm} 
\frac{\mbox{card}(\deltanbmnu)}{\max_{\nu\in\anbm} \mbox{card}(\deltanbmnu)}\right) \leq \frac{1}{N} \log \mbox{card}(A_{K,N,m}) .
\]
It follows from part (a) that $\delta_N \goto 0$ as $N \goto \infty$.

We continue with the estimation of $\mbox{card}(\omeganbm)$. By Lemma \ref{lem:deltankmnu} and the fact that logarithm is an increasing function
\beas
\lefteqn{
-\min_{\nu\in\anbm} R(\thetanbnu | \rhoalphab) + f(\alpha,b,c,K) - \max_{\nu\in \anbm} |\zeta_N(\nu)|} \\
&& \leq \max_{\nu\in\anbm} \left(\frac{1}{N} \log \mbox{card}(\deltanbmnu)\right) \\
&& = \frac{1}{N} \log \left(\max_{\nu\in\anbm} \mbox{card}(\deltanbmnu)\right) \\
&& \leq -\min_{\nu\in\anbm} R(\thetanbnu | \rhoalphab) + f(\alpha,b,c,K) + \max_{\nu\in \anbm} |\zeta_N(\nu)|.
\eeas
As proved in Lemma \ref{lem:deltankmnu}, $\max_{\nu\in \anbm} |\zeta_N(\nu)| \goto 0$ as $N \goto \infty$. 
Hence by (\ref{eqn:lineone})
\bea
\label{eqn:lefteqntwo}
\lefteqn{
-\min_{\nu\in\anbm} R(\thetanbnu | \rhoalphab) + f(\alpha,b,c,K) - \max_{\nu\in \anbm} |\zeta_N(\nu)| + \delta_N } \\
\nonumber 
&& \leq \frac{1}{N} \log \mbox{card}(\omeganbm) \\
\nonumber
&& \leq -\min_{\nu\in\anbm} R(\thetanbnu | \rhoalphab) + f(\alpha,b,c,K) + \max_{\nu\in \anbm} |\zeta_N(\nu)| 
+ \delta_N .
\eea

Under the assumption that $R(\cdot | \rhoalphab)$ attains its infimum over $\pnbc$, we define
\[
\eta_N = \frac{1}{N} \log \mbox{card}(\omeganbm) - f(\alpha,b,c,K) + \min_{\theta \in \pnbc} R(\theta | \rhoalphab).
\]
In the last two paragraphs of this proof, we show that $\eta_N \goto 0$ as $N \goto \infty$. Given this fact, 
the last equation
yields the asymptotic formula (\ref{eqn:frac1n}) in part (b).

We now prove that $\eta_N \goto 0$ as $N \goto \infty$. To do this, we use (\ref{eqn:lefteqntwo}) to write
\[ 
|\eta_N| \leq \left(\min_{\nu\in\anbm} R(\thetanbnu | \rhoalphab) - \min_{\theta \in \pnbc} R(\theta | \rhoalphab)\right)
+ \max_{\nu\in \anbm} |\zeta_N(\nu)| + \delta_N .
\]
Like the second and third terms on the right side, the first term on the right side is nonnnegative because $\anbm$ is a subset of $\pnbc$.
Since $\max_{\nu\in \anbm} |\zeta_N(\nu)| \goto 0$ and $\delta_N \goto 0$ as $N \goto \infty$, 
it will follow that $\eta_N \goto 0$ if we can show that $R(\cdot | \rhoalphab)$ attains its infimum over $\pnbc$ and that 
\be 
\label{eqn:complicated}
\lim_{N \goto \infty} \min_{\nu \in \anbm} R(\thetanbnu | \rhoalphab) = \min_{\theta \in \pnbc} R(\theta | \rhoalphab).
\ee
Given the existence of $\min_{\theta \in \pnbc} R(\theta | \rhoalphab)$, 
this assertion is certainly plausible since as shown in Corollary \ref{cor:dense}, the measures $\thetanbnu$ are 
dense in $\pnbc$ for $\nu \in \cup_{N \in \N}\anbm$.

We start the proof of (\ref{eqn:complicated}) by noting that since $R(\cdot | \rhoalphab)$ has compact level sets in $\pnbc$ [Thm.\ \ref{thm:relentropy}(d)],
$R(\cdot | \rhoalphab)$ attains its infimum over $\pnbc$ at some measure $\theta^\star$. In assertion (i) in part (f) of Theorem \ref{thm:relentropy}, we show
that $\theta^\star = \rhoalphabc$. However, this detail is not needed in the present proof, which we would like
to keep as self-contained as possible. We prove (\ref{eqn:complicated}) by applying Theorem \ref{thm:approximate} to $\theta = \theta^\star$,
obtaining a sequence $\thetan$ with the following properties: 
\begin{itemize}
  \item For $N \in \N$, $\thetan \in \bnbm$ has components $\thetanj = \nunj/N$ for $j \in \Nb$, where $\nun$ is an appropriate sequence in $\anbm$.
  \item $\thetan \Rightarrow \theta^\star$ as $N \goto \infty$.
  \item $R(\thetan | \rhoalphab) \goto R(\theta^\star | \rhoalphab)$ as $N \goto \infty$.
\end{itemize}
The limit in (\ref{eqn:complicated}) follows from the inequalities
\[
\min_{\theta \in \pnbc} R(\theta | \rhoalphab) \leq \min_{\nu \in \anbm} R(\thetanbnu | \rhoalphab) \leq R(\thetan | \rhoalphab)
\]
and the limit  \\ 
\[
R(\thetan | \rhoalphab) \goto R(\theta^\star | \rhoalphab) = \min_{\theta \in \pnbc} R(\theta | \rhoalphab) \mbox{ as } N \goto \infty.
\]
This completes the proof of Lemma \ref{lem:omegankm} and 
thus the proof of the local estimate in part (b) of Theorem \ref{thm:mainestimate}. \ \ink

\skp

We end this section by explaining the insight behind the key step in the proof of Lemma \ref{lem:deltankmnu}. This key step is to rewrite 
the sum in line 2 of (\ref{eqn:twocoeffagain}) as shown in line 3. This allows us to express
the sum in line 3 as the relative entropy $R(\thetanbnu | \rhoalphabc)$ plus terms that are independent of $\thetanbnu$. 
We now motivate this step. In order to streamline this motivation, we drop all error terms and avoid rigor.

Our starting point is line 2 of (\ref{eqn:twocoeffagain}). 
If we do not rewrite the sum as shown in line 3 of that display, then we have the following
modification of the conclusion of Lemma \ref{lem:deltankmnu}:
\be 
\label{eqn:star}
\frac{1}{N} \log \mbox{card}(\deltanbmnu) \approx - \sum_{j \in \Nb}\thetanbnuj \log(\thetanbnuj j!) + c \log K - c.
\ee
This in turn leads to the following modification of Lemma \ref{lem:omegankm}:
\[
\frac{1}{N} \log \mbox{ card }(\omeganbm) \approx c \log K - c - \min_{\nu \in \anbm}\left(\sum_{j \in \Nb}\thetanbnuj \log(\thetanbnuj j!)\right).
\]
For $\nu \in \cup_{N \in \N}\anbm$ the probability measures $\thetanbnu$ are dense in $\pnc$ [Cor.\ \ref{cor:dense}]. 
Hence it is plausible that as $N \goto \infty$ the minimum in the last display can be replaced by
\be 
\label{eqn:needthisone}
\min_{\theta\in\pnbc}\left(\sum_{j \in \Nb} \thetaj \log(\thetaj j!)\right).
\ee
To determine this minimum, we introduce two Lagrange multipliers corresponding to the two equality constraints $\sum_{j\in\Nb} \thetaj = 1$
and $\sum_{j \in \Nb} j\thetaj = c$ satisfied by $\theta \in \pnbc$. A formal calculation, which we omit, suggests that the minimum is attained at the unique
$\theta \in \pnbc$ having components 
\[
\thetaj = \frac{1}{Z_b(\alpha)} \cdot \frac{\alpha^j}{j!} \mbox{ for } j \in \Nb,
\]
where $\alpha = \alphabc$ and $Z_b(\alpha) = Z_b(\alphabc)$ are chosen so that 
$\sum_{j \in \Nb} \thetaj = 1$ and $\sum_{j\in\Nb} j \thetaj = c$ [Thm. \ref{thm:mainestimate}(a)]. The measure $\theta$ with $\alpha = \alphabc$ 
coincides with the Poisson distribution $\rhoalphabc$ 
appearing in the local large deviation estimate in part (b) of Theorem \ref{thm:mainestimate}. One easily checks that the value of the minimum 
in (\ref{eqn:needthisone}) is $c\log\alphabc - \log\zalphabc$. These calculations suggest that 
\be 
\label{eqn:starstar}
\frac{1}{N} \log \mbox{card}(\omegankm) \approx   c \log K - c - c\log\alphabc + \log\zalphabc.
\ee

When (\ref{eqn:starstar}) is combined with (\ref{eqn:star}), we have by (\ref{eqn:choref}) 
\beas
\lefteqn{
\frac{1}{N} \log \pnkm(\Thetanb = \thetanbnu) } \\ 
& & = \frac{1}{N} \log \pnkm(\deltanbmnu) \\
& & = \frac{1}{N} \log \mbox{card}(\deltanbmnu) - \frac{1}{N} \log \mbox{card}(\omeganbm) \\
& & \approx  -\sum_{j \in \Nb}\thetanbnuj \log (\thetanbnuj j!) + c \log K - c \\
& & \hspace{.5in} - (c \log K - c - c\log\alphabc + \log\zalphabc \\
& & \approx -\sum_{j \in \Nb}\thetanbnuj \log (\thetanbnuj j!) + c\log\alphabc - \log\zalphabc.
\eeas
The last line of this display can be rewritten as
\beas
\lefteqn{
-\sum_{j \in \Nb}\thetanbnuj \log\left(\frac{\thetanbnuj}{[\alphabc]^j/(\zalphabc \cdot j!)}\right)} \\
&& = -\sum_{j \in \Nb}\thetanbnuj \log(\thetaj/\rhoalphabcj) = - R(\thetanbnuj | \rhoalphabc).
\eeas
It follows that
\[
\frac{1}{N} \log \pnkm(\Thetanb = \thetanbnu) \approx -R(\thetanbnuj|\rhoalphabc).
\]
Except for the error terms, this coincides with the conclusion of part (b) of Theorem \ref{thm:mainestimate}.

The calculation just presented was our first attempt to prove Lemmas \ref{lem:deltankmnu} and \ref{lem:omegankm}.
It also guided us to the much more efficient current proofs both of Lemma \ref{lem:deltankmnu} --- where the sum
in line 2 of (\ref{eqn:twocoeffagain}) is written directly in terms of the relative entropy --- and of Lemma \ref{lem:omegankm}. 
An analogous but much simpler calculation motivates the solution of a finite dimensional problem involving the minimum of a relative
entropy over a set of probability measures having fixed mean. This simpler calculation is directly related to the present paper
because it gives the form of the Maxwell--Boltzmann distribution for a random ideal gas. 
For details see section 6.4 of \cite{LesHouches}, sections 4-5 of \cite{Ellis-Boltzmann}, and section 4 of \cite{EllisTaasan2},
each of which emphasizes different aspects of the calculation. This completes
the motivation of the proof of Lemma \ref{lem:deltankmnu}. 

In the next section we show how the local large deviation estimate in part (b) of Theorem
\ref{thm:mainestimate} yields the LDP in Theorem \ref{thm:ldpthetankm}.

\section{Proof of Theorem \ref{thm:ldpthetankm} from Part (b) of Theorem \ref{thm:mainestimate}}
\label{section:proof2}
\beginsec

In Theorem \ref{thm:ldpthetankm} we state the
LDP for the sequence $\Thetanb$ of number-density measures. This sequence takes values
in $\pnbc$, which is the set of probability measures on $\N$ having mean $c \in (b,\infty)$. 
The purpose of the present section is to 
show how the local large deviation estimate in part (b) of Theorem \ref{thm:mainestimate}
yields the LDP for $\Thetanb$. The basic idea is first to prove the large deviation limit for $\thetanbnu$ lying in open balls in $\pnbc$ and in other subsets defined in terms of open balls
and then to use this large deviation limit
to prove the LDP in Theorem \ref{thm:ldpthetankm}. Both of these steps are implemented as applications
of the general formulation in Theorems \ref{thm:balls} and \ref{thm:ballstoldp}.

In Theorem \ref{thm:ldlimitballs} we state the large deviation limit for open balls and other subsets defined
in terms of open balls. 
Two types of open balls are considered. Let $\theta$ be a measure in $\pnbc$, and take $r > 0$. 
Part (a) states the large deviation limit for open balls in $\pnbc$ defined by
\[
B_\pi(\theta, r) = \{\mu \in \pnbc: \pi(\theta,\mu) < r\} ,
\]
where $\pi$ denotes the Prohorov metric on $\pnbc$ \cite[\S 3.1]{EthierKurtz}. This limit will be used to prove the large deviation
upper bound for compact subsets of $\pnbc$ in part (b) of Theorem \ref{thm:ldpthetankm} and the large deviation lower bound for open subsets of $\pnbc$
in part (d) of Theorem \ref{thm:ldpthetankm}. Now let $\theta$ be a measure in $\pnbbc$. Part (b) states the large deviation limit for
sets of the form $\hatb_\pi(\theta,r) \cap \pnbc$, where $\hatb_\pi(\theta,r)$ is the open ball in $\pnbbc$ defined by
\[
\hatb_\pi(\theta,r) = \{\mu \in \pnbbc: \pi(\theta,\mu) < r\}.
\]
This limit will be used to prove the large deviation upper bound for closed subsets in part (c) of Theorem \ref{thm:ldpthetankm}. 
Since $\pnbc$ is a dense subset of $\pnbbc$ [Thm.\ \ref{thm:pnc}(b)], 
$\hatb_\pi(\theta,r) \cap \pnbc$ is nonempty. If $\theta \in \pnbc$, then $B_\pi(\theta, r) = \widehat{B}_\pi(\theta, r) \cap \pnbc$,
and the conclusions of parts (a) and (b) of the next theorem coincide.
For $A$ a subset of $\pnbc$ or $\pnbbc$ we denote by $R(A | \rhoalphabc)$ the infimum of $R(\theta | \rhoalphabc)$ over $\theta \in A$.

\begin{thm}
\label{thm:ldlimitballs}
Fix a nonnegative integer $b$ and 
a rational number $c \in (b,\infty)$. Let $m$ be the function $m(N)$ appearing in the definitions of $\omeganbm$ in {\em (\ref{eqn:omegankm})}
and satisfying $m(N) \goto \infty$ and $m(N)^2/N \goto 0$ as $N \goto \infty$.
The following conclusions hold.

{\em (a)} Let $\theta$ be a measure in $\pnbc$ and take $r > 0$. Then for any open ball $B_\pi(\theta,r)$ in $\pnbc$,
$R(B_\pi(\theta,r) | \rhoalphabc)$ is finite, and we have the large deviation limit
\[
\lim_{N \goto \infty} \frac{1}{N} \log \Pnbm(\Thetanb \in B_\pi(\theta,r)) = -  R(B_\pi(\theta,r) | \rhoalphabc) .
\]

{\em (b)} Let $\theta$ be a measure in $\pnbbc$ and take $r > 0$. Then the set $\hatb_\pi(\theta,r) \cap \pnbc$ is nonempty,
$R(\hatb_\pi(\theta,r) \cap \pnbc | \rhoalphabc)$ is finite, and we have the large deviation limit
\[
\lim_{N \goto \infty} \frac{1}{N} \log \Pnbm(\Thetanb \in \hatb_\pi(\theta,r) \cap \pnbc) = - R(\hatb_\pi(\theta,r) \cap \pnbc | \rhoalphabc) .
\]
\end{thm}

We prove Theorem \ref{thm:ldlimitballs} by applying the local large deviation estimate in Lemma \ref{lem:deltankmnu}. 
A key step is to approximate probability measures in $B_\pi(\theta,\ve)$ and in $\hatb_\pi(\theta,r) \cap \pnbc$
by appropriate sequences of 
probability measures in the range of $\Thetanb$. This procedure allows one to show in part (a) that 
the infimum $R(B_\pi(\theta, \ve)| \rhoalphabc)$ can be approximated
by the infimum of $R(\theta | \rhoalphabc)$ over $\theta$ lying in the intersection of $B_\pi(\theta,\ve)$ and the range of $\Thetanb$; a similar
statement holds for the infimum in part (b). A set of 
hypotheses that allow one to carry out this approximation procedure is given in Theorem \ref{thm:balls}, a general formulation
that yields Theorem \ref{thm:ldlimitballs} as a special case. 

Theorem \ref{thm:balls} is formulated for a complete, separable metric space $\X$ containing a relatively compact subset $\W$ that is not closed. We define
$\mathZ$ to be the closure of $\W$ in $\X$. In the application to Theorem \ref{thm:ldlimitballs} $\X$ equals $\pnb$, the set
of probability measures on $\N$; $\W$ equals $\pnbc$, the subset of $\pnb$ containing probability measures with mean $c$; 
and $\mathZ$ equals $\pnbbc$, the subset of $\pnb$ containing probability measures with mean lying in the closed interval $[b,c]$.
If $\tau$ denotes the metric on $\X$, then for $x \in \W$ and $r > 0$ open balls in $\W$ have the form
\[
B_\tau(x,r) = \{y \in \W : \tau(x,y) < r\}.
\]
For $x \in \mathZ$ and $r > 0$ open balls in $\mathZ$ have the form
\[
\widehat{B}_\tau(x,r) = \{y \in \mathZ : \tau(x,y) < r\} .
\]

\begin{thm} 
\label{thm:balls}
For $N \in \N$ let $(\Omega_N,\mathcal{F}_N,Q_N)$ be a 
sequence of probability spaces. Let $\X$ be a complete, separable metric space, $\W$ a relatively compact subset of $\X$
that is not closed and thus not compact, and $\mathZ$
the closure of $\W$ in $\X$; thus $\mathZ$ is compact. Also let $Y_N$ be a sequence of random vectors mapping 
$\Omega_N$ into $\W$, and
let $I$ be a function mapping $\X$ into $[0,\infty]$. For $A$ a subset of $\X$ we denote the infimum of $I$ over $A$ by $I(A)$. 
We assume the following four hypotheses. 

{\em (i)} For $\omega \in \Omega$ the range of $Y_N(\omega)$ is a finite subset $\W_N$ of $\W$, and the cardinality of $\W_N$ satisfies
\[
\lim_{N \goto \infty} \frac{1}{N} \log \mbox{\em card}(\W_N) = 0 .
\]

{\em (ii)} For each $y \in \W_N$ we have $I(y) < \infty$ and the local large deviation estimate
\[
\frac{1}{N} \log Q_N(Y_N = y) = -I(y) + \ve_N(y) , 
\]
where $\ve_N(y) \goto 0$ as $N \goto \infty$ uniformly for $y \in \W_N$.

{\em (iii)} There exists a dense subset $\D$ of $\W$ such that $I(y) < \infty$ for all $y \in \D$.

{\em (iv)} For any $y \in \W$ satisfying $I(y) < \infty$, there exists a sequence $y_N \in \W_N$
for which $y_N \goto y$ and $I(y_N) \goto I(y)$ as $N \goto \infty$. 

Under these hypotheses the following conclusions hold.

{\em (a)} For any open ball $B$ in $\W$, $I(B)$ is finite, and we have the large deviation limit
\[
\lim_{N \goto \infty} \frac{1}{N} \log Q_N(Y_N \in B) = -I(B) .
\]

{\em (b)} For any open ball $\hatb$ in $\mathZ$, $\hatb \cap \W$ is nonempty, $I(\hatb \cap \W)$ is finite, and we have the large deviation limit
\[
\lim_{N \goto \infty} \frac{1}{N} \log Q_N(Y_N \in \hatb \cap \W) = -I(\hatb \cap \W) .
\]
\end{thm}

\noindent
{\bf Proof.}  (a) By hypothesis (iii), for any open ball $B$ in $\W$ there exists $x \in B \cap \D$ such that $I(x) < \infty$. 
Thus $I(B) \leq I(x) < \infty$. By the local large deviation estimate in hypothesis (ii)
\[
Q_N(Y_N \in B) = \sum_{y \in B \cap \W_N} Q_N(Y_N = y)
= \sum_{y \in B \cap \W_N} \exp[-N( I(y) - \ve_N(y))] .
\]
For the last sum in this equation we have the bounds
\beas
\max_{y \in B \cap W_N} \exp[-N( I(y) - \ve_N(y))] & \leq & \sum_{y \in B \cap \W_N} \exp[-N( I(y) - \ve_N(y))]  \\
& \leq & \mbox{card}(\W_N) \cdot \max_{y \in B \cap \W_N} \exp[-N( I(y) - \ve_N(y))].
\eeas
In addition, for the term $\max_{y \in B \cap \W_N} \exp[-N( I(y) - \ve_N(y))]$ we have the bounds
\beas
\lefteqn{
\exp\!\left[-N \left(I(B \cap W_N) + \max_{y \in B \cap \W_N}\ve_N(y)\right)\right]} \\
& = & \exp\!\left[-N \left(\min_{y \in B \cap \W_N} I(y) + \max_{y \in B \cap \W_N}\ve_N(y)\right)\right] \\
& \leq & \max_{y \in B \cap \W_N} \exp[-N( I(y) - \ve_N(y))] \\
& \leq & \exp\!\left[-N \left(\min_{y \in B \cap \W_N} I(y) - \max_{y \in B \cap \W_N}\ve_N(y)\right)\right] \\
& = & \exp\!\left[-N \left(I(B \cap W_N) - \max_{y \in B \cap \W_N}\ve_N(y)\right)\right].
\eeas
It follows that 
\beas
\lefteqn{
- I(B \cap \W_N) - \max_{y \in B \cap \W_N}\ve_N(y)} \\
& \leq & \frac{1}{N} \log Q_N(Y_N \in B) \\
& \leq & - I(B \cap \W_N) + \max_{y \in B \cap \W_N}\ve_N(y) + \frac{\log(\mbox{card}(\W_N))}{N} .
\eeas

Since $\ve_N(y) \goto 0$ uniformly for $y \in \W_N$, by hypothesis (i) the proof is done once we show that 
\be 
\label{eqn:limk}
\lim_{N \goto \infty} I(B \cap \W_N) = I(B) .
\ee
Since $B \cap \W_N \subset B$, we have $I(B) \leq I(B \cap \W_N)$, which implies that
\[
I(B) \leq \liminf_{N \goto \infty} I(B \cap \W_N).
\]
The limit in (\ref{eqn:limk}) is proved if we can show that 
\be 
\label{eqn:limsup}
\limsup_{N \goto \infty} I(B \cap \W_N) \leq I(B).
\ee
For any $\delta > 0$ there exists $y^\star \in B$ such that $I(y^\star) \leq I(B) + \delta < \infty$. 
Hypothesis (iv) guarantees the existence of a sequence $y_N \in \W_N$
such that $y_N \goto y^\star$ and $I(y_N) \goto I(y^\star)$. Since for all sufficiently large $N$ we have $y_N \in B \cap \W_N$, it follows that
$I(B \cap \W_N) \leq I(y_N)$. Hence 
\[
\limsup_{N \goto \infty} I(B \cap \W_N) \leq \lim_{N \goto \infty} I(y_N) = I(y^\star) \leq I(B) + \delta .
\]
Taking $\delta \goto 0$ gives (\ref{eqn:limsup}) and thus proves the limit (\ref{eqn:limk}). This completes the proof of part (a). 

(b) Let $\hatb$ be any open ball in $\mathZ$. Since $\W$ is dense in $\mathZ$, $\hatb \cap \W$ is nonempty. 
By hypothesis (iii) there exists $x \in \hatb \cap \D$ such that $I(x) < \infty$. Thus $I(\hatb \cap \W) \leq I(\hatb \cap \D) \leq I(x) < \infty$. 
To prove the limit in part (b), we proceed as in the proof of the limit in part (a), replacing the set $B$ in part (a) by the set $\hatb \cap \W$.
Since $\W_N \subset \W$, we have $\hatb \cap \W \cap \W_n = \hatb \cap \W_N$. By the local large deviation estimate in hypothesis (ii)
\beas
Q_N(Y_N \in \hatb \cap \W_N) & = & \sum_{y \in \hatb \cap \W \cap \W_N} Q_N(Y_N = y) \\
& = & \sum_{y \in \hatb \cap \W_N} Q_N(Y_N = y) = \sum_{y \in \hatb \cap \W_N} \exp[-N( I(y) - \ve_N(y))].
\eeas
Exactly as in the proof of part (a), it follows that
\beas
\lefteqn{
- I(\hatb \cap \W_N) - \max_{y \in \hatb \cap \W_N}\ve_N(y)} \\
& \leq & \frac{1}{N} \log Q_N(Y_N \in \hatb \cap \W_N) \\
& \leq & - I(\hatb \cap \W_N) + \max_{y \in \hatb \cap \W_N}\ve_N(y) + \frac{\log(\mbox{card}(\W_N))}{N} .
\eeas

Since $\ve_N(y) \goto 0$ uniformly for $y \in W_N$, by hypothesis (i) the proof is done once we show that 
\be 
\label{eqn:limkk}
\lim_{N \goto \infty} I(\hatb \cap \W_N) = I(\hatb \cap \W) .
\ee
Since $\hatb \cap \W_N \subset \hatb \cap \W$, we have $I(\hatb \cap \W) \leq I(\hatb \cap \W_N)$, which implies that
\[
I(\hatb \cap \W) \leq \liminf_{N \goto \infty} I(\hatb \cap \W_N).
\]
The limit in (\ref{eqn:limk}) is proved if we can show that 
\be 
\label{eqn:limsupp}
\limsup_{N \goto \infty} I(\hatb \cap \W_N) \leq I(\hatb \cap \W).
\ee
For any $\delta > 0$ there exists $y^\star \in \hatb \cap \W$ such that 
$I(y^\star) \leq I(\hatb \cap \W) + \delta < \infty$. 
Hypothesis (iv) guarantees the existence of a sequence $y_N \in \W_N$
such that $y_N \goto y^\star$ and $I(y_N) \goto I(y^\star)$. Since for all sufficiently large $N$ 
we have $y_N \in \hatb \cap \W_N$, it follows that
$I(\hatb \cap \W_N) \leq I(y_N)$. Hence 
\[
\limsup_{N \goto \infty} I(\hatb \cap \W_N) \leq \lim_{N \goto \infty} I(y_N) = I(y^\star) \leq I(\hatb \cap \W) + \delta .
\]
Taking $\delta \goto 0$ gives (\ref{eqn:limsupp}) and thus proves the limit (\ref{eqn:limkk}). This completes the proof of part (b) and thus the proof of the theorem. \ \ink

\skp
We now prove Theorem \ref{thm:ldlimitballs} as an application of Theorem \ref{thm:balls}. 
 In Theorem \ref{thm:balls} we make the following 
identifications for $N \in \N$.
\begin{itemize}
  \item The probability
spaces $(\Omega_N,\mathcal{F}_N,Q_N)$ are $(\omeganbm,\mathcal{F}_{N,b,m}, \Pnbm)$, where $\omeganbm$ is the set defined
in (\ref{eqn:omegankm}), $\mathcal{F}_{N,b,m}$ is the $\sigma$-algebra of all subsets of $\omeganbm$, and $\Pnbm$ is the conditional probability
defined in (\ref{eqn:condprob}). 
    \item $\X$ equals $\pnb$, $\W$ equals $\pnbc$, and $\mathZ$ equals $\pnbbc$. These spaces have the properties postulated in Theorem \ref{thm:balls}:  
$\pnb$ is a complete, separable metric space; $\pnbc$ is relatively compact subset of $\pnb$ that is not
closed; and $\pnbbc$ is the closure of $\pnbc$ in $\pnb$. The properties of $\pnb$ are proved in Theorems 3.3.1 and Theorem 3.1.7
of \cite{EthierKurtz}, and the properties
of $\pnbc$ and $\pnbbc$ are proved in Theorem \ref{thm:pnc}. 
    \item The random vectors $Y_N$ equal $\Thetanb$, where $\Thetanb$ is the number-density measure defined in (\ref{eqn:thetankm}). 
$\Thetanb$ maps $\omeganbm$ into the subspace $\W = \pnbc$ of $\pnb$. 
  \item The function $I$ is the relative entropy $R(\cdot | \rhoalphabc)$ on $\pnb$.
$R(\cdot | \rhoalphabc)$ maps $\pnb$ into $[0,\infty]$ [Thm.\ \ref{thm:relentropy}(a)], as specified in the third sentence of Theorem \ref{thm:balls}.
  \item The range $\W_N$ of $Y_N = \Thetanb$ is the set of probability measures $\thetanbnu \in \bnbm$,
the components of which are specified in (\ref{eqn:thetankmj}). The set $\bnbm \subset \pnbc$ is defined in (\ref{eqn:bnkm}).
\end{itemize}

We now verify that the four hypotheses of Theorem \ref{thm:balls} are valid in the setting of Theorem \ref{thm:ldlimitballs}.

\skp
\noi
{\it Verification of hypothesis {\em (i)} in Theorem {\em \ref{thm:balls}}.} In the setting of Theorem \ref{thm:ldlimitballs}
$\W_N$ is range of $\Thetanb(\omega)$ for $\omega \in \omeganbm$. This range is $\bnbm$, the elements of
which are in one-to-one correspondence with the elements of the set $\anbm$
defined in (\ref{eqn:nuj3}). 
As shown in part (a) of Lemma \ref{lem:omegankm}
\[
0 \leq \frac{\log \mbox{card}(\W_N)}{N} = \frac{\log \mbox{card}(\anbm)}{N} 
\goto 0 \ \mbox{ as } N \goto \infty .
\]
This completes the verification of hypothesis (i) in Theorem \ref{thm:balls}.

\skp
\noindent
{\it Verification of hypothesis {\em (ii)} in Theorem {\em \ref{thm:balls}}.} In the setting of Theorem \ref{thm:ldlimitballs} 
hypothesis (ii) in Theorem \ref{thm:balls} is
given by the local estimate in part (b) of Theorem \ref{thm:mainestimate}. 
As shown there, the error $\ve_N(\nu) \goto 0$ as $N \goto \infty$ 
uniformly for $\nu \in \anbm$.
Since there is a one-to-one correspondence between $\nu \in \anbm$ and $\theta\in \bnbm$,
the error in part (b) of Theorem \ref{thm:mainestimate} converges to 0 uniformly for $\theta \in \bnbm$, 
which is the range of $\Thetanb(\omega)$ for $\omega \in \omeganbm$. 
This completes the verification of hypothesis (ii) in Theorem \ref{thm:balls}.

\skp
\noi 
{\it Verification of hypothesis {\em (iii)} in Theorem {\em \ref{thm:balls}}.}  The fact that there exists a dense subset 
of $\theta \in \pnbc$ for which $R(\theta | \rho) < \infty$ is proved in Corollary \ref{cor:dense}.
This completes the verification of hypothesis (iii) in Theorem \ref{thm:balls}.

\skp
\noi 
{\it Verification of hypothesis {\em (iv)} in Theorem {\em \ref{thm:balls}}.} In Theorem \ref{thm:approximate} we prove that any $\alpha \in (0,\infty)$ and any $\theta \in \pnbc$
satisfying $R(\theta | \rhoalphab) < \infty$
there exists a sequence $\thetan \in \bnbm$ for which $\thetan \Rightarrow \theta$ and $R(\thetan | \rhoalphab) \goto R(\theta | \rhoalphab)$
as $N \goto \infty$. In particular, this property holds for $\alpha = \alphabc$. This completes the verification of hypothesis (iv) in Theorem \ref{thm:balls}. 

\skp
Having verified the four hypotheses of Theorem \ref{thm:balls} in the context of Theorem \ref{thm:ldlimitballs}, we have finished the proof of the latter theorem
from the former theorem.

\skp
Theorem \ref{thm:ldpthetankm} states the LDP for the 
number-density measures $\Thetanb$ in the droplet model.
In order to complete the proof of Theorem \ref{thm:ldpthetankm}, we show how to lift the large
deviation limit for open balls in Theorem \ref{thm:ldlimitballs} to the large deviation upper bound for 
compact sets and for closed sets in $\pnbc$
and the large deviation lower bound for open sets in $\pnbc$. This procedure is carried out as an application of
Theorem \ref{thm:ballstoldp}, a general result formulated in a setting close to that of Theorem \ref{thm:balls}. 
In Theorem \ref{thm:ballstoldp} the assumption in Theorem \ref{thm:balls} on the function $I$ is strengthened
to the assumption that $I$ is lower semicontinuous on $\X$.

The LDP in the next theorem has a number of unique features
because $\W$ is not a closed subset of $\X$. 
The large deviation upper bound takes two forms depending on whether the subset $F$ of $\W$ is compact or 
whether $F$ is closed. When $F$ is compact, in part (b)
we obtain the standard large deviation bound for $F$ with $-I(F)$ on the
right hand side. When $F$ is closed, in part (c) we obtain a different form of the standard large deviation upper bound; $-I(F)$ on the right hand side is
replaced by $-I(\overline{F})$, where $\overline{F}$ is the closure of $F$ in the compact space $\Y$.
When $F$ is compact, its closure in the compact space $\pnbbc$ is $F$ itself. In this case the large deviation upper bounds in parts (c) 
and (d) coincide.

\begin{thm}
\label{thm:ballstoldp}
For $N \in \N$ let $(\Omega_N,\mathcal{F}_N,Q_N)$ be a 
sequence of probability spaces. Let $\X$ be a complete, separable metric space, $\W$ a relatively compact subset of $\X$
that is not closed and thus not compact, and $\mathZ$
the closure of $\W$ in $\X$; thus $\mathZ$ is compact. Also let $Y_N$ be a sequence of random vectors mapping 
$\Omega_N$ into $\W$, and $I$ be a lower semicontinuous function mapping $\X$ into $[0,\infty]$. 
We assume the following two limits: for any open ball $B$ in $\W$
\be 
\label{eqn:repeatballs}
\lim_{N \goto \infty} \frac{1}{N} \log Q_N(Y_N \in B) = -I(B) 
\ee
and for any open ball $\hatb$ in $\mathZ$
\be 
\label{eqn:repeatballsup}
\lim_{N \goto \infty} \frac{1}{N} \log Q_N(Y_N \in \hatb \cap \W) = -I(\hatb \cap \W).
\ee

Then, as $N \goto \infty$, with respect to the measures $Q_N$,
the sequence $Y_N$ satisfies the LDP on $\W$ with rate function $I$ in the following sense.

{\em (a)} For any compact subset $F$ of $\W$ we have the large deviation upper bound
\[
\limsup_{N \goto \infty} \frac{1}{N} \log Q_N\{Y_N \in F\} \leq - I(F).
\]

{\em (b)} For any closed subset $F$ of $\W$ we have the large deviation upper bound
\[
\limsup_{N \goto \infty} \frac{1}{N} \log Q_N\{Y_N \in F\} \leq - I(\overline{F}),
\]
where $\overline{F}$ denotes the closure of $F$ in $\mathZ$.

{\em (c)} For any open subset $G$ of $\W$ we have the large deviation lower bound
\[
\liminf_{N \goto \infty} \frac{1}{N} \log Q_N(Y_N \in G\} \geq - I(G) .
\]
\end{thm}

Theorem \ref{thm:ldpthetankm} is an immediate consequence of this theorem, Theorem \ref{thm:ldlimitballs},
and Theorem \ref{thm:relentropy}. Part (a) of Theorem \ref{thm:ldlimitballs}
proves the large deviation limit
for any open ball in $\pnbc$, which corresponds to the limit (\ref{eqn:repeatballs}) in Theorem \ref{thm:ballstoldp}. Part (b)
 of Theorem \ref{thm:ldlimitballs} proves the large deviation limit
for $\hatb \cap \mathZ$, where $\hatb$ is any open ball in $\pnbbc$. This corresponds to the limit (\ref{eqn:repeatballsup}) 
in Theorem \ref{thm:ballstoldp}. In the application to 
Theorem \ref{thm:ldpthetankm} $\W$ is the relatively compact, nonclosed subset $\pnbc$ of $\X = \pnb$ and $\mathZ$ is the compact subset
$\pnbbc$ of $\pnb$. According to parts (a) and (b) of Theorem \ref{thm:relentropy},
$R(\cdot | \rhoalphabc)$ maps $\pnbc$ into $[0,\infty]$ and is lower semicontinuous on $\pnb$,
while part (d) of that theorem proves that $R(\cdot | \rhoalphabc)$ has compact level sets in $\pnbc$.
This last property of the relative entropy is needed for part (a) of Theorem \ref{thm:ldpthetankm}.

\skp
\noindent
{\bf Proof of Theorem \ref{thm:ballstoldp}.} We prove the three large deviation bounds in the order (c), (a), and (b).

\skp
(c) Let $G$ be any open subset of $\W$. 
We denote by $\tau$ the metric on $\X$. For any point $x \in G$ there exists $\ve > 0$ such that the open
ball $B_\tau(x,\ve) = \{y \in \W : \tau(x,y) < \ve\}$ is a subset of $G$. The limit (\ref{eqn:repeatballs}) implies that
\beas
\liminf_{N \goto \infty} \frac{1}{N} \log Q_N(Y_N \in G) & \geq &
\lim_{N \goto \infty} \frac{1}{N} \log Q_N(Y_N \in B_\tau(x,\ve)) \\
& = & -I(B_\tau(x,\ve)) \geq -I(x) .
\eeas
Since $x$ is an arbitrary point in $G$, it follows that
\[
\liminf_{N \goto \infty} \frac{1}{N} \log Q_N(Y_N \in G) \geq  -\inf_{x \in G} I(x) = -I(G) .
\]
This completes the proof of the large deviation lower bound for any open set $G$ in $\W$.

\skp
(a) Let $F$ be any compact subset of $\W$. 
We first prove the large deviation upper bound for $F$ under the assumption that $I(F) < \infty$.
The proof when $I(F) = \infty$ is given afterward.
We start by showing that for each $x \in F$
\be 
\label{eqn:ibxve}
\liminf_{\ve \goto 0^+} I(B_\tau(x,\ve)) \geq I(F) .
\ee
Let $\ve_n$ be any positive sequence converging to $0$, and take any $\delta > 0$. 
For any $n \in \N$ there exists $x_n \in B_\tau(x,\ve_n)$ such that $I(B_\tau(x,\ve_n)) + \delta \geq I(x_n)$. Since $x_n \goto x$, 
the lower semicontinuity of $I$ on $\W$ and the fact that $x \in F$ imply that
\[
\liminf_{n \goto \infty} I(B_\tau(x,\ve_n)) + \delta \geq \liminf_{n \goto \infty}I(x_n) \geq I(x) \geq I(F).
\]
Sending $\delta \goto 0$ yields (\ref{eqn:ibxve}) because $\ve_n$ is an arbitrary positive
sequence converging to 0. 

We now prove the large deviation upper bound in part (a). Take any $\eta > 0$. 
By (\ref{eqn:ibxve}) for each $x \in F$ there exists $\ve_x > 0$ such that
\[
I(B_\tau(x,\ve_x)) \geq I(F) - \eta .
\]
The open balls $\{B_\tau(x,\ve_x), x \in F\}$ cover $F$.
Since $F$ is compact, there
exist $T < \infty$ and finitely many points $x_i \in F, i = 1,2,\ldots,T$, such that $F\subset \bigcup_{i=1}^T
B_\tau(x_i,\ve_i)$, where $\ve_i = \ve_{x_i}$. It follows that
\[
\min_{i=1,2,\ldots,T} I(B_\tau(x_i,\ve_i)) \geq I(F) - \eta .
\]
By Lemma 1.2.15 in \cite{DemboZeitouni} and by the limit (\ref{eqn:repeatballs}) applied to $B = B_\tau(x_i,\ve_i)$ 
\bea
\label{eqn:morestuff}
\lefteqn{
 \limsup_{N \goto \infty}\frac{1}{N}\log Q_N\{Y_N \in F\}} \\ 
\nonumber
& & \leq \limsup_{N \goto \infty}\frac{1}{N}\log Q_N\left(Y_N \in \bigcup_{i=1}^T
B_\tau(x_i, \ve_i)\right) \\ 
\nonumber
&& \leq \limsup_{N \goto \infty}\frac{1}{N}\log \left(\sum_{i=1}^{T}Q_N(Y_N \in {B_\tau}(x_i,\ve_i))\right) \\
\nonumber & & = \max_{i=1,2,\ldots,T} \left( \limsup_{N \goto \infty}\frac{1}{N}\log Q_N(Y_N \in {B_\tau}(x_i,\ve_i)) \right) \\
\nonumber & & = - \min_{i=1,2,\ldots,T} I({B_\tau}(x_i,\ve_i)) \leq - I(F) + \eta .
\eea
Sending $\eta \goto 0$, we obtain
\[ 
\limsup_{N\goto \infty}\frac{1}{N}\log Q_N\{Y_N \in F\}\le -I(F) .
\]
This completes the proof of the large deviation upper bound for any compact subset $F$ of $\W$ under the assumption that $I(F) < \infty$.

We now assume that $I(F) = \infty$, which implies that $I(x) = \infty$ for each $x \in F$. The proof of the large deviation upper bound
when $I(F) = \infty$ rests on the assertion that for each $x \in F$ there exists $\ve_x > 0$ such that $I({B_\tau}(x,\ve_x)) = \infty$. 
Indeed, if this assertion were false, then there would exist a sequence $x_n \in \W$ satisfying $I(x_n) < \infty$ and $x_n \goto x$. Since $I$
is lower semicontinuous on $\W$, it would follow that $\liminf_{n \goto \infty} I(x_n) \geq I(x) = \infty$, which in turn would imply that $I(x_n) = \infty$. 
This contradiction completes the proof that for each $x \in F$ there exists $\ve_x > 0$ such that $I({B_\tau}(x,\ve_x)) = \infty$. As in the case 
when $I(F) < \infty$, the open balls $\{B_\tau(x,\ve_x), x \in F\}$ cover $F$.
Since $F$ is compact, there
exist $T < \infty$ and finitely many points $x_i \in F, i = 1,2,\ldots,T$, such that $F\subset \bigcup_{i=1}^T
B_\tau(x_i,\ve_i)$, where $\ve_i = \ve_{x_i}$. It follows that 
\[
\min_{i=1,2,\ldots,T} I(B_\tau(x_i,\ve_i)) = \infty = I(F).
\]
By the same steps as in (\ref{eqn:morestuff})
\[
\limsup_{N \goto \infty}\frac{1}{N}\log Q_N\{Y_N \in F\} \leq - \min_{i=1,2,\ldots,T} I({B_\tau}(x_i,\ve_i)) = -\infty = -I(F).
\]
This completes the proof of the large deviation upper bound for any compact subset $F$ of $\W$ when $I(F) = \infty$. The proof of part (a) is complete.

\skp
(b) Let $F$ be any closed subset of $\W$. We claim that $F$ equals $\overline{F} \cap \W$, where $\overline{F}$ is the closure of $F$ in $\mathZ$. Since $\mathZ$ is compact,
the closed subset $\overline{F}$ is also compact. 
Clearly $F \subset \overline{F} \cap \W$. On the other hand, any $x \in \overline{F} \cap \W$ 
is a limit point lying in $\W$ of a sequence $x_n$ in $F$. Since $F$ is closed in $\W$, any $x \in \hatf \cap \W$ lies in $F$. This completes the proof
that $F = \overline{F} \cap \W$. This is a special case of a general result in topology stated in Theorem 17.2 of \cite{Munkres}.

We first prove the large deviation upper bound for $F$ under the assumption that $I(\hatf) < \infty$.
The proof when $I(\hatf) = \infty$ is given afterward. The proof proceeds as in part (a), essentially
by replacing the balls $B_\tau(x,\ve)$ for $x \in \W$ by $\hatb_\tau(x,\ve) \cap \W$ for $x \in \mathZ$,
where $\hatb_\tau(x,\ve) = \{y \in \mathZ : \tau(x,y) < \ve\}$. 
As in the proof of part (a), we start by showing that for each $x \in \hatf$
\be 
\label{eqn:ibxvf}
\liminf_{\ve \goto 0^+}I(\hatb_\tau(x,\ve) \cap \W) \geq I(\hatf).
\ee
Let $\ve_n$ be any positive sequence converging to $0$, and take any $\delta > 0$. 
For any $n \in \N$ there exists $x_n \in \hatb_\tau(x,\ve_n) \cap \W$ such that 
$I(\hatb_\tau(x,\ve_n) \cap \W) + \delta \geq I(x_n)$. Since $x_n \goto x$, the lower semicontinuity of $I$ 
and the fact that $x \in \hatf$ imply that
\[
\liminf_{n \goto \infty} I(\hatb_\tau(x,\ve_n) \cap \W) 
+ \delta \geq \liminf_{n \goto \infty}I(x_n) \geq I(x) \geq I(\hatf).
\]
Sending $\delta \goto 0$ yields (\ref{eqn:ibxvf}) because $\ve_n$ is an arbitrary positive
sequence converging to 0. 

We now prove the large deviation upper bound in part (b). Take any $\eta > 0$. 
By (\ref{eqn:ibxvf}) for each $x \in \hatf$ there exists $\ve_x > 0$ such that
\[
I(\hatb_\tau(x,\ve_x) \cap \W) \geq I(\hatf) - \eta .
\]
The open balls $\{\hatb_\tau(x,\ve_x), x \in \hatf\}$ cover $\hatf$.
Since $\hatf$ is compact, there
exist $T < \infty$ and finitely many points $x_i \in \hatf, i = 1,2,\ldots,T$, such that $\hatf \subset \bigcup_{i=1}^T
\hatb_\tau(x_i,\ve_i)$, where $\ve_i = \ve_{x_i}$. 
It follows that
\[\min_{i=1,2,\ldots,T} I(\hatb_\tau(x_i,\ve_i) \cap \W) \geq I(\hatf) - \eta 
\]
and 
\[
\hatf \cap \W \subset \bigcup_{i=1}^T \left(\hatb_\tau(x_i,\ve_i) \cap \W\right).
\]
Since $F = \hatf \cap \W$, we have again by Lemma 1.2.15 in \cite{DemboZeitouni} 
\bea
\label{eqn:moremorestuff}
\lefteqn{
 \limsup_{N \goto \infty}\frac{1}{N}\log Q_N\{Y_N \in F\}} \\
\nonumber
 && = \limsup_{N \goto \infty}\frac{1}{N}\log Q_N\{Y_N \in \hatf \cap \W\} \\
\nonumber
&& \le \limsup_{N \goto \infty}\frac{1}{N}\log Q_N\left(Y_N \in \bigcup_{i=1}^T
\left(\hatb_\tau(x_i, \ve_i) \cap \W \right)\right) \\ 
\nonumber 
&&\le \limsup_{N \goto \infty}\frac{1}{N}\log \left(\sum_{i=1}^{T}Q_N(Y_N \in {\hatb_\tau}(x_i,\ve_i) \cap \W)\right) \\
\nonumber
&& = \max_{i=1,2,\ldots,T} \left( \limsup_{N \goto \infty}\frac{1}{N}\log Q_N(Y_N \in {\hatb_\tau}(x_i,\ve_i) \cap \W) \right).
\eea
We now apply the limit (\ref{eqn:repeatballsup}) to $\hatb \cap \W = \hatb_\tau(x_i,\ve_i) \cap \W$, obtaining
\bea
\label{eqn:moremoremorestuff}
\lefteqn{
 \limsup_{N \goto \infty}\frac{1}{N}\log Q_N\{Y_N \in F\}} \\
\nonumber 
&& \leq \max_{i=1,2,\ldots,T} \left( \limsup_{N \goto \infty}\frac{1}{N}\log Q_N(Y_N \in {\hatb_\tau}(x_i,\ve_i) \cap \W) \right) \\
\nonumber
&& = - \min_{i=1,2,\ldots,T} I(\hatb_\tau(x_i,\ve_i) \cap \W) \leq -I(\hatf) + \eta.
\eea
Sending $\eta \goto 0$, we obtain
\[ 
\limsup_{N\goto \infty}\frac{1}{N}\log Q_N\{Y_N \in F\}\le -I(\hatf) .
\]
This completes the proof of the large deviation upper bound for any closed subset $F$ of $\W$ under the assumption that $I(\hatf) < \infty$.

We now assume that $I(\hatf) = \infty$, which implies that $I(x) = \infty$ for each $x \in \hatf$. The proof of the large deviation upper bound
when $I(\hatf) = \infty$ rests on the assertion that for each $x \in \hatf$ there exists $\ve_x > 0$ such that $I({\hatb_\tau}(x,\ve_x) \cap \W) = \infty$. 
As in the proof of part (b), this assertion is a consequence of the lower semicontinuity of $I$. As in the proof of the large deviation
upper bound when $I(\hatf) < \infty$, the open balls $\{\hatb_\tau(x,\ve_x), x \in \hatf\}$ cover $\hatf$.
Since $\hatf$ is compact, there
exist $T < \infty$ and finitely many points $x_i \in \hatf, i = 1,2,\ldots,T$, such that 
$\hatf \subset  \bigcup_{i=1}^T \hatb_\tau(x_i,\ve_i)$, where $\ve_i = \ve_{x_i}$. It follows that 
\[
\min_{i=1,2,\ldots,T} I(\hatb_\tau(x_i,\ve_i)) = \infty = I(\hatf)
\]
and
\[
\hatf \cap \W \subset \bigcup_{i=1}^T
\hatb_\tau(x_i,\ve_i) \cap \W.
\]
By the same steps as in (\ref{eqn:moremorestuff}) and (\ref{eqn:moremoremorestuff})
\beas
\lefteqn{
\limsup_{N \goto \infty}\frac{1}{N}\log Q_N\{Y_N \in F\}} \\ 
& & = \limsup_{N \goto \infty}\frac{1}{N}\log Q_N\{Y_N \in \hatf \cap \W\} \\
& & \leq - \min_{i=1,2,\ldots,T} I({\hatb_\tau}(x_i,\ve_i) \cap \W) = -\infty = -I(\hatf).
\eeas
This completes the proof of the large deviation upper bound for any closed subset $F$ of $\W$ when $I(\hatf) = \infty$. The proof of part (b) as well
as the proof of the theorem are done. \ \ink

\skp
This paper contains four appendices. In appendix A we prove properties of the relative entropy needed in
the paper. Theorem \ref{thm:approximate} in appendix B states a basic approximation result 
that is applied in two crucial places in the paper. In appendix C we study a number of properties of the quantity $\alphabc$
appearing in part (a) of Theorem \ref{thm:mainestimate}. 
In appendix D we discuss why we impose the constraint involving $m = m(N)$ 
in the definitions of $\omeganbm$ in (\ref{eqn:omegankm}) and $\pnbm$ in (\ref{eqn:condprob}) and how, if this constraint could be eliminated, 
then our results could be formulated in a more natural way. 

\skp
\skp


\noi
\LARGE {\bf Appendices}
\vspace{-.2in}
\normalsize
\appendix

\renewcommand{\thesection}{\Alph{section}}
\renewcommand{\theequation}
{\Alph{section}.\arabic{equation}}
\renewcommand{\thedefn}
{\Alph{section}.\arabic{defn}}
\renewcommand{\theass}
{\Alph{section}.\arabic{ass}}

\section{Properties of Relative Entropy} 
\beginsec

We fix a nonnegative integer $b$ and a real number $c \in (b,\infty)$. Given $\theta$ a probability measure on $\Nb = \{n \in \Z : n \geq b\}$, the mean 
$\int_\N x \theta(dx)$ of $\theta$ is denoted by $\langle \theta \rangle$. In Theorem \ref{thm:relentropy}
we study properties of the relative entropy $R(\theta | \rhoalphab)$ and $R(\theta | \rhoalphabc)$ for $\theta$
in each of the following three spaces:
$\pnb$, the set of probability measures on $\N$; $\pnbc$, the set of $\theta \in \pnb$ satisfying
$\langle\theta\rangle = c$;
and $\pnbbc$, the set of $\theta \in \pnb$ satisfying $\langle\theta\rangle \in [b,c]$. The Prohorov metric
introduces a topology on $\pnb$ that is equivalent to the topology of weak convergence. These three
spaces have the following properties: $\pnb$ is a complete, separable metric space; $\pnbc$ is relatively compact, separable subset of $\pnb$ that is not
closed in $\pnb$ and therefore is not complete;
$\pnbbc$ is the closure of $\pnbc$ in $\pnb$ and is a compact, separable metric space. The properties of $\pnb$ are proved in Theorems 3.3.1 and Theorem 3.1.7
of \cite{EthierKurtz}, and the properties
of $\pnbc$ and $\pnbbc$ are proved in Theorem \ref{thm:pnc}. 

We recall that for $\alpha \in (0,\infty)$, $\rhoalphab$ denotes the Poisson distribution on $\N_b$ having components
\[
\rhoalphabj = \frac{1}{\zbalpha} \cdot \frac{\alpha^j}{j!} \mbox{ for } j \in \Nb,
\]
where $Z_0(\alpha) = e^\alpha$, and for $b \in \N$, $\zbalpha = e^\alpha - \sum_{j=0}^{b-1} \alpha^j/j!$.
According to part (a) of
Theorem \ref{thm:mainestimate} there exists a unique value $\alpha = \alphabc$ for which 
$\langle\rhoalphabc\rangle = c$; thus $\rhoalphabc$ lies in $\pnbc$. Assertion (ii) in part (f) of the next theorem
plays an important role in the main part of the paper. 
After the statement of Lemma \ref{lem:omegankm} we use this assertion 
to show that the arbitrary parameter $\alpha$ in Lemmas \ref{lem:deltankmnu} and \ref{lem:omegankm}
must have the value $\alphabc$ in Theorem \ref{thm:mainestimate}.

\begin{thm}
\label{thm:relentropy}

Fix a nonnegative integer $b$ and a real number $c \in (b,\infty)$. For any $\alpha \in (0,\infty)$ the relative entropy $R(\theta | \rhoalphab) = \sum_{j \in \Nb} \theta_j \log(\theta_j/\rhoalphabj)$ has the following 
properties.

{\em (a)} $R( \cdot | \rhoalphab)$ maps $\pnb$ into $[0,\infty]$, and for $\theta \in \pnb$,
$R(\theta | \rhoalphab) = 0$ if and only if $\theta = \rhoalphab$.  

{\em (b)} $R(\cdot | \rhoalphab)$ is a convex, lower semicontinuous function on $\pnb$. In other words, for $\theta$ and $\sigma$ in $\pnb$, $\lambda \in (0,1)$,
and $\thetan$ a sequence in $\pnb$ converging weakly to $\theta$
\[
R(\lambda \theta + (1-\lambda) \sigma | \rhoalphab) \leq \lambda R(\theta | \rhoalphab) + (1-\lambda) R(\sigma | \rhoalphab) 
\]
and
\[
\liminf_{N \goto \infty} R(\thetan | \rhoalphab) \geq R(\theta | \rhoalphab).
\]

{\em (c)} $R(\cdot |  \rhoalphab)$ is a strictly convex function on the set $A = \{\theta \in \pnb : R(\theta | \rhoalphab) < \infty\}$. In other words, if $\theta
\not = \sigma$ are two measures in $A$, then for $\lambda \in (0,1)$
\[
R(\lambda \theta + (1-\lambda) \sigma | \rhoalphab) < \lambda R(\theta | \rhoalphab) + (1-\lambda) R(\sigma | \rhoalphab).
\]

{\em (d)} $R(\cdot | \rhoalphab)$ has compact level sets in $\pnb$, in $\pnbbc$ and in $\pnbc$. In other words, 
for $\Y$ equal to any of these three spaces and any
$M < \infty$, the set $\{\theta \in \Y : R(\theta | \rhoalphab) \leq M\}$ is a compact subset of $\Y$. 

{\em (e)} Define 
\[
g(\alpha,b,c) = \log \zalphab - c \log \alpha - (\log \zalphabc - c \log \alphabc),
\] 
where $Z_0(\alpha) = e^\alpha$, and for $b \in \N$, $\zbalpha = e^\alpha - \sum_{j=0}^{b-1} \alpha^j/j!$.
Then for any $\theta \in \pnbc$
\[
R(\theta | \rhoalphab) = R(\theta | \rhoalphabc) + g(\alpha,b,c).
\]

{\em (f)} The following two assertions hold. 
\begin{itemize}
  \item[{\em (i)}] $R(\theta | \rhoalphab)$ attains its infimum over $\theta \in \pnbc$ at the unique
measure $\theta = \rhoalphabc$, and 
\[
\min_{\theta \in \pnbc} R(\theta | \rhoalphab) = R(\rhoalphabc | \rhoalphab) = g(\alpha,b,c).
\]
  \item[{\em (ii)}] For any $\theta \in \pnbc$, $R(\theta|\rhoalphab)$ is related to $R(\theta|\rhoalphabc)$ by the formula 
\[
R(\theta | \rhoalphab) - \min_{\theta \in \pnbc} R(\theta | \rhoalphab) = R(\theta | \rhoalphabc).
\]
\end{itemize}
\end{thm}

\noindent
{\bf Proof}.  
(a)--(c) These properties are proved in Lemma 1.4.1 and in part (b) of Lemma 1.4.3 in \cite{DupuisEllis}. 

\skp
(d) 
The fact that $R(\cdot|\rhoalphabc)$ has compact level sets in $\pn$ is proved in part (c) of Lemma 1.4.3 in \cite{DupuisEllis}. 
 According to 
part (b) of Theorem \ref{thm:pnc}, $\pnbbc$ is a compact subset of $\pnb$. Hence for any $M < \infty$
\[
\{\theta \in \pnbbc : R(\theta | \rhoalphab) \leq M\} = \{\theta \in \pnb : R(\theta | \rhoalphab) \leq M\} \cap \pnbbc
\]
is a compact subset of $\pnbbc$. This completes the proof that $R(\cdot | \rhoalphab)$ has compact level sets in $\pnbbc$. 

Because $\pnbc$ is not a closed subset of $\pnbbc$ [Thm.\ \ref{thm:pnc}(a)], the proof that $R(\cdot | \rhoalphab)$ has compact level sets in $\pnbc$
is more subtle. If $\thetanlower$ is any sequence in $\pnbc$ satisfying $R(\thetanlower | \rhoalphab) \leq M$, then since $\thetanlower \in \pnb$ and
$R(\cdot | \rhoalphab)$ has compact level sets in $\pnb$, there exists $\theta \in \pnb$ and a subsequence $\thetanlowerprime$ such that $\thetanlowerprime
\Rightarrow \theta$ and $R(\theta | \rhoalphab) \leq M$. To complete the proof that $R(\cdot | \rhoalphab)$ has compact level sets in $\pnbc$,
we must show that $\theta \in \pnbc$; i.e., that $\langle\theta\rangle = c$. By Fatou's Lemma
\[
\langle\theta\rangle \leq \liminf_{N \goto \infty} \langle\thetanlowerprime\rangle = c.
\]
In addition, for any $w\in (0,\infty)$
\[
\int_{\N_b} e^{wx} \rhoalphab(dx) = \sum_{j \in \Nb} e^{wj} \rhoalphabj = \frac{1}{\zalphab} \cdot \sum_{j \in \Nb} e^{wj} \frac{\alpha^j}{j!} 
\leq \frac{1}{\zalphab} \cdot \exp(\alpha e^w) < \infty.
\]
Lemma 5.1 in \cite{DonVar3} implies that the sequence $\thetanlowerprime$ is uniformly integrable; i.e., 
\[
\lim_{D \goto \infty} \sup_{n \in \N} \int_{\{x \in \N : x \geq D\}} x \thetanlowerprime(dx) = 0.
\]
These properties of $\theta$ and $\thetanlowerprime$
imply that $c = \lim_{n^\prime \goto \infty} \langle\thetanlowerprime\rangle = \langle\theta\rangle$ \cite[Appendix, Prop.\ 2.3]{EthierKurtz}. 
This completes the proof that $R(\cdot | \rhoalphab)$ has compact level sets in $\pnbc$. The proof of part (d) is finished. 

\skp
(e) For any $\theta \in \pnbc$ we have $\sum_{j\in\Nb}\theta_j = 1$ and $\sum_{j\in\Nb} j \theta_j = c$. Hence
\beas
R(\theta | \rhoalphab) & = & \sum_{j \in \Nb} \thetaj \log(\thetaj/\rhoalphabj) \\
& = &  \sum_{j \in \Nb} \thetaj \log(\thetaj/\rhoalphabcj) + \sum_{j \in \Nb} \thetaj \log(\rhoalphabcj/\rhoalphabj) \\
& = & R(\theta | \rhoalphabc) + \sum_{j \in \Nb} \thetaj \log\left(\frac{[\alphabc]^j}{\zalphabc j!} \cdot \frac{\zalphab j!}{\alpha^j} \right) \\
& = & R(\theta | \rhoalphabc) + \sum_{j \in \Nb} \thetaj \log(\zalphab/\zalphabc) + \sum_{j \in \Nb} j \thetaj \log(\alphabc/\alpha) \\
& = & R(\theta | \rhoalphabc) + \log(\zalphab/\zalphabc) + c \log(\alphabc/\alpha) \\
\\ & = & R(\theta | \rhoalphabc) + g(\alpha,b,c).
\eeas
This completes the proof of part (e). 

\skp
(f) (i) Since $R(\cdot|\rhoalphab)$ has compact level sets in $\pnbc$, it attains its infimum over $\pnbc$. By part (a) $R(\cdot | \rhoalphabc)$
attains its minimum value of 0 over $\pnbc$ at the unique measure $\rhoalphabc$. Hence part (e) implies that the minimum value
of $R(\cdot | \rhoalphab)$ over $\pnbc$ equals
\beas
\min_{\theta \in \pnbc} R(\theta | \rhoalphab) & = & \min_{\theta \in \pnbc} R(\theta | \rhoalphabc) + g(\alpha,b,c) \\
& = & g(\alpha,b,c) = R(\rhoalphabc|\rhoalphabc) + g(\alpha,b,c) = R(\rhoalphabc | \rhoalphab).
\eeas
The last equality follows by applying part (e) with $\theta = \rhoalphabc$. This display shows that
$R(\cdot | \rhoalphab)$ attains its infimum over $\pnbc$ at $\rhoalphabc$. 
Let us assume that $R(\cdot | \rhoalphab)$ attains its infimum over $\pnbc$ at another measure $\theta^\star \not = \rhoalphabc$.
Then for any $\lambda \in (0,1)$, we have $\lambda \rhoalphabc + (1-\lambda) \theta^\star \in \pnbc$. The strict convexity of $R(\cdot|\rhoalphab)$
in part (c) yields 
\beas
\min_{\theta \in \pnbc} R(\theta|\rhoalphab) & \leq & R(\lambda \rhoalphabc + (1-\lambda) \theta^\star | \rhoalphab)\\ 
& < & \lambda R(\rhoalphabc | \rhoalphab) + (1-\lambda) R(\theta^\star | \rhoalphab) =
\min_{\theta \in \pnbc} R(\theta|\rhoalphab).
\eeas
The equality of the extreme terms contradicts the strict inequality, proving that $R(\cdot | \rhoalphab)$ attains its infimum over $\pnbc$ at 
the unique measure $\rhoalphabc$. 
This completes the proof of assertion (i) in part (f). 

(ii) By assertion (i) $\min_{\theta \in \pnbc} R(\theta|\rhoalphab) = g(\alpha,b,c)$. Substituting this into part (e)
yields assertion (ii). This completes the proof of part (f). The proof of Theorem \ref{thm:relentropy} is done. \ \ink

\skp

This completes our discussion of properties of the relative entropy. The main theorem in appendix B is a basic approximation result that is 
applied in two crucial places in the paper.


\section{Approximating \boldmath $\theta \in \pnbc$ \unboldmath by \boldmath $\thetan \in \bnbm$ \unboldmath} 
\beginsec

Fix a nonnegative integer $b$ and a rational number $c \in (b,\infty)$. $\pnbc$ is the set of probability measures on $\Nb =
\{n \in \Z : n \geq b\}$ having mean $c$. 
We recall the definitions of the sets $\anbm$ and $\bnbm$, which are introduced at the beginning of section \ref{section:proof1}:
\[
\anbm = \left\{\nu = \{\nuj, j \in \Nb\} \in \N_0^\N : \sum_{j \in \Nb}\nu_j = N, \ \sum_{j \in \Nb}j \nu_j = K, 
\ \mbox{and} \ |\nu|_+ \leq m = m(N)\right\} 
\]
and 
\[
\bnbm = \{\theta \in \pnbc : \theta_j = \nu_j/N \mbox{ for } j \in \Nb \mbox{ for some } \nu \in \anbm\} .
\]
In the formula defining $\anbm$, $\N_0$ is the set of nonnegative integers and $|\nu|_+ = \mbox{card}\{j \in \Nb : \nu_j \geq 1\}$.
The quantities $K$ and $m$ are functions of $N$ as $N \goto \infty$: $K = Nc$, and $m$ is the function $m(N)$ appearing in the definition
of $\omeganbm$ in (\ref{eqn:omegankm}) and satisfying
$m(N) \goto \infty$ and $m(N)^2/N \goto 0$ as $N \goto \infty$. 

Our goal in this appendix is to prove the approximation theorem, Theorem \ref{thm:approximate}, and Corollary \ref{cor:dense}. The theorem 
is applied in two crucial places in the paper. It is first applied 
near the end of the proof of Lemma \ref{lem:omegankm} to prove the limit in (\ref{eqn:complicated}) 
and thus to complete the proof of that lemma. 
Theorem \ref{thm:approximate} is also needed to verify hypothesis (iv) in Theorem \ref{thm:balls} in the setting of Theorem \ref{thm:ldlimitballs}.
Theorem \ref{thm:balls} is applied to lift the local large deviation
estimate in part (b) of Theorem \ref{thm:mainestimate} to the large deviation limit for open balls and certain other subsets in
Theorem \ref{thm:ldlimitballs}. 

Because $R(\cdot | \rhoalphab)$ is lower semicontinuous on $\pnb$ [Thm.\ \ref{thm:relentropy}(b)], the weak
convergence in part (a) of the next theorem implies that $\liminf_{N \goto \infty} R(\thetan | \rhoalphab) \geq R(\theta|\rhoalphab)$. 
The proof of the convergence
$R(\thetan|\rhoalphab) \goto R(\theta|\rhoalphab)$ in part (b) requires the finiteness of $R(\theta|\rhoalphab)$
and special properties of the sequence $\thetan$ proved in Lemma \ref{lem:nujstar}.

\begin{thm}
\label{thm:approximate}
Fix a nonnegative integer $b$ and a rational number $c \in (b,\infty)$, and let $\theta$ be any probability measure in $\pnbc$.
Let $m$ be the function $m(N)$ appearing in the definition of $\omeganbm$ in {\em (\ref{eqn:omegankm})} and satisfying 
$m(N) \goto \infty$ and $m(n)^2/N \goto 0$ as $N \goto \infty$.  
Then for any $\alpha \in (0,\infty)$ there exists a sequence $\thetan \in \bnbm$ for which the following properties
hold.

{\em (a)} $\thetan \Rightarrow \theta$ as $N \goto \infty$.

{\em (b)} If $R(\theta | \rhoalphab) < \infty$, then $R(\thetan | \rhoalphab) \goto R(\theta | \rhoalphab)$
as $N \goto \infty$.
\end{thm}

We also need the following corollary, which is applied to verify hypothesis (iii) in Theorem \ref{thm:balls} in the setting of Theorem \ref{thm:ldlimitballs}.
It also shows that $\pnbc$ is separable, a fact needed in parts (a) and (b) of Theorem \ref{thm:pnc}. 

\begin{cor}
\label{cor:dense}
Fix a nonnegative integer $b$ and a rational number $c \in (b,\infty)$. Let $m$ be the function $m(N)$ 
appearing in the definition of $\omeganbm$ in {\em (\ref{eqn:omegankm})} and satisfying
$m(N) \goto \infty$ and $m(N)^2/N \goto 0$ as $N \goto \infty$. 
Then there exists a countable dense subset of $\pnbc$ consisting of $\theta \in \pnbc$ for which $R(\theta | \rhoalphabc) < \infty$. 
This countable dense subset is $\cup_{N \in \N} \bnbm$, where $\bnbm$ is defined at the beginning of this section. It follows that $\pnbc$ is separable.
\end{cor}

\noindent
{\bf Proof.} Given any $\theta \in \pnbc$ and any $\ve > 0$, let $B_\pi(\theta,\ve)$ denote the open ball 
with center $\theta$ and radius $\ve$ defined in terms of the Prohorov metric $\pi$.
We apply part (a) of Theorem \ref{thm:approximate} with $\alpha = \alphabc$. Since the measures $\thetan$ constructed in part (a) of that theorem 
converge weakly to $\theta$,
for all sufficiently large $N$ we have $\thetan \in B_\pi(\theta,\ve)$. 
The fact that only finitely many of the components $\theta^{(N)}_j$ are nonzero implies that $R(\thetan | \rhoalphabc) < \infty$ for all $N$.
Since $\cup_{N \in \N} \bnbm$ is a countable set, the proof is complete. \ \ink

\skp
Given $\theta \in \pnbc$, we determine a sequence $\nu^{(N)} \in \anbm$ such that the probability measures $\thetan$ 
with components $\thetanj = \nunj/N$ have the properties stated in parts
(a) and (b) of Theorem \ref{thm:approximate}. We start by defining 
\[
\jstar = \min\{j \in \Nb : \theta_j > 0\} .
\]
For example, for the Poisson distribution $\rhoalphabc$ defined in part (a) of 
Theorem \ref{thm:ldpthetankm}, $\jstar = b$ since for $j \in \nb$ all the components $\rhoalphabcj$ are positive.

We next define the components $\nunj$ of $\nun$ for all $j \in \Nb$ except for the two values $j = \jstar$ and $j = \jstarplusone$. The two
components corresponding to these two values of $j$ will then be defined so that $\nun$ satisfies the two summation constraints in the definition of $\anbm$. In order to simplify the notation, the components $\nunj$ are written as $\nuj$. For $x \in \R$ we denote by $\lfloor{x}\rfloor$ the largest integer less than or equal to $x$. The definition of the components 
is the following:
\be 
\label{eqn:nunj} 
\nuj = \left\{ \begin{array}{cl} 0 & \mbox{ if } b \leq j \leq \jstarminusone \\
\lfloor N\thetaj \rfloor & \mbox{ if } \jstarplustwo \leq j \leq \jstarmminusone \\
0 & \mbox{ if } j \geq \jstarm. 
\end{array}
\right.
\ee

We make a few simple observations. If $\jstar = b$, then the first line of this definition is vacuous. For $\jstarplustwo \leq j \leq \jstarmminusone$
\be 
\label{eqn:nujinequality}
\max\!\left(\thetaj - \frac{1}{N},0\right) \leq \frac{\nuj}{N} \leq \thetaj 
\mbox{ for all } N \ \mbox{ and } \ \lim_{N \goto \infty} \frac{\nuj}{N} = \thetaj.
\ee
In addition, for $b \leq j \leq j^\star -1$, we have $\nuj/N = 0 = \thetaj$. 
If for some $j$ satisfying $\jstarplustwo \leq j \leq \jstarmminusone$ we have $\thetaj = 0$, then $\nuj = 0$.

We now define $\nuj$ for $j = \jstar$ and $j = \jstarplusone$ so that $\nuj/N \goto \thetaj$ for these two values and so that the following two
summation constraints in the definition of $\anbm$ are valid:
\be 
\label{eqn:2formulas}
\sum_{j \in \Nb}\nu_j = N \ \mbox{ and } \ \sum_{j \in \Nb}j \nu_j = K .
\ee
With these definitions 
of $\nujstar$ and $\nujstarplusone$, we have $|\nu|_+ \leq m$. According to part (d) of Lemma \ref{lem:nujstar}, 
the resulting vector $\nu$ lies in $\anbm$ for all sufficiently large $N$. 

In order to keep the notation manageable, we introduce the set of $m-2$ indices
\[
\Phijstarm = \{j \in \Nb : \jstarplustwo \leq j \leq \jstarmminusone\} .
\]
Since $\nuj = 0$ for $b \leq j \leq \jstarminusone$ and 
for $j \geq \jstarm$,
the two equalities in (\ref{eqn:2formulas}) can be rewritten in the form
\be 
\label{eqn:sum1}
\nujstar + \nujstarplusone  =  N - \sum_{j \in \Phijstarm} \nuj 
\ee
and
\be 
\label{eqn:sum2}
\jstar \nujstar + (\jstarplusone) \nujstarplusone  =  K - \sum_{j \in \Phijstarm} j \nuj .
\ee
These are two linear equations for the two unknowns $\nujstar$ and $\nujstarplusone$. Solving them for the two unknowns and inserting 
$\nuj = \lfloor N\thetaj \rfloor$ for $j \in \Phijstarm$, we obtain the following definitions
of $\nujstar$ and $\nujstarplusone$:
\bea 
\label{eqn:nujstar}
\nujstar & = & (\jstarplusone)N - K + \sum_{j \in \Phijstarm} j \nuj 
- (\jstarplusone) \sum_{j \in \Phijstarm} \nuj \\
\nonumber
& = & (\jstarplusone)N - K + \sum_{j \in \Phijstarm} j \lfloor N\thetaj \rfloor 
- (\jstarplusone) \sum_{j \in \Phijstarm} \lfloor N \thetaj \rfloor
\eea
and 
\bea
\label{eqn:nujstarplusone}
\nujstarplusone & = & K - \jstar N - \sum_{j \in \Phijstarm} j \nuj  + 
\jstar \sum_{\Phijstarm} \nuj \\
\nonumber
& = & K - \jstar N - \sum_{j \in \Phijstarm} j \lfloor N \thetaj \rfloor  + 
\jstar \sum_{\Phijstarm} \lfloor N \thetaj \rfloor .
\eea

The next lemma states a number of facts about $\nu_j$ for $j \in \Nb$ 
that are needed to prove Theorem \ref{thm:approximate}. Parts (a) and (b) give upper and lower bounds on $\nujstar$ and $\nujstarplusone$
that follow from (\ref{eqn:nujstar}) and (\ref{eqn:nujstarplusone}). 
The reason for imposing the condition that $m^2/N \goto 0$ as $N \goto \infty$ in Theorem \ref{thm:approximate} is the
appearance of this quantity as an error term in parts (a) and (b).
Part (c) focuses on the convergence of $\nu_j/N$ to $\theta_j$
for $\jstar \leq j \leq \jstar + m -1$. Part (d) shows that 
for all sufficiently large $N$ the vector $\nu^{(N)}$ with components $\nu_j$ is an element of $\anbm$ and 
the measure $\theta^{(N)}$ with components $\theta^{(N)}_j = \nu_j/N$ for $j \in \Nb$
is an element of $\bnbm \subset \pnbc$.
In order to prove part (b) of Theorem \ref{thm:approximate}
concerning the convergence $R(\thetan | \rhoalphab) \goto R(\theta | \rhoalphab)$,
we will use the fact, stated in part (e), that for all $j \in \Nb$ satisfying $j \not = \jstarplusone$
we have $\theta^{(N)}_j = \nu_j/N \leq \theta_j$ for all $N$.
The conclusion of part (f) is that such a bound does not exist for $j = \jstarplusone$ and that in general there does not exist $M < \infty$ such that for any $N \in \N$, $\nu_{\jstarplusone}/N \leq M \theta_{\jstarplusone}$.

\begin{lem}
\label{lem:nujstar}
Fix a nonnegative integer $b$ and a rational number $c \in (b,\infty)$, and let $\theta$ be any probability measure in $\pnbc$. 
Let $m$ be the function $m(N)$ appearing in the definition of $\omeganbm$ in {\em (\ref{eqn:omegankm})} 
and satisfying $m(N) \goto \infty$ and $m(n)^2/N \goto 0$ as $N \goto \infty$. 
We define $\betam = \sum_{j \geq \jstarm} \thetaj$ and $\gammam = \sum_{j \geq \jstarm} j \thetaj$;
since $\theta \in \pnbc$, $\betam \goto 0$ and $\gammam \goto 0$ as $N \goto \infty$. The following conclusions hold.

{\em (a)} $\nujstar$ satisfies the inequalities
\[
N \thetajstar \geq \nujstar \geq N \left(\thetajstar + (\jstarplusone) \betam - \gammam - \frac{m^2}{N}\right) .
\]

{\em (b)} $\nujstarplusone$ satisfies the inequalities
\[
N\left( \thetajstarplusone + \gammam - \jstar\betam + \frac{m^2}{N}\right) \geq \nujstarplusone \geq 
N ( \thetajstarplusone + \gammam - \jstar\betam) \geq N \thetajstarplusone .
\]

{\em (c)} For all $j \in \Nb$ we have $\lim_{N \goto \infty} \thetanj = \lim_{N \goto \infty} \nu_j/N = \theta_j$.

{\em (d)} For all sufficiently large $N$ the vector $\nun$ with components $\nuj$ defined in 
{\em (\ref{eqn:nunj})}, {\em (\ref{eqn:nujstar})}, and {\em (\ref{eqn:nujstarplusone})} is an element of $\anbm$.
Hence for all sufficiently large $N$ the  
measure $\theta^{(N)}$ with components $\theta^{(N)}_j = \nu_j/N$ for $j \in \Nb$
is an element of $\bnbm \subset \pnbc$.

{\em (e)} 
For all $j \in \Nb$ satisfying $j \not  = \jstarplusone$ we have $\theta^{(N)}_j = \nu_j/N \leq \theta_j$ for all $N \in \N$.

{\em (f)} The upper bound $\theta^{(N)}_{\jstarplusone} = \nu_{\jstarplusone}/N \leq \theta_{\jstarplusone}$
does not hold for any $N$. 
On the other hand, 
if $\thetajstarplusone > 0$, then for all sufficiently large $N$ we have $\nujstarplusone/N \leq 2 \thetajstarplusone$. 
However, if $\thetajstarplusone = 0$, then in general there does not exist $M < \infty$ such that for any $N \in \N$, $\nu_{\jstarplusone}/N \leq M \theta_{\jstarplusone}$.

\end{lem}

\noindent
{\bf Proof.} (a) We first prove the lower bound.  According to (\ref{eqn:nujinequality}), $\nuj \geq N (\thetaj - 1/N)$ for all $j \in \Phijstarm$.
Since for all $j \in \Phijstarm$ we have $j > \jstarplusone$, the first line of (\ref{eqn:nujstar}) implies that
\bea
\label{eqn:3lines}
\nujstar & = & N\left[\jstarplusone - c + \sum_{j \in \Phijstarm} (j - \jstar -1) \frac{\nuj}{N}\right] \\
\nonumber
& \geq & N\left[\jstarplusone - c + \sum_{j \in \Phijstarm} (j - \jstar -1) \left(\thetaj - \fraconen\right)\right] \\
\nonumber
& = & \, N\left[(\jstarplusone)\left(1 - \sum_{j \in \Phijstarm} \thetaj\right) - c + \sum_{j \in \Phijstarm} j\thetaj
- \sum_{j \in \Phijstarm} (j - \jstar -1) \frac{1}{N} \right] . 
\eea
We now use the facts that $\theta_j = 0$ for $b \leq j \leq \jstarminusone$, $\sum_{j \in \Nb} \thetaj = 1$, $\sum_{j \in \Nb} j \thetaj = c$
to calculate
\bea
\label{eqn:tocalculate1}
\sum_{j \in \Phijstarm} j\thetaj & = & \sum_{j = \jstarplustwo}^{\jstarmminusone} j \thetaj \\
\nonumber
& = & \sum_{j \in \N} j \thetaj - \jstar\thetajstar - (\jstarplusone)\thetajstarplusone - \gamma_m \\
\nonumber
& = & c - \jstar\thetajstar - (\jstarplusone)\thetajstarplusone - \gamma_m
\eea
and
\bea
\label{eqn:tocalculate2}
\sum_{j \in \Phijstarm} \thetaj & = & \sum_{j = \jstarplustwo}^{\jstarmminusone} \thetaj \\
\nonumber 
& = & \sum_{j \in \N} \thetaj - \thetajstar - \thetajstarplusone - \beta_m \\
\nonumber
& = & 1- \thetajstar - \thetajstarplusone -  \beta_m .
\eea
In addition
\be 
\label{eqn:tocalculate3}
\sum_{j \in \Phijstarm} (j - \jstar -1) \frac{1}{N} = \frac{1}{N} \sum_{j = 1}^{m-2} j 
= \frac{(m-2)(m-1)}{2N} \leq \frac{m^2}{2N} .
\ee
Substituting (\ref{eqn:tocalculate1}), (\ref{eqn:tocalculate2}), and (\ref{eqn:tocalculate3}) into the last expression in (\ref{eqn:3lines}), we conclude that
\[
\nujstar \geq N \left[\thetajstar + (\jstarplusone) \betam - \gammam - \frac{m^2}{N}\right] .
\]
This is the lower bound in part (a). 

We now prove the upper bound in part (a). According to (\ref{eqn:nujinequality}), $\nuj \leq N\thetaj$ for all $j \in \Phijstarm$.
Since for all $j \in \Phijstarm$ we have $j > \jstarplusone$, the first line of (\ref{eqn:nujstar}) implies that
\beas
\nujstar & = & N\left(\jstarplusone - c + \sum_{j \in \Phijstarm} (j - \jstar -1) \frac{\nuj}{N}\right) \\
& \leq & N\left(\jstarplusone - c + \sum_{j \in \Phijstarm} (j - \jstar -1)\thetaj\right) .
\eeas
Except for the absence of the term containing $1/N$, this is the same expression that appears in the second line of (\ref{eqn:3lines}). Hence
by a calculation similar to that yielding the lower bound in part (a)
\[
\nujstar \leq N (\thetajstar + (\jstarplusone) \betam - \gammam) .
\]
We now use the fact that
\[
(\jstarplusone) \betam - \gammam = (\jstarplusone)\sum_{j \geq \jstarplusone} \thetaj - \sum_{j \geq \jstarplusone} j \thetaj \leq 0 . 
\]
Substituting this inequality into the preceding display shows that $\nujstar \leq N\thetajstar$ \,.
This is the upper bound in part (a). The proof of part (a) is complete.

\skp
(b) We first prove the upper bound. According to (\ref{eqn:nujinequality}), $\nuj \geq N(\thetaj - 1/N)$ for all $j \in \Phijstarm$. Since for all
$j \in \Phijstarm$ we have $j > \jstar$, the first line of (\ref{eqn:nujstarplusone}) implies that
\bea
\label{eqn:3linesx}
\nujstarplusone & = & N\left[c - \jstar - \sum_{j \in \Phijstarm} (j - \jstar) \frac{\nuj}{N}\right] \\
\nonumber
& \leq & N\left[c - \jstar - \sum_{j \in \Phijstarm} (j - \jstar) \left(\thetaj - \fraconen\right)\right] \\
\nonumber
& = & \, N\left[c - \sum_{j \in \Phijstarm} j\thetaj - \jstar\left(1 - \sum_{j \in \Phijstarm} \thetaj\right)
+ \sum_{j \in \Phijstarm} (j - \jstar) \frac{1}{N} \right] . 
\eea
As in the proof of (\ref{eqn:tocalculate3}),
\[
\sum_{j \in \Phijstarm} (j - \jstar) \frac{1}{N} \leq \frac{m^2}{2N} .
\]
Substituting this inequality as well as the equalities in (\ref{eqn:tocalculate1}) and (\ref{eqn:tocalculate2}) into the last expression in (\ref{eqn:3linesx}),
we conclude that
\[
\nujstarplusone \leq N\left(\thetajstarplusone + \gammam - \jstar \betam + \frac{m^2}{2N}\right) .
\]
This is the upper bound in part (b). 

We now prove the lower bound in part (b). According to (\ref{eqn:nujinequality}), $\nuj \leq N\thetaj$ for all $j \in \Phijstarm$. Since
for all $j \in \Phijstarm$ we have $j > \jstar$, the first line of (\ref{eqn:3linesx}) implies that
\beas
\nujstarplusone & = & N\left[c - \jstar - \sum_{j \in \Phijstarm} (j - \jstar) \frac{\nuj}{N}\right] \\
\nonumber
& \geq & N\left[c - \jstar - \sum_{j \in \Phijstarm} (j - \jstar)\thetaj\right] \\
\eeas
Except for the absence of the term containing $1/N$, this is the same expression that appears in the second line of (\ref{eqn:3linesx}). 
Hence by a calculation similar to that yielding the upper bound in part (b)
\[
\nujstarplusone \geq N(\thetajstarplusone + \gammam - \jstar \betam) .
\]
This is the second inequality in part (b). We now use the fact that
\[
N(\thetajstarplusone + \gammam - \jstar \betam) = N \thetajstarplusone + 
N \sum_{j \geq \jstarplusone} (j - \jstar) \thetaj \geq N \thetajstarplusone .
\]
This is the third inequality in part (b). The proof of part (b) is complete.

\skp
(c) For $j = \jstar$ and $j = \jstarplusone$ the limits $\lim_{N \goto \infty} \nuj/N = \thetaj$ 
are immediate consequences of parts (a) and (b) since each of the quantities $\betam$, $\gammam$, and $m^2/N$ converge to 0
as $N \goto \infty$. For $j \in \N$ satisfying $j \geq \jstar + 2$ the limit 
$\lim_{N \goto \infty} \nuj/N = \thetaj$ follows from (\ref{eqn:nujinequality}) and the fact that $m \goto \infty$ as $N \goto \infty$. 
Finally, for $j \in \Nb$ satisfying $b \leq j \leq \jstar - 1$, $\nuj/N = 0 = \thetaj$. The proof of part (c) is complete.

\skp
(d) According to (\ref{eqn:nunj}), for all $j \in \Nb$ satisfying $j \not = \jstar, \jstar+1$
we have $\nuj \in \N_0$ for all $N$. We now consider
$\nujstar$. As $N \goto \infty$, each of the quantities $\betam$, $\gammam$, and $m^2/N$ converge to 0. Since $\thetajstar > 0$, it follows 
from the lower bound in part (a) of this lemma that $\nujstar > 0$
for all sufficiently large $N$. The definition of $\nujstar$ in (\ref{eqn:nujstar}) shows that $\nuj$ is an integer for all $N$. It follows that
$\nujstar \in \N$ for all sufficiently large $N$. Finally we consider $\nujstarplusone$. 
The lower bound in part (b) of this lemma shows that $\nujstarplusone \geq 0$. The definition of  
$\nujstarplusone$ in (\ref{eqn:nujstarplusone}) shows that $\nujstarplusone$ is an integer for all $N$. It follows that $\nujstarplusone \in \N_0$ for all $N$.
We conclude that for all sufficiently large $N$
the vector $\nun$ is an element of $\N_0^N$. In addition, since $\nuj = 0$ for all $j \in \Nb$ satisfying
$b \leq j \leq \jstar - 1$ and $j \geq \jstar + m$, we have $|\nun|_+ \leq m$; i.e., at most of the components $\nuj$ are positive.
These correspond to the indices $j \in \Nb$ satisfying $\jstar \leq j \leq \jstar + m -1$. 
If the definitions of $\nujstar$ and $\nujstarplusone$ in (\ref{eqn:nujstar}) and (\ref{eqn:nujstarplusone}) 
are substituted into (\ref{eqn:sum1}) and (\ref{eqn:sum2}), then we see that the components
$\nunj$ satisfy the two equality constraints in the definition of $\ankm$ for all $N$. It follows
that $\nun \in \ankm$ for all sufficiently large $N$. We also conclude that the measure $\thetan$
having components $\thetanj = \nuj/N$ for $j \in \Nb$ is an element of $\bnkm \subset \pnc$ for all sufficiently large 
$N$. The proof of part (d) is complete. 

\skp
(e) For $j = \jstar$ and all $N$, we have $\nujstar/N \leq \theta_{\jstar}$ by the upper bound in part (a) 
of Lemma \ref{lem:nujstar}.
For all $j \in \Nb$ satisfying $\jstar + 2 \leq j \leq \jstar + m -1$ and for all $N$, we have 
$\nuj/N \leq \thetaj$ by (\ref{eqn:nujinequality}). Finally, by (\ref{eqn:nunj}) 
for all $j \in \Nb$ satisfying $b \leq j \leq \jstar - 1$ and $j \geq \jstar + m$ and for all $N$ we have
$\nunj/N  = 0 \leq \thetaj$. The proof of part (e) is complete.

\skp
(f) Assume that $\thetajstarplusone > 0$. By the upper bound in part (b) of this lemma, $\gammam - \jstar \betam + m^2/N \goto 0$ 
as $N \goto \infty$. Hence for all sufficiently large $N$, $\nujstarplusone/N \leq 2 \thetajstarplusone$. 
However, even if $\thetajstarplusone > 0$. the upper bound $\nujstarplusone/N \leq \thetajstarplusone$ cannot hold for
any $N$ because of the three additional terms in the upper bound in part (b); while $\gammam$ and $\betam$ can be 0 for 
sufficiently large $N$, the term $m^2/N > 0$ for all $N$. This proves the first two assertions in part (f).  
Concerning the third assertion, let us see how the bound $\nu_{\jstarplusone}/N \leq M \theta_{\jstarplusone}$ can fail. We assume that $\theta_{\jstarplusone} = 0$
and that there exists a subsequence $j^\prime \goto \infty$ such that $\theta_{j^\prime} > 0$ along this subsequence. By the lower bound in 
part (a) of this lemma 
\[
\nujstarplusone \geq N (\gammam - \jstar\betam) = N\left(\sum_{j \leq \jstar + m} (j - \jstar) \thetaj\right).
\]
Since $\theta_{j^\prime} > 0$ along the subsequence $j^\prime \goto \infty$, it follows that for all $N \in \N$ and all 
$j^\prime$
\[
\nujstarplusone \geq N(j^\prime - \jstar) \theta_{j^\prime} > 0.
\]
Since $\theta_{\jstarplusone} = 0$ and $\nu_{\jstarplusone}/N > 0$ for all $N \in \N$, the bound $\nu_{\jstarplusone}/N \leq M\theta_{\jstarplusone}$
cannot hold for any $M < \infty$. 
This completes the proof of part (f). The proof of Lemma \ref{lem:nujstar} is done. \ink

\skp
We are now ready to prove Theorem \ref{thm:approximate}. 
Given $\theta \in \pnbc$, $\thetan$ in this theorem is the sequence with components 
$\thetan_j = \nuj/N$ for $j \in \Nb$.  
 The quantities $\nuj = \nunj$ are defined in (\ref{eqn:nunj}),
(\ref{eqn:nujstar}), and (\ref{eqn:nujstarplusone}). In the proof of the theorem we work
with sufficiently large $N \in \N$ guaranteeing, 
according to part (d) of Lemma \ref{lem:nujstar}, that 
$\thetan$ is a probability measure lying in $\bnbm \subset \pnbc$.

\skp
\noi 
{\bf Proof of part (a) of Theorem \ref{thm:approximate}.}  We prove that $\thetan \Rightarrow \theta$ by showing that for any bounded function $f$ mapping $\Nb$
into $\R$
\[
\lim_{N \goto \infty} \int_{\Nb} f d\thetan = \lim_{N \goto \infty} \sum_{j \in \Nb} f(j) \thetan_j = \sum_{j \in \Nb} f(j) \thetaj = \int_{\Nb} f d\theta .
\]
We use the facts that $\nuj = 0 = \thetaj$ for $b \leq j \leq \jstarminusone$, $\nuj = 0$ for $j \geq \jstar + m$, and 
\[
\max_{\jstarplustwo \leq j \leq \jstarmminusone} \left|\frac{\nuj}{N} - \thetaj\right| \leq \fraconen \,.
\]
These facts, which follow from (\ref{eqn:nunj}) and (\ref{eqn:nujinequality}), give the upper bound
\beas
\lefteqn{
\left|\sum_{j \in \Nb} f(j) \thetan_j - \sum_{j \in \Nb} f(j) \thetaj\right| } \\
&& \leq |f(\jstar) \left|\frac{\nujstar}{N} - \thetajstar\right| + |f(\jstarplusone)| \left|\frac{\nujstarplusone}{N} - \thetajstarplusone\right| \\
&& \hspace{.5in} + \|f\|_\infty \sum_{j=\jstarplustwo}^{\jstarmminusone} \left|\frac{\nuj}{N} - \thetaj\right| + \|f\|_\infty \sum_{j \geq \jstarm} \thetaj \\
&& \leq |f(\jstar) \left|\frac{\nujstar}{N} - \thetajstar\right| + |f(\jstarplusone)| \left|\frac{\nujstarplusone}{N} - \thetajstarplusone\right| \\
&& \hspace{.5in} + \|f\|_\infty (m-2) \left(\max_{\jstarplustwo \leq j \leq \jstarmminusone} 
\left|\frac{\nuj}{N} - \thetaj\right|\right) + \|f\|_\infty \sum_{j \geq \jstarm} \thetaj \\
&& \leq |f(\jstar) \left|\frac{\nujstar}{N} - \thetajstar\right| + |f(\jstarplusone)| \left|\frac{\nujstarplusone}{N} - \thetajstarplusone\right| \\
&& \hspace{.5in} + \|f\|_\infty \frac{m}{N} + \|f\|_\infty \sum_{j \geq \jstarm} \thetaj .
\eeas
By part (c) of Lemma \ref{lem:nujstar} $\nujstar/N \goto \thetajstar$ and ${\nujstarplusone}/{N} \goto \thetajstarplusone$ as $N \goto \infty$. Since 
$m/N \goto 0$ and $\sum_{j \geq \jstarm} \thetaj \goto 0$ as $N \goto \infty$, it follows that
\[
\lim_{N \goto \infty} \left|\sum_{j \in \Nb} f(j) \thetan_j - \sum_{j \in \Nb} f(j) \thetaj\right| = 0 .
\]
This completes the proof of part (a) of Theorem \ref{thm:approximate}. \ \ink

\skp
\noi 
{\bf Proof of part (b) of Theorem \ref{thm:approximate}.}  Let $\theta$ be a probability measure in $\pnbc$. 
We prove that if $R(\theta | \rhoalphab) < \infty$, then 
\[
\lim_{N \goto \infty} R(\thetan | \rhoalphab) = R(\theta | \rhoalphab).
\]

We use the following facts. 
\begin{enumerate}
  \item For all $j \in \Nb$ we have $\lim_{N \goto \infty} \thetanj = \thetaj$.
  \item For all $j \in \nb$ satisfying $j \not = \jstarplusone$, we have $\thetanj \leq \theta$.
\end{enumerate}
Item 1, which is stated in part (c) of Lemma \ref{lem:nujstar}, follows from the weak convergence $\thetan \Rightarrow \theta$ proved in part (a)
of Theorem \ref{thm:approximate}. 
Item 2, which is stated in part (e) of Lemma \ref{lem:nujstar}, is easily verified. For $j=\jstar$ the upper bound  $\thetan_{\jstar} \leq \theta_{\jstar}$ is valid by 
part (a) of Lemma \ref{lem:nujstar}. For all other $j \in \nb$ satisfying $j \not = \jstarplusone$, the upper
bound $\thetanj \leq \thetaj$ is a consequence of (\ref{eqn:nunj}) and (\ref{eqn:nujinequality}). 
According to  part (f) of Lemma \ref{lem:nujstar} the upper bound $\thetan_{\jstarplusone} \leq \theta_{\jstarplusone}$
is not valid for any $N$, and in general there does not exist $M < \infty$ such that for any $N \in \N$,
$\thetan_{\jstarplusone} \leq M \theta_{\jstarplusone}$. Because of this anomaly the term in $R(\thetan | \rhoalphab)$ 
corresponding to $j = \jstarplusone$ must be handled separately.

Define 
$\varphi(x) = x \log x$ for $x \in [0,\infty)$; if $x = 0$, then $\varphi(x) = 0$. This function is continuous on $[0,\infty)$.  For each $j \in \Nb$, since
$\thetanj \goto \thetaj$ as $N \goto \infty$, it follows that $\varphi(\thetanj/\rhoalphabj) \goto \varphi(\thetaj/\rhoalphabj)$ 
as $N \goto \infty$. To prove part (b)
of Theorem \ref{thm:approximate} we must justify the following interchange of the limit $N \goto \infty$ and the sum over $j \in \nbjstar$:
\beas
\lefteqn{
\lim_{N \goto \infty} R(\thetan | \rhoalphab)} \\ 
& & = \lim_{N \goto \infty} \rhoalphabjstarplusone
\varphi(\thetan_{\jstarplusone}/\rhoalphabjstarplusone) 
+ \lim_{N \goto \infty} \sum_{j \in \nbjstar} \rhoalphabj \varphi(\thetanj/\rhoalphabj) 
\\
& & =\rhoalphabjstarplusone \varphi(\theta_{\jstarplusone}/\rhoalphabjstarplusone)  
+ \sum_{j \in \nbjstar} 
\rhoalphabj \left(\lim_{N \goto \infty} \varphi(\thetanj/\rhoalphabj)\right) \\
& & = \rhoalphabjstarplusone \varphi(\theta_{\jstarplusone}/\rhoalphabjstarplusone)  + \sum_{j \in \nbjstar} \rhoalphabj \varphi(\thetaj/\rhoalphabj) = R(\theta | \rhoalphab) .
\eeas
We justify the interchange of the limit and the sum over $j \in \nbjstar$ by applying the Dominated Convergence Theorem. This procedure requires finding constants $a_j$ for $j \in \nbjstar$ such that
for all sufficiently large $N \in \N$ 
\[
\rhoalphabj |\varphi(\thetanj/\rhoalphabj)| \leq a_j \ \mbox{ and } \ \sum_{j \in \nbjstar} a_j < \infty .
\]

The key to applying the Dominated Convergence Theorem is to use two properties of $\varphi(x) = x \log x$: its boundedness on the interval
$[0,1)$ and its monotonicity on the interval $[1,\infty)$.

\skp
\noi 
{\bf Property 1.} For $x \in [0,1)$, $0 \geq \varphi(x) \geq -e^{-1}$.

\noi 
{\bf Property 2.} For $x \in [1,\infty)$, $\varphi(x) \geq 0$, $\varphi(x) \goto \infty$ as $x \goto \infty$, and $\varphi$ is monotone in the sense that
for $1 \leq x < y$, $0 \leq \varphi(x) < \varphi(y)$.

\skp
Let $\Psi = \{\jstarplusone\}$. 
We write $\varphi(x) = \varphiplus(x) - \varphiminus(x)$, where $\varphiplus(x) = \varphi(x) \cdot 1_{[1,\infty)}(x)$ and 
$\varphiminus(x) = - \varphi(x) \cdot 1_{[0,1)}(x)$. For $N \in \N$ define
\[
C_N = \{j \in \Nb \setminus \Psi : \thetanj/\rhoalphabj \in [0,1)\} \ \mbox{ and } \ D_N = \{j \in \Nb \setminus \Psi: \thetanj/\rhoalphabj \in [1,\infty)\} .
\]
In terms of these sets we write
\[
\sum_{j \in \Nb \setminus \Psi} \rhoalphabj |\varphi(\thetanj/\rhoalphabj)| = 
\sum_{j \in C_N} \rhoalphabj \varphiminus(\thetanj/\rhoalphabj) 
+ \sum_{j \in D_N} \rhoalphabj \varphiplus(\thetanj/\rhoalphabj) .
\]
For $j \in C_N$ the boundedness of $\varphi$ on $[0,1)$ implies that
\[
0 \leq \rhoalphabj \varphiminus(\thetanj/\rhoalphabj) \leq e^{-1} \rhoalphabj .
\]
For $j \in D_N$ the monotonicity of $\varphi$ on $[1,\infty)$ and the bound $\thetanj \leq \thetaj$ imply that
\[
0 \leq \rhoalphabj \varphiplus(\thetanj/\rhoalphabj) \leq \rhoalphabj \varphiplus(\thetaj/\rhoalphabj) \leq \rhoalphabj |\varphi(\thetaj/\rhoalphabj)|.
\]
Thus for all $j \in \Nb \setminus \Psi$
\[
\rhoalphabj |\varphi(\thetanj/\rhoalphabj)| \leq a_j = e^{-1} \rhoalphabj + \rhoalphabj |\varphi(\thetaj/\rhoalphabj)|.
\]

Using the fact that $R(\theta | \rhoalphab) < \infty$, we prove that $\sum_{j \in \Nb \setminus \Psi} a_j < \infty$. We have
\bea
\label{eqn:sumaj}
\sum_{j \in \Nb \setminus \Psi} a_j & \leq & e^{-1} \sum_{j \in \Nb \setminus \Psi} \rhoalphabj 
+ \sum_{j \in \Nb \setminus \Psi} \rhoalphabj |\varphi(\thetaj/\rhoalphabj)| \\
\nonumber & \leq & e^{-1} + \sum_{j \in \Nb \setminus \Psi} \rhoalphabj |\varphi(\thetaj/\rhoalphabj)| .
\eea
Define
\[
C = \{j \in \Nb \setminus \Psi: \thetaj/\rhoalphabj \in [0,1)\} \ \mbox{ and } \ D = \{j \in \Nb \setminus \Psi: \thetaj/\rhoalphabj \in [1,\infty)\} .
\]
In terms of these sets we write
\beas
\lefteqn{R(\theta | \rhoalphab)} \\ 
& & = \rhoalphabjstarplusone \varphi(\thetajstarplusone/\rhoalphabjstarplusone) + 
\sum_{j \in \Nb \setminus \Psi} \rhoalphabj \varphi(\thetaj/\rhoalphabj)\\
& & = \rhoalphabjstarplusone \varphi(\thetajstarplusone/\rhoalphabjstarplusone) 
- \sum_{j \in C} \rhoalphabj \varphiminus(\thetaj/\rhoalphabj) + 
\sum_{j \in D} \rhoalphabj \varphiplus(\thetaj/\rhoalphabj) .
\eeas
For $j \in C \cup \Psi$ we have $0 \leq \rhoalphabj \varphiminus(\thetaj/\rhoalphabj) \leq e^{-1} \rhoalphabj$. 
Hence 
\[
\rhoalphabjstarplusone \varphiminus(\thetajstarplusone/\rhoalphabjstarplusone) \leq e^{-1} \rhoalphabjstarplusone 
\leq e^{-1}
\]
and 
\[
\sum_{j \in C} \rhoalphabj \varphiminus(\thetaj/\rhoalphabj) \leq e^{-1} \sum_{j \in C} \rhoalphabj \leq e^{-1}.
\]
It follows that 
\beas
\lefteqn{
\sum_{j \in \Nb \setminus \Psi} \rhoalphabj |\varphi(\thetaj/\rhoalphabj)|} \\
& & = \sum_{j \in C} \rhoalphabj \varphiminus(\thetaj/\rhoalphabj) 
+ \sum_{j \in D} \rhoalphabj \varphiplus(\thetaj/\rhoalphabj) \\
& & \leq e^{-1} + \sum_{j \in D} \rhoalphabj \varphiplus(\thetaj/\rhoalphabj) \\
& & = e^{-1} + R(\theta | \rhoalphab) - \rhoalphabjstarplusone \varphi(\thetajstarplusone/\rhoalphabjstarplusone)  
+ \sum_{j \in C}\rhoalphabj \varphiminus(\thetaj/\rhoalphabj)\\ 
&& \leq e^{-1} + R(\theta | \rhoalphab) + \rhoalphabjstarplusone \varphiminus(\thetajstarplusone/\rhoalphabjstarplusone)  
+ \sum_{j \in C}\rhoalphabj \varphiminus(\thetaj/\rhoalphabj)\\ 
&& \leq 3e^{-1} + R(\theta | \rhoalphab) < \infty .
\eeas
Substituting the last display into (\ref{eqn:sumaj}), we conclude that
\[
\sum_{j \in \Nb \setminus \Psi} a_j < 4e^{-1} + R(\theta|\rhoalphab) < \infty .
\]
This completes the proof of part (b). The proof of Theorem \ref{thm:approximate} is done. \ \ink

\skp
In appendix C we study prove part (a) of Theorem \ref{thm:mainestimate} as well as a number of other properties of the parameter $\alphabc$ that defines the Poisson equilibrium distribution $\rhoalphabc$.


\section{Proof of Part (a) of Theorem \ref{thm:mainestimate} re \boldmath$\alphabc$\unboldmath}
\beginsec

The goal of this appendix is to prove Theorem \ref{thm:alphabc}. 
Part (a) restates part (a) of Theorem \ref{thm:mainestimate} concerning the existence of $\alphabc$.
This parameter defines the Poisson distribution $\rhoalphabc$ appearing in the local large deviation estimate in part (b) of Theorem
\ref{thm:mainestimate}. In part (b) we derive two sets of bounds on $\alphabc$ and use these bounds to show that 
$\alphabc$ is asymptotic to $c$ as $c \goto \infty$. Part (c) 
shows an interesting monotonic relationship between $\alphabc$ and $\alpha_{b+1}(c)$ while part (d) makes precise the relationship between
$\rhoalphabc$ and a Poisson random variable having parameter $\alphabc$. Parts (a), (b), and (d) of the next theorem appear in Theorem C.1 in \cite{EllisTaasan1} as parts (a), (b), and (c). Part (c) of the next theorem is new.

The fact that $\alphabc$ is asymptotic to $c$ as $c \goto \infty$
is certainly plausible. If $c$ is large, then the mean of $\rhoalphabc$, which equals $c$, is not changed appreciably if $\rhoalphabc$ is replaced by
a standard Poisson distribution on $\N \cup \{0\}$
with parameter $\alpha_b(c)$. Since the mean of an actual Poisson distribution on $\N \cup \{0\}$ with parameter $\alpha_b(c)$ is $\alpha_b(c)$, we expect that if $c$ is large, then $\alpha_b(c)$ should be close to $c$. 

\begin{thm}
\label{thm:alphabc}
Fix a nonnegative integer $b$ and a real number $c \in (b,\infty)$. 
For $\alpha \in (0,\infty)$ define $Z_0(\alpha)$, and for $b \in \N$ define $\zalphab = e^\alpha - \sum_{j=0}^{b-1}\alpha^j/j!$. 
Let $\rhoalphab$ be the probability measure on $\Nb$ whose components are defined by 
\[
\rhoalphabj = \frac{1}{\zalphab} \cdot \frac{\alpha^j}{j!} \ \mbox{ for } j \in \Nb .
\]
The following conclusions hold.

{\em (a)} There exists a unique value $\alphabc \in (0,\infty)$ such that 
$\rhoalphabc$ lies in the set $\pnbc$ of probability measures on $\Nb$ having mean $c$. If $b=0$, then 
$\alpha_0(c) = c$. If $b \in \N$, then $\alphabc$ is the unique solution in $(0,\infty)$ of $\alpha \zalphabminusone/\zalphab = c$.
 
{\em (b)} For $b \in \N$ 
\[
c > \alphabc > c-b \ \mbox{ and } c > \alphabc > c(1 - 2^b e^{-(c-b)/2}).
\]
Either of these bounds imply that $\alphabc$ is asymptotic to $c$ as $c \goto \infty$; i.e., $\lim_{c \goto \infty} \alphabc/c = 1$.

{\em (c)} For all $b \in \N \cup \{0\}$ and $c > b+1$,
$\alpha_{b+1}(c) < \alphabc$.

{\em (d)} For $b \in \N$, if $\Xi_{\alphabc}$ is a Poisson random variable with parameter $\alphabc$,
then $\rhoalphabc$ is the distribution of $\Xi_{\alphabc}$ conditioned on $\Xi_{\alphabc} \in \Nb$.
\end{thm}

Before we prove Theorem \ref{thm:alphabc}, we state a second theorem that focuses on the case $b=1$. In this 
case the equilibrium distribution $\rhoalphaonec$ is a probability measure on $\N_1 =\N$. In part (a) 
we give the proof of the existence of $\alphabc$ for $b=1$, which is much more straightforward than the proof for general $b$. 
In parts (b) and (c) we give two iterative procedures for calculating $\alphaonec$ while in part (d) we derive two sets of 
inequalities that are tighter than the inequalities for $\alphabc$ for general $b$ given in part (b) of Theorem \ref{thm:alphabc}.
Like the inequalities in part (b) of Theorem \ref{thm:alphabc}, the inequalities in part (d) of the next theorem imply that 
$\alphaonec$ is asymptotic to $c$ as $c \goto \infty$. 

\begin{thm}
\label{thm:alphac}
Fix a real number $c \in (1,\infty)$. The following results are valid. 

{\em (a)} There exists a unique value $\alphaonec \in (0,\infty)$ such that $\rho_{1,\alphaonec}$ lies in the set $\pnc$ of probability
measures on $\N$. The quantity $\alphaonec$ is the unique solution in $(0,\infty)$ of $\alpha e^\alpha = c(e^{\alpha} - 1)$. 

{\em (b)} Let $\alpha_1 = c$ and consider the following iterative procedure defined for $n \in \N, n \geq 2$:
\[
\alpha_{n+1} = c(1 - e^{-\alpha_{n}}) .
\]
Then the sequence $\{\alpha_n, n \in \N\}$ is monotonically decreasing and $\lim_{n \goto \infty} \alpha_n = \alphaonec$.

{\em (c)} Let $\beta_1 = \log c$ and consider the following iterative procedure defined for $n \in \N, n \geq 2$:
\[
\beta_{n+1} = c(1 - e^{-\beta_{n}}) .
\]
Then the sequence $\{\beta_n, n \in \N\}$ is monotonically increasing and $\lim_{n \goto \infty} \beta_n = \alphaonec$.

{\em (d)} We have the following two bounds on $\alphaonec$:
\[ 
c(1 - e^{-c}) > \alphaonec > c-1 \ \mbox{ and } \ c(1 - e^{-c}) > \alphaonec > c(1 - e^{-c + 1}) .
\]
Either of these bounds implies that $\alphaonec$ is asymptotic to $c$ as $c \goto \infty$; i.e., $\lim_{c \goto \infty} \alphaonec/c = 1$. 
\end{thm}

\noindent 
{\bf Proof.} (a) The measure $\rho_{1,\alpha}$ is a probability measure on $\N$ having mean 
\beas
\sum_{j \in \N} j \rho_{1,\alpha;j} & = & \frac{1}{e^{\alpha -1}} \cdot \sum_{j \in \N} \frac{\alpha^j}{(j-1)!} \\
& = & \frac{1}{e^{\alpha} -1} \cdot 
\alpha \sum_{j=0}^\infty \frac{\alpha^j}{j!} = \frac{1}{e^{\alpha} -1} \cdot \alpha e^{\alpha} .
\eeas
Thus $\rho_{1,\alpha}$ has mean $c$ if and only if $\alpha$ satisfies $\alpha e^{\alpha} = c(e^{\alpha} - 1)$.
We prove part (a) by showing that this equation has a unique solution $\alphaonec \in (0,\infty)$ for any $c > 1$.

The proof that $\alpha e^{\alpha} = c(e^{\alpha} - 1)$ has a unique 
solution $\alphaonec \in (0,\infty)$ for any $c > 1$ is straightforward. A positive real number $\alpha$ solves $\alpha e^\alpha = c(e^\alpha - 1)$ if and only if
\[
\gamma_1(\alpha) = c, \mbox{ where } \gamma_1(\alpha) = \frac{\alpha}{1 - e^{-\alpha}} . 
\]
The function $\gamma_1$ is continuously differentiable on $(0,\infty)$ and $\lim_{\alpha \goto 0^+} \gamma_1(\alpha) = 1$. In addition, for $\alpha \in (0,\infty)$
\[
\gamma_1'(\alpha) = \frac{1 - (1 + \alpha)e^{-\alpha}}{(1 - e^{-\alpha})^2} = e^{-\alpha} \cdot \frac{e^\alpha - 1 - \alpha}{(1 - e^{-\alpha})^2} > 0.
\]
The inequality holds since for $\alpha > 0$, $e^\alpha - 1 - \alpha > 0$. It follows that there exists a sufficiently small value of $\ve > 0$ such that
$1 < \gamma_1(\ve) < c$ and $\gamma_1$ is monotonically increasing on $(\ve,\infty)$. Since $\gamma_1(\alpha) \goto \infty$ as $\alpha \goto \infty$, we conclude that
there exists a unique value $\alpha = \alphaonec \in (0,\infty)$ solving $\gamma_1(\alphaonec) = c$ and thus solving $\alphaonec e^{\alphaonec} = c(e^{\alphaonec} - 1)$. 
This completes the proof of part (a).

(b) Since $e^{-c} < 1$, we have the inequality 
\[
\alpha_2 = c(1 - e^{-\alpha_{1}}) = c(1 - e^{-c}) < c = \alpha_1 .
\]
We use induction to prove that the sequence $\alpha_n$ is monotonically decreasing. For $n \in \N, n \geq 2$, under the assumption
that $\alpha_n < \alpha_{n-1}$, this property of the sequence is a consequence of the following calculation: 
\[
\alpha_{n+1} - \alpha_n = c(e^{-\alpha_{n-1}} - e^{-\alpha_n}) < 0 .
\]
We now use induction to prove that the sequence $\alpha_n$ is bounded below by $\log c$. For $n=1$, $a_1 = c > \log c$. Assuming that
$\alpha_n > \log c$, we have
\[
\alpha_{n+1} = c(1 - e^{-\alpha_n}) > c(1 - e^{-\log c}) = c-1 > \log c .
\]
The last inequality follows from the facts that when $c=1$, $c - 1 = 0 = \log c$ and that for $c \in (1,\infty)$, $(c-1)' = 1 > 1/c = (\log c)'$.
This completes the proof that $\alpha_n > \log c$ for all $n \in \N$. Since $\alpha_n$ is a monotonically decreasing sequence bounded 
above by $c$ and below by $\log c$, we conclude 
$\alpha^\star = \lim_{n \goto \infty} \alpha_n$ exists and satisfies both $\alpha^\star \in (\log c, c)$ and 
$\alpha^\star = c(1 - e^{-\alpha^\star})$. Because $\alphaonec$ is the unique positive solution of this equation,
it follows that $\lim_{n \goto \infty} \alpha_n = \alphaonec$. This completes the proof of part (b). 

(c) Since $\beta_1 = \log c$, we have the inequality 
\[
\beta_2 = c(1 - e^{-\beta_{1}}) = c(1 - e^{-\log c}) = c - 1 > \log c.
\]
We use induction to prove that the sequence $\beta_n$ is monotonically increasing. For $n \in \N, n \geq 2$, under the assumption
that $\beta_{n-1} < \beta_{n}$, this is a consequence of the following calculation: 
\[
\beta_{n+1} - \beta_n = c(e^{-\beta_{n-1}} - e^{-\beta_n}) > 0 .
\]
We now use induction to prove that the sequence $\beta_n$ is bounded above by $c$. 
For $n = 1$, $\beta_1 = \log c < c$. Assuming that $\beta_n < c$, we have
\[
\beta_{n+1} = c(1 - e^{-\beta_n}) < c(1 - e^{-c}) < c ,
\]
This completes the proof that $\beta_n$ is bounded above by $c$. 
Since $\beta_n$ is a monotonically increasing sequence bounded 
above by $c$ and below by $\log c$, we conclude 
$\beta^\star = \lim_{n \goto \infty} \beta_n$ exists and satisfies both $\beta^\star \in (\log c, c)$ and 
$\beta^\star = c(1 - e^{-\beta^\star})$. Because $\alphaonec$ is the unique positive solution of this equation,
it follows that $\lim_{n \goto \infty} \beta_n = \alpha(c)$. This completes the proof of part (c). 

(d) We first prove that $c(1 - e^{-c}) > \alphaonec$. This follows immediately from the iterative procedure discussed in part (a),
which implies that $c = \alpha_1 > \alpha_2 = c(1 - e^{-c}) > \alphaonec$. One can obtain the weaker upper bound $c > \alphaonec$
directly if one writes the equation solved by $\alphaonec$ in the form
\be 
\label{eqn:solved}
\alphaonec = c (1 - e^{-\alphaonec}) 
\ee
and uses the fact that $e^{-\alphaonec} \in (0,1)$.

We now prove a series of three lower bounds,
the last two of which, in combination with the upper bound $c(1 - e^{-c}) > \alphaonec$, imply that $\alphaonec \sim c$
as $c \goto \infty$. The first lower bound is $\alphaonec > \log c$. 
To prove this, we use the fact that $\alphaonec > 0$ to write $e^{\alphaonec} - 1 \geq \alphaonec$. 
It follows that 
\[
\alphaonec = c (1 - e^{-\alphaonec}) = c e^{-\alphaonec} (e^{\alphaonec} - 1) > c e^{-\alphaonec} \alphaonec ,
\]
or equivalently that $e^{\alphaonec} > c$. This 
implies that $\alphaonec > \log c$, as claimed. 

We now bootstrap this lower bound into a tighter lower bound by substituting 
$\alphaonec > \log c$ into the right hand side of (\ref{eqn:solved}), obtaining the second lower bound
\be 
\label{eqn:asympalphac1}
\alpha_1(c) = c (1 - e^{-\alpha_1(c)}) > c(1 - e^{-\log c}) = c\left(1 - \frac{1}{c}\right) = c-1 .
\ee
It follows that 
\[
1 - e^{-c} > \frac{\alphaonec}{c} > 1 - \frac{1}{c} .
\]
This implies that $\lim_{c \goto \infty} \alpha_1(c)/c = 1$ or that $\alphaonec$ is asymptotic to $c$ as $c \goto \infty$.

By bootstrapping the lower bound in (\ref{eqn:asympalphac1}), we obtain yet a tighter lower bound on $\alphaonec$ which gives a second
proof that $\alphaonec \sim c$.  To do this, 
we substitute $\alpha_1(c) > c-1$ into the right hand side of (\ref{eqn:solved}), obtaining the third
lower bound $\alpha_1(c) > c(1 - e^{-c +1})$. It follows that
\be 
\label{eqn:asympalphac}
1 - e^{-c} > \frac{\alpha_1(c)}{c} > 1 - e^{-c +1} .
\ee
This implies $\lim_{c \goto \infty} \alphaonec/c = 1$ at a rate that is at least exponentially fast.
By contrast, (\ref{eqn:asympalphac1}) shows a much slower rate of convergence to 1 that is only of the order $1/c$. Interestingly,
iterating this procedure
again does not give a tighter lower bound than that in (\ref{eqn:asympalphac}). This completes the proof of Theorem \ref{thm:alphac}. \ \ink

\skp
We now turn to the proof of Theorem \ref{thm:alphabc}.
According to part (a) of this theorem, for $b \in \N$, $\alphabc$ is the unique solution of $\alpha \zalphabminusone/\zalphab = c$.
The heart of the proof of Theorem \ref{thm:alphabc}, and its most subtle step, is to prove that the function 
$\gamma_b(\alpha) = \alpha \zalphabminusone/\zalphab$ satisfies $\gamma_b'(\alpha) > 0$ for $\alpha \in (0,\infty)$ and thus
is monotonically increasing on this interval. This fact is proved in the next lemma.

\begin{lem}
\label{lem:alphabc}
Fix a positive integer $b$ and a real number $c \in (b,\infty)$. For $\alpha \in (0,\infty)$ the function 
$\gamma_b(\alpha) = \alpha \zalphabminusone/\zalphab$ satisfies $\gamma_b'(\alpha) > 0$.
\end{lem}

{\noindent}
{\bf Proof.} For $b \in \N$ and for $\alpha \in (0,\infty)$, we have $Z_{b}'(\alpha) = \zbminusonealpha$. Thus
\[
\gamma_b(\alpha) = \frac{\alpha \zbminusonealpha}{\zbalpha} = \alpha (\log \zbalpha)'.
\]
The key to proving that $\gamma_b'(\alpha) > 0$ is to represent $\log \zbalpha$ in terms of the moment generating function of a probability measure.
We do this by first expressing $\zbalpha$ in terms of the upper incomplete gamma function via the formula
\be 
\label{eqn:gammafunction}
\zbalpha = \frac{e^\alpha}{(b-1)!} \int_0^\alpha x^{b-1} e^{-x} dx.
\ee
This formula is easily proved by induction. For $b=1$ the right side equals $e^\alpha - 1 = Z_1(\alpha)$. Assuming that it is true for $b=n$, we prove that it
is true for $b=n+1$ by integrating by parts, which gives
\beas
\frac{e^\alpha}{n!} \int_0^\alpha x^{n} e^{-x} dx & = & \frac{e^\alpha}{(n-1)!} \int_0^\alpha x^{n-1} e^{-x} dx - \frac{\alpha^{n}}{n!} \\
& = & Z_{n}(\alpha) - \frac{\alpha^{n}}{n!} = Z_{n+1}(\alpha).
\eeas
This completes the proof of (\ref{eqn:gammafunction}) for all $b \in \N$.

As suggested in \cite{Neuman}, we now make the change of variables $x = y \alpha$, obtaining the representation
\be 
\label{eqn:changevariables}
\zbalpha = \frac{e^\alpha}{b!} \alpha^b g_b(\alpha), \mbox{ where } g_b(\alpha) = \int_{-1}^0 e^{\alpha y} b(-y)^{b-1} dy.
\ee
The function $g_b$ is the moment generating function of the probability measure on $\R$ having the density $h_b(y) = b(-y)^{b-1}$ on 
$[-1,0]$. For $\alpha \in (0,\infty)$ let $\sigma_{b,\alpha}$ be the probability measure on $\R$ having the density $e^{\alpha y} h_b(y)/g_b(\alpha)$ on $[-1,0]$.
A straightforward calculation shows that 
\[
(\log g_b)'(\alpha) = \int_{\R} y \sigma_{b,\alpha}(dy) \ \mbox{ and } \ 
(\log g_b)''(\alpha) = \int_{\R} [y - g_b'(\alpha)]^2 \sigma_{b,\alpha}(dy).
\]
As the variance of the nontrivial probability measure $\sigma_{b,\alpha}$, we conclude that $(\log g_b)''(\alpha) > 0$ for all $\alpha \in (0,\infty)$.

Using (\ref{eqn:changevariables}) and the power series representations
\[
\zbminusonealpha = \sum_{j=b-1}^\infty \frac{\alpha^j}{j!} \ \mbox{ and } \zbalpha = \sum_{j=b}^\infty \frac{\alpha^j}{j!},
\] 
we calculate
\beas
\gamma_b'(\alpha) & = & (\log \zbalpha)' + \alpha (\log \zbalpha)'' \\
& = & (\log \zbalpha)' + \alpha \left[\log\left(\frac{e^\alpha}{b!} \alpha^b g_b(\alpha)\right)\right]'' \\
& = & \frac{\zbminusonealpha}{\zbalpha} + \alpha [\alpha  - \log(b!) + b \log \alpha + \log g_b(\alpha)]'' \\
& = &  \frac{\zbminusonealpha}{\zbalpha} - \frac{b}{\alpha} + \alpha (\log g_b(\alpha))'' \\
& = & \frac{\alpha \zbminusonealpha - b \zbalpha}{\alpha\zbalpha} + \alpha (\log g_b(\alpha))'' \\
& = & \frac{1}{\alpha \zbalpha} \cdot \sum_{j=b}^\infty \left(\frac{1}{(j-1)!} - \frac{b}{j!}\right)\alpha^j + \alpha (\log g_b(\alpha))'' \\
& = & \frac{1}{\zbalpha} \cdot \sum_{j=b}^\infty \frac{j-b}{j!}\alpha^{j-1} + \alpha (\log g_b(\alpha))'' > 0.
\eeas
This completes the proof of the lemma. \ink

\skp
We are now ready to prove Theorem \ref{thm:alphabc}. 

\skp
\noindent
{\bf Proof of Theorem \ref{thm:alphabc}.} (a) We first consider $b=0$. In this case $\rho_{0,\alpha}$ is a standard Poisson distribution
on $\N_0$ having mean $\alpha$. It follows that $\alpha_0(c) = c$ is the unique value for which $\rho_{0,\alpha_0(c)}$ has mean $c$ and thus lies 
in ${\mathcal P}_{\N_0,c}$.
This completes the proof of part (a) for $b = 0$.

We now consider $b \in \N$.  In this case $\rhoalphab$ is a probability measure on $\Nb$ having mean 
\bea
\label{eqn:meanrho}
\sum_{j \in \Nb} j \rhoalphabj & = & \frac{1}{\zalphab} \cdot \sum_{j \in \Nb} \frac{\alpha^j}{(j-1)!} \\ 
\nonumber
& = & \frac{1}{\zalphab} \cdot 
\alpha \sum_{j=b-1}^\infty \frac{\alpha^j}{j!} = \frac{1}{\zalphab} \cdot \alpha Z_{b-1}(\alpha).
\eea
Thus $\rhoalphab$ has mean $c$ if and only if $\alpha$ satisfies $\gamma_b(\alpha) = c$,
where $\gamma_b(\alpha) = \alpha Z_{b-1}(\alpha)/\zalphab$. We prove part (a) by showing 
that $\gamma_b(\alpha) = c$ has a unique solution $\alphabc \in (0,\infty)$ for all $b \in \N$ and any $c > b$.

The proof depends on the following three steps: 
\begin{enumerate}
  \item $\lim_{\alpha \goto 0^+}\gamma_b(\alpha) = b$;
  \item $\lim_{\alpha \goto \infty} \gammabalpha = \infty$;
  \item for all $\alpha \in (0,\infty)$, $\gamma_b'(\alpha) > 0$.
\end{enumerate} 
These three steps yield part (a). Indeed, by steps 1 and 3 there exists a sufficiently small value of $\ve > 0$ such that
$b < \gamma_b(\ve) < c$, and by step 3 $\gamma_b$ is monotonically increasing on $(\ve,\infty)$. 
Since by step 2 $\gamma_b(\alpha) \goto \infty$ as $\alpha \goto \infty$, we conclude that
there exists a unique value $\alpha = \alphabc \in (0,\infty)$ solving $\gamma_b(\alphabc) = c$ and thus guaranteeing that
$\rhoalphabc \in \pnbc$.

Step 3 is proved in Lemma \ref{lem:alphabc}. We now prove steps 1 and 2.

\skp
\noi 
{\bf Step 1.} For $b \in \N$ and for $\alpha \in (0,\infty)$ satisfying $\alpha \goto 0$
\beas
\gammabalpha & = & \frac{\alpha \zalphabminusone}{\zalphab} = \frac{\alpha \sum_{j=b-1}^\infty \alpha^j/j!}{\sum_{j=b}^\infty \alpha^j/j!} \\
& = & \frac{\sum_{j=b}^\infty \alpha^j/(j-1)!}{\sum_{j=b}^\infty \alpha^j/j!} = \frac{\alpha^b/(b-1)! + \mbox{o}(1)}{\alpha^b/b! + \mbox{o}(1)}
= b + \mbox{o}(1).
\eeas
The terms denoted by $\mbox{o}(1)$ converge to 0 as $\alpha \goto 0$. It follows that $\lim_{\alpha \goto 0^+}\gamma_b(\alpha) = b$. This completes the proof 
of step 1.

\skp
\noi 
{\bf Step 2.} For $b \in \N$ and for $\alpha \in (0,\infty)$
\beas
\gammabalpha & = & \frac{\alpha \zalphabminusone}{\zalphab} = \frac{\alpha \left(e^\alpha - \sum_{j=0}^{b-2} \alpha^j/j!\right)}
{e^\alpha - \sum_{j=0}^{b-1} \alpha^j/j!} \\ 
& = & \frac{\alpha \left(1 - e^{-\alpha}\sum_{j=0}^{b-2} \alpha^j/j!\right)}
{1 - e^{-\alpha} \sum_{j=0}^{b-1} \alpha^j/j!} = \alpha(1 + \mbox{o}(1)).
\eeas
The term denoted by $\mbox{o}(1)$ converges to 0 as $\alpha \goto \infty$. It follows that $\lim_{\alpha \goto \infty} \gammabalpha = \infty$.
This completes the proof of step 2.

Having completed steps 1, 2, and 3, we have proved part (a) for all $b \in \N$. Since we also validated
part (a) for $b=0$, the proof of part (a) for all nonnegative integers $b$ is done.

\skp
\noi 
(b) We first prove that $\alphabc < c$ 
for $b \in \N$ by observing that for any $\alpha \in (0,\infty)$ we have $\zbminusonealpha > \zalphab$.
Thus $\gamma_b(\alpha) = \alpha\zbminusonealpha/\zbalpha > \alpha$, 
which implies that $\alphabc < \gamma_b(\alphabc) = c$. To prove that $\alphabc > c-b$, we use the inequality 
\[
\zbalpha = \sum_{j=b}^\infty \frac{\alpha^j}{j!} > \frac{\alpha^b}{b!}
\]
to write 
\beas
\gammabalpha & = & \frac{\alpha \zbminusonealpha}{\zalphab} = \alpha + \frac{\alpha (\zalphabminusone - \zalphab)}{\zalphab} \\
& = & \alpha + \frac{\alpha^b/(b-1)!}{\zalphab} < \alpha + \frac{\alpha^b/(b-1)!}{\alpha^b/b!} = \alpha + b.
\eeas
It follows that $c = \gamma_b(\alphabc) < \alphabc + b$, which gives the desired lower bound $\alphabc > c-b$.

We now bootstrap this lower bound into the tighter lower bound indicated in part (b). To do this we note that for
any $\alpha \in (0,\infty)$ 
\beas
\alpha e^\alpha & > & \alpha \zbminusonealpha = \gammabalpha \zbalpha \\
& = & \gammabalpha \left(e^\alpha - \sum_{j=0}^{b-1} \frac{\alpha^j}{j!}\right) > \gammabalpha (e^\alpha - 2^b e^{\alpha/2}).
\eeas
The first lower bound $\alphabc > c-b$ now yields the tighter lower bound
\[
\alphabc > \gamma_b(\alphabc) (1 - 2^b e^{-\alphabc/2}) = c(1 - 2^b e^{-\alphabc/2}) > c(1-2^b e^{-(c-b)/2}).
\]
This completes the proof of the bounds in part (b).
Either of these bounds imply that $\lim_{c \goto \infty} \alphabc/c = 1$. This proves that $\alphabc$ is asymptotic to $c$
as $c \goto \infty$, completing the proof of part (b).

\skp
\noi 
(c) According to part (b), for $c > 1$ we have 
$\alpha_1(c) < c = \alpha_0(c)$. In order to prove that for $b \in \N$ and $c > b+1$ we have
$\alpha_{b+1}(c) < \alphabc$ , we first prove that for $b \in \N$ and any $\alpha \in (0,\infty)$
we have $\gamma_{b}(\alpha) < \gamma_{b+1}(\alpha)$. As shown in the proof of part (b), for all $\alpha \in (0,\infty)$
\[
\gammabalpha = \alpha + \frac{\alpha^b/(b-1)!}{\zalphab} = \alpha + \frac{\alpha^b}{(b-1)! \cdot \zalphab}.
\]
By substituting the power series representation for $\zalphab$, we find that
\[
\frac{(b-1)! \cdot \zalphab}{\alpha^b} = (b-1)! \cdot \sum_{j=0}^\infty \frac{\alpha^j}{(j+b)!} =
\sum_{j=0}^\infty \frac{\alpha^j}{\prod_{i=0}^j (b+i)}.
\] 
Since the product $\prod_{i=0}^j (b+i)$ is a strictly increasing function of $b \in \N$, it follows that for fixed $\alpha \in (0,\infty)$
\[
\gammabalpha = \alpha + \frac{\alpha^b}{(b-1)! \cdot \zalphab} =
\alpha + \left(\sum_{j=0}^\infty \frac{\alpha^j}{\prod_{i=0}^j (b+i)}\right)^{-1}
\]
is a strictly increasing function of $b \in \N$. This proves that $\gamma_{b}(\alpha) < \gamma_{b+1}(\alpha)$ for $b \in \N$. 
We now choose $c > b+1$. Then $c = \gamma_b(\alphabc) < \gamma_{b+1}(\alphabc)$. In step 3 in the proof of part (a) we showed that 
 $\gamma_b'(\alpha) > 0$ for $\alpha \in (0,\infty)$ and thus that $\gamma_b$ is strictly increasing on $(0,\infty)$. If $\alpha_{b+1}(c) \geq \alphabc$, it would then 
follow that $c < \gamma_{b+1}(\alphabc) \leq \gamma_{b+1}(\alpha_{b+1}(c))$. This contradicts
the fact that $\gamma_{b+1}(\alpha_{b+1}(c)) = c$ and completes the proof of assertion (c).

\skp
\noi 
(d) For $b \in \N$ we identify $\rhoalphabc$ as the distribution of $\Xi_{\alphabc}$ conditioned on $\Xi_{\alphabc} \in \Nb$. Let $\Xi_{\alphabc}$ be defined on a probability
space having measure $P$. For any $j \in \Nb$
\beas
P(\Xi_{\alphabc} = j \, | \, \Xi_{\alphabc} \in \Nb) & = & \frac{1}{P(\Xi_{\alphabc} \in \Nb)} \cdot P(\Xi_{\alphabc} = j) \\ 
& = & \frac{1}{1 -  e^{\alphabc} \textstyle \sum_{i=0}^{b-1} [\alphabc]^i/i!} \cdot e^{-\alphabc} \frac{[\alphabc]^j}{j!} \\
& = & \frac{1}{\zalphabc} \cdot \frac{[\alphabc]^j}{j!} = \rhoalphabcj.
\eeas
This completes the proof of part (d). The proof of Theorem \ref{thm:alphabc} is done as is the proof of part (a) of Theorem
\ref{thm:mainestimate}. \ink

\skp
In the next and final appendix we explore how the restriction involving $m = m(N)$ could be avoided in the definition of the set of configurations
$\omeganbm$ in (\ref{eqn:omegankm}) and in the definition of the microcanonical ensemble $\pnbm$ in (\ref{eqn:condprob}). Avoiding this restriction
would enable us to present our results in a more natural form.

\section{Avoiding Restriction Involving \boldmath$m = m(N)$\unboldmath}
\beginsec

In this appendix we explore a more natural formulation of our results, and we explain the issues that make such a formulation so challenging.
Among these issues there is a limitation that seems to be inherent in the approximation procedure we use to prove our results. This discussion
makes contact with several interesting ideas including Stirling numbers of the second kind and associated Stirling numbers of the second kind.

Let us review the notation. 
We start with the configuration space $\Omega_N = \Lambda_N^K$. 
For $\omega \in \Omega_N$, $K_\ell(\omega)$ is the droplet-size random variable denoting the number of particles occupying the site 
$\ell \in \Lambda_N$, and $N_j(\omega)$ is
the number of sites for which $K_\ell(\omega) = j$. We also introduce $|N(\omega)|_+$, which is the number of indices $j$ for which $N_j(\omega) \geq 1$. 
Given $b$ a nonnegative integer, we focus on the configuration space $\omeganbm$ consisting of all 
$\omega \in \Omega_N$ for which every site of $\Lambda_N$ is occupied by at least $b$ particles 
and for which $|N(\omega)|_+ \leq m$. The quantity $m$ is a function $m(N)$ satisfying $m(N) \goto \infty$ and $m(N)^2/N \goto 0$ as $N \goto 
\infty$. In symbols
\be 
\label{eqn:needomeganbm}
\omeganbm = \{\omega \in \Omega_N : K_\ell(\omega) \geq b \ \forall \ell \in \Lambda_N \mbox{ and } |N(\omega)|_+ \leq m\}.
\ee

The first constraint involving $K_\ell$ is intrinsic to the definition of the model. By contrast, the second constraint involving 
$m$ is not intrinsic to the definition of the model, but rather is a useful technical device that enables us to control the errors that arise at various stages of the analysis. 
A more natural configuration space would be the set $\omeganb$ consisting of all $\omega \in \Omega_N$ 
for which every site of $\Lambda_N$ is occupied by at least $b$ particles but for which there is no restriction on the 
number of positive quantities $N_j(\omega)$. In symbols
\be 
\label{eqn:needomeganb}
\omeganb = \{\omega \in \Omega_N : K_\ell(\omega) \geq b \ \forall \ell \in \Lambda_N\}.
\ee

We now come to the main point. Let $P_N$ be the uniform probability measure on $\Omega_N$ that assigns equal probability $1/N^K$ to each of the $N^K$ configurations in $\Omega_N$. 
All of the results in the paper are formulated for the probability measure $\Pnbm$, defined as the restriction of $P_N$ to $\omeganbm$. 
However, because the second constraint in the definition
of $\omeganbm$ involving $m$ is not intrinsic to the definition of the model, it would be more natural to formulate 
our results for the probability measure $\Pnb$, defined as
the restriction of $P_N$ to the larger and more natural configuration space $\omeganb$.

In order to understand why our results are formulated for $\Pnbm$ and not for $\Pnb$, we explain how the constraint involving $m$ arises in the paper. There are three sources. First, in Lemma \ref{lem:deltankmnu} we require that $m \log N/N \goto 0$ as $N \goto \infty$ to prove that the error
$\zeta^{(2)}_N(\nu)$ in (\ref{eqn:firststirling}) converges to 0 uniformly for $\nu \in \anbm$. 
Second, we require that $m/N 
\goto 0$ as $N \goto \infty$ to prove part (a) of Lemma \ref{lem:omegankm} and the weak convergence
$\thetan \Rightarrow \theta$ in part (a) of Theorem \ref{thm:approximate}.
Part (a) of Lemma \ref{lem:omegankm} is used to prove part (b) of the lemma and to 
verify hypothesis (i) in Theorem \ref{thm:balls} when applied to Theorem \ref{thm:ldlimitballs}. Third, to prove part (b) of Lemma \ref{lem:omegankm} and to 
verify hypothesis (iv) in Theorem \ref{thm:balls} when applied to Theorem \ref{thm:ldlimitballs}, the stronger condition that $m^2/N \goto 0$ 
as $N \goto \infty$ is required. The source of this error is Lemma \ref{lem:nujstar}, which is used to prove the approximation result in Theorem 
\ref{thm:approximate}. This stronger condition on $m$ is optimal in the sense that it is a minimal assumption guaranteeing that an error term in the lower bound in part (a) of Lemma \ref{lem:nujstar} and in the upper bound in part (b) of the lemma converge to 0.

The stronger condition that $m^2/N \goto 0$ as $N \goto \infty$ means that $m \goto \infty$ at a slower rate than $\sqrt{N}$. What we find fascinating is the fact that the relationship between $m$ and $\sqrt{N}$ is also central to another component of our analysis. As we show in the next theorem, 
if $m \goto \infty$ at a faster rate than $\sqrt{N}$, then for all sufficiently large $N$ the configuration spaces $\omeganbm$ and $\omeganb$ coincide as do the conditional probability measures $\Pnbm$ and $\Pnb$. 

\begin{thm}
\label{thm:msqrtn}
Fix a nonnegative integer $b$ and a rational number $c \in (b,\infty)$. 
Define $\omeganb$ as in {\em (\ref{eqn:needomeganb})} and $\omeganbm$ as in {\em (\ref{eqn:needomeganbm})}, 
where $m = m(N)$ any function satisfying $m(N) \goto \infty$ as $N \goto \infty$.
The following conclusions hold.

{\em (a)} $\max_{\omega \in \omeganb} |N(\omega)|_+ = \sqrt{2(cN + 1/8)} - 1/2$.

{\em (b)} If $m/\sqrt{N} \goto \infty$ as $N \goto \infty$, then for all sufficiently large $N$, $\omeganbm = \omeganb$ and $\Pnbm = \Pnb$.
\end{thm} 

\noi 
{\bf Proof.} (a) For $\omega \in \omeganb$, $N(\omega)$ denotes the sequence $\{N_j(\omega), j \in \Nb\}$. Let $\kappa(\omega) = |N(\omega)|_+$,
and let $1 \leq j_1 < j_2 < \ldots < j_{\kappa(\omega)}$ denote 
the indices for which $N_j(\omega) \geq 1$. We have strict inequality since the $|N(\omega)|_+$ droplet classes have different
sizes. Since for each of these indices we have $j_k \geq k$, the second conservation law in (\ref{eqn:conserve}) 
implies that
\[
K = cN = \sum_{k=1}^{\kappa(\omega)} j_k N_{j_k}(\omega) \geq \sum_{k=1}^{\kappa(\omega)} k N_{j_k}(\omega) 
\geq \sum_{k=1}^{\kappa(\omega)} k = \frac{\kappa(\omega)(\kappa(\omega) + 1)}{2}.
\]
It follows that
\[
2cN \geq \kappa(\omega)(\kappa(\omega) + 1) = (\kappa(\omega) + 1/2)^2 - 1/4,
\]
which in turn implies that 
\be 
\label{eqn:kappa}
\kappa(\omega) = |N(\omega)|_+ \leq \sqrt{2(cN + 1/8)} - 1/2.
\ee
Now let $\omega$ be any configuration in $\omeganb$
for which $N_k(\omega) = 1$ for $k = 1,2,\ldots,|N(\omega)|_+$. In this case
\[
K = cN = \sum_{k=1}^{\kappa(\omega)} k N_{k}(\omega) = \frac{\kappa(\omega)(\kappa(\omega) + 1)}{2},
\]
which in turn implies that $|N(\omega)|_+ = \sqrt{2(cN + 1/8)} - 1/2$. Since this gives equality in (\ref{eqn:kappa}), the proof of part (a) is complete.

(b) Since $m/\sqrt{N} \goto \infty$, part (a) implies that for any $\omega \in \omeganbm$ we have $|N(\omega)|_+ \leq m$ for all sufficiently
large $N$. It follows that 
for all sufficiently large $N$, $\omeganbm = \omeganb$. Since $\Pnbm$ and $\Pnb$ are the respective restrictions of $P_N$ to $\omeganbm$ and $\omeganb$, it also
follows that these two probability measures coincide for all sufficiently large $N$. The proof of the lemma is complete. \ink

\skp

Theorem D.1 motivated us to seek a new approximation procedure. The new procedure would replace the condition $m^2/N \goto 0$, needed to prove Lemma \ref{lem:nujstar}, with a function $m = m(N)$ satisfying 
$m/\sqrt{N} \goto \infty$, needed to prove Theorem D.1,
and satisfying the conditions needed to prove Lemma \ref{lem:deltankmnu}, part (a)
of Lemma \ref{lem:omegankm}, and part (a) of Theorem \ref{thm:approximate}, 
which are $m \log N/N \goto 0$ and $m/N \goto 0$; an example of such a function would be
$m = N^\delta$ for some $\delta \in (1/2,1)$. If we could find such an approximation
procedure,
then all our results formulated for $\Pnbm$ would automatically hold for the more 
natural measure $\Pnb$. Unfortunately, despite great effort, we were unsuccessful. 

Because of this situation it is worthwhile to look more closely at the two components of the approximation procedure presented in appendix B.
Given any measure $\theta \in \pnbc$, this procedure constructs a sequence $\theta^{(N)}$ lying in the range
$\bnbm$ of $\Thetanb$ and having the following
two properties: 
\begin{itemize}
  \item[(a)] $\theta^{(N)} \Rightarrow \theta$ as $N \goto \infty$;
  \item[(b)] if $R(\theta|\rhoalphab) < \infty$, then $R(\theta^{(N)}|\rhoalphab) \goto R(\theta|\rhoalphab)$ as $N \goto \infty$. 
\end{itemize} 
We are able to construct a number of sequences $\theta^{(N)} \in \bnbm$ that satisfy property (a) 
under the hypothesis that $m/N \goto 0$. However, none of these satisfy property (b) 
with a function $m$ satisfying $m/\sqrt{N} \goto \infty$. On the basis of this experience,
we conjecture that there exists no sequence $\theta^{(N)} \in \bnbm$ satisfying both properties (a) and (b) under a hypothesis that is weaker than the current
condition that $m^2/N \goto 0$.  

This setback motivated us to seek an alternate approach that would allow us to replace the probability measure $\Pnbm$, which is the restriction of the uniform 
measure $P_N$ to $\omeganbm$, with the probability measure $\Pnb$, which is the restriction of $P_N$ to $\omeganb$. The alternate approach is based on 
equation (\ref{eqn:giveitaname}) relating the probability measures $\Pnb$ and $\Pnbm$. This approach 
is successful for $b=0$ and $b=1$ in transferring to $\Pnb$ the large deviation lower bound proved in part (d) of Theorem \ref{thm:ldpthetankm} for $\Pnbm$. However, so far it has been not successful for any value of $b$ in transferring to $\Pnb$ either of the large deviation upper bounds proved in parts (b) and (c) 
of Theorem \ref{thm:ldpthetankm} for $\Pnbm$.

The starting point of the alternate approach is the following relationship between $\Pnb$ and $\Pnbm$. For $A$ any subset of $\omeganb$
\bea
\label{eqn:giveitaname}
\Pnb(A) & = & \frac{P_N(A \cap \omeganb)}{P_N(\omeganb)} \\
\nonumber 
& = & \frac{P_N(A \cap \omeganbm)}{P_N(\omeganb)} + \frac{P_N(A \cap (\omeganb \setminus \omeganbm))}{P_N(\omeganb)} \\
\nonumber
& = & \frac{\mbox{card}(\omeganbm)}{\mbox{card}(\omeganb)} \cdot \Pnbm(A) + \Pnb(A \cap (\omeganb \setminus \omeganbm)).
\eea 

Part (a) of the next theorem gives a hypothesis
that allows us to transfer the large deviation lower bound for open subsets of $\pnbc$
from $\Pnbm$ to $\Pnb$. According to part (b), this hypothesis is satisfied for $b=0$ and $b=1$. We prove part (a) after the statement of the theorem.
The proof of part (b) for $b=0$ is based on Proposition \ref{prop:omegankmagain} while the proof for $b=1$ is based on 
Proposition \ref{prop:omegankmagain} and Theorem \ref{thm:bender}.

\begin{thm}
\label{thm:transfer}
Fix a nonnegative integer $b$ and
a rational number $c \in (b,\infty)$. Let $m$ be the function $m(N)$ appearing in the definition of $\omeganbm$ in {\em (\ref{eqn:omegankm})}
and satisfying $m(N) \goto \infty$ and $m(N)^2/N \goto 0$ as $N \goto \infty$.
Let $\rhoalphabc \in \pnbc$ be the distribution having the components defined in {\em (\ref{eqn:rhoj})}. The following conclusions hold.

{\em (a)} Assume that 
\be 
\label{eqn:wantliminf} 
\lim_{N \goto \infty} \frac{1}{N} \log \left(\frac{\mbox{\em card}(\omeganbm)}{\mbox{\em card}(\omeganb)} \right) = 0.
\ee
Then for any open subset $G$ of $\pnbc$ we have the large deviation lower bound
\be 
\label{eqn:llbound}
\liminf_{N \goto \infty} \frac{1}{N} \Pnb(\omega \in \omeganb : \Thetanb(\omega) \in G) \geq -R(G | \rhoalphabc).
\ee

{\em (b)} The hypothesis in part {\em (a)} is satisfied for $b=0$ and $b=1$. Thus for these values of $b$ the large deviation lower bound
{\em (\ref{eqn:llbound})} holds.
\end{thm}

\noi 
{\bf Proof of part (a).} Let $A = \{\omega \in \omeganb : \Thetanb(\omega) \in G\}$. It follows from (\ref{eqn:giveitaname}) that 
\[
\Pnb(A) \geq \frac{\mbox{card}(\omeganbm)}{\mbox{card}(\omeganb)} \cdot \Pnbm(A).
\]
Hence by the hypothesis in part (a) and the large deviation lower bound in part (d) of Theorem \ref{thm:ldpthetankm}
\beas
\lefteqn{
\liminf_{N \goto \infty} \frac{1}{N} \log \Pnb(A)} \\ 
&& \geq  \liminf_{N \goto \infty} \frac{1}{N} \log \left(\frac{\mbox{card}(\omeganbm)}{\mbox{card}(\omeganb)}\right) + 
\liminf_{N \goto \infty} \frac{1}{N} \log \Pnbm(A) \\
&& = \liminf_{N \goto \infty} \frac{1}{N} \log \Pnbm(A) \\
&& = \liminf_{N \goto \infty} \frac{1}{N} \log \Pnbm(\omega \in \omeganbm : \Thetanb(\omega) \in G) \geq -R(G | \rhoalphabc). 
\eeas
This completes the proof of part (a). \ink

\skp

In order to prove part (b) of Theorem \ref{thm:transfer}, we now show that condition (\ref{eqn:wantliminf}) holds if $b=0$ or $b=1$. To prove this we compare the asymptotic behavior of $\mbox{card}(\Omega_{N,b,m})$
with that of $\mbox{card}(\Omega_{N,b})$ for these values of $b$. A formula for the asymptotic behavior of $\mbox{card}(\Omega_{N,b,m})$ 
for any nonnegative integer $b$ is derived in part (b) 
of Lemma \ref{lem:omegankm}. In the next proposition we express this formula in a different and more useful form for $b=0$ and $b=1$. 
Although we do not apply it here, in part (c) we give the analogous formula for $b \in \N$ satisfying $b \geq 2$. 

\begin{prop}
\label{prop:omegankmagain}
Let $b=0$ or $b=1$, and fix a rational number $c \in (b,\infty)$. Let $m$ be the function $m(N)$ appearing in the definition of $\omeganbm$ in {\em (\ref{eqn:omegankm})}
and satisfying $m(N) \goto \infty$ and $m(N)^2/N \goto 0$ as $N \goto \infty$.
Let $\alphabc$ be the quantity defined in part {\em (a)} of Theorem {\em \ref{thm:mainestimate}}.
The following conclusions hold.

{\em (a)} For $b=0$
\[
\frac{1}{N} \log \mbox{\em card}(\omeganzerom) = c \log N + \eta_N,
\]
where $\eta_N \goto 0$ as $N \goto \infty$. 

{\em (b)} For $b=1$
\[
\frac{1}{N} \log \mbox{\em card}(\omeganonem) = c \log N + (c-1) \log[c/\alphaonec] + \alphaonec - c + \eta_N, 
\]
where $\eta_N \goto 0$ as $N \goto \infty$. 

{\em (c)} For $b \in \N$ satisfying $b \geq 2$
\[
\frac{1}{N} \log \mbox{\em card}(\omeganbm) = c \log N + c\log[c/\alphabc] + \log \zalphabc - c + \eta_N.
\]
where $\zalphabc = e^{\alphabc} - \sum_{j=0}^{b-1} [\alphabc]^j/j!$ and $\eta_N \goto 0$ as $N \goto \infty$.
\end{prop}

\skp
\noindent
{\bf Proof.} We start by considering any nonnegative integer $b$. 
Let $\alpha$ be the positive real number in Lemma {\ref{lem:deltankmnu}},
and define $f(\alpha,b,c,K) = \log \zbalpha - c\log \alpha + c \log K - c$. According to part (b) of Lemma
\ref{lem:omegankm}
\[
\frac{1}{N} \log \mbox{card}(\omeganbm) =
f(\alpha,b,c,K) - \min_{\theta \in \pnbc}R(\theta | \rhoalphab) + \eta_N,
\]
where $\eta_N \goto 0$ as $N \goto \infty$.  We now appeal to item (i) in part (f) of Theorem \ref{thm:relentropy}, which shows
that 
\[
\min_{\theta \in \pnbc}R(\theta | \rhoalphab) = g(\alpha,b,c) = \log \zbalpha - c \log \alpha - (\log \zbalphac - c \log \alphabc).
\]
Substituting this formula into the preceding display, we obtain
\bea
\label{eqn:ohno}
\lefteqn{
\frac{1}{N} \log \mbox{card}(\omeganbm)} \\ 
\nonumber
& & = c \log K + \log \zalphabc - c \log \alphabc + c \log K - c + \eta_N \\
\nonumber 
& & =  c \log N + c\log[c/\alphabc] + \log \zalphabc - c + \eta_N.
\eea
where $\eta_N \goto 0$ as $N \goto \infty$.

We next use (\ref{eqn:ohno}) to prove part (a) for $b=0$ and part (b) for $b=1$. Part (c) for $b \in \N$ satisfying $b \geq 2$ is obtained by specializing 
(\ref{eqn:ohno}) to these values.

(a) As pointed out in part (a) of Theorem \ref{thm:mainestimate}, if $b=0$, then $\alpha_0(c) = c$. In this case (\ref{eqn:ohno}) becomes
\[
\frac{1}{N} \log \mbox{card}(\omeganzerom) = c \log N + \eta_N,
\]
where $\eta_N \goto 0$ as $N \goto \infty$. This completes the proof of part (a).

(b) For $b=1$, $\alphaonec$ is the unique solution in $(0,\infty)$ of the equation 
\[
c = \frac{\alphaonec Z_{0}(\alphaonec)}{Z_{1}(\alphaonec)} = \frac{\alphaonec e^{\alphaonec}}{Z_{1}(\alphaonec)}.
\] 
It follows that
\[
\log Z_1(\alphaonec) = \alphaonec + \log \alphaonec - \log c.
\]
Substituting this back into (\ref{eqn:ohno}) yields
\[
\frac{1}{N} \log \mbox{card}(\omeganonem) = c \log N + (c-1) \log[c/\alphaonec] + \alphaonec - c + \eta_N, 
\]
where $\eta_N \goto 0$ as $N \goto \infty$. The last equation coincides with the conclusion of part (b) for $b=1$. This completes the proof of the theorem.
\ink

\skp
We now prove part (b) of Theorem \ref{thm:transfer} first for $b=0$ and then for $b=1$.

\skp
\noindent
{\bf Proof of part (b) of Theorem \ref{thm:transfer} for \boldmath$b=0$\unboldmath.} 
We verify condition (\ref{eqn:wantliminf}) for $b=0$. 
According to part (a) of Proposition \ref{prop:omegankmagain}
\[
\frac{1}{N} \log \mbox{card}(\Omega_{N,0,m}) = c \log N + \eta_N, 
\]
where $\eta_N \goto 0$ as $N \goto \infty$. On the other hand, when $b=0$, $\omeganb$ equals $\Omega_N = \Lambda_N^K$. Therefore
\[
\frac{1}{N} \log \mbox{card}(\Omega_{N,0}) = \frac{1}{N} \cdot K \log N = c \log N.
\]
We conclude that 
\[
\frac{1}{N} \log \left(\frac{\mbox{card}(\Omega_{N,0,m})}{\mbox{card}(\Omega_{N,0})} \right) = \eta_N \goto 0 \mbox{ as } N \goto \infty.
\]
We conclude that condition (\ref{eqn:wantliminf}) holds for $b=0$ and thus that the large deviation lower bound
(\ref{eqn:llbound}) is valid for $b=0$.  This completes the proof. \ink

\skp
The verification of condition (\ref{eqn:wantliminf}) for $b=1$ is much deeper than that for $b=0$.

\skp
\noi 
{\bf Proof of part (b) of Theorem \ref{thm:transfer} for \boldmath$b = 1$\unboldmath.} 
This proof depends on the relationship between $\mbox{card}(\omeganone)$
and Stirling numbers of the second kind. Given $c$ a rational number in $(1,\infty)$, let $K$ and $N$ be positive integers satisfying $K/N = c$. 
We denote by $S(K,N)$ the Stirling number of the second kind, which is the number of ways to partition a set of $K$ elements into $N$ nonempty subsets
\cite[pp.\ 96--97]{Charal}. 
The $N!$ permutations of the class of all such partitions correspond to all the ways of placing the $K$ particles in the droplet model onto the $N$ sites 
of $\Lambda_N$ and therefore are in one-to-one correspondence with the elements of $\omeganone$. 
It follows that $\mbox{card}(\omeganone) = N! \cdot S(K,N)$. 

The computation of $N^{-1} \log \mbox{card}(\omeganone)$ is given in part (b) of the next theorem.
This computation is based on a deep, classical result on the asymptotic behavior of $S(K,N)$ that is
derived in Example 5.4 in \cite{Bender} and is stated in part (a) of the next theorem in our notation. The quantities in \cite{Bender} denoted by $n$,
$k$, and $r$ correspond respectively to our $K$, $N$, and $\alphaonec$. 

We now apply part (b) of Proposition \ref{prop:omegankmagain} and the conclusion of the next theorem;
the former involves the error term $\eta_N \goto 0$ as $N \goto \infty$, and the latter involves the error term $\ve_N \goto 0$ as $N \goto \infty$. 
Except for the error terms the asymptotic formulas are identical. Hence we obtain
\[
\frac{1}{N} \log \left(\frac{\mbox{card}(\Omega_{N,1,m})}{\mbox{card}(\Omega_{N,1})} \right) = \eta_N - \ve_N \goto 0 \mbox{ as } N \goto \infty.
\]
This shows that the hypothesis in part (a) of Theorem \ref{thm:transfer} is satisfied for $b=1$. The proof of part (b) of this theorem for $b=1$
will be complete after we prove the next result.

\begin{thm}
\label{thm:bender} Let
$S(K,N)$ denote the Stirling number of the second kind. Fix a rational number $c \in (1,\infty)$,
any $\delta \in (1,\infty)$, and any $M \in (\delta,\infty)$. Then 
as $K \goto \infty$ and $N \goto \infty$ with $K/N = c$
\beas
\frac{1}{N} \log \mbox{\em card}(\omeganonem) & = & \frac{1}{N} \log (N! \cdot S(K,N)) \\
& = & c \log N + (c-1) \log[c/\alphaonec] + \alphaonec - c + \ve_N, 
\eeas
where $\ve_N \goto 0$ as $N \goto \infty$. 
\end{thm}

\noi 
{\bf Proof.} We start with the asymptotic formula for $S(K,N)$ 
derived in Example 5.4 in \cite{Bender} and stated here in our notation. For any $\delta \in (0,1)$ and any $M < \infty$, uniformly for $c \in (1+\delta,M)$
the asymptotic behavior of $S(K,N)$ is given by
\[
S(K,N) = \frac{K! e^{N\alpha_1(c)}}{N! c^{N-1}\alpha_1(c)^{K-N-2} [1- c e^{-\alphaonec}] \sqrt{2\pi K}}.
\]
The quantities in \cite{Bender} denoted by $n$, $k$, and $r$ correspond respectively to our $K$, $N$, and $\alphaonec$. It follows that
\beas
\lefteqn{
\frac{1}{N} \log \mbox{card}(\omeganonem)} \\ 
& & =  \frac{1}{N} \log (N! \cdot S(K,N)) \\
& & =  \frac{K!}{N} + \alphaonec - \log c - \frac{K-N}{N} \log \alphaonec + \ve_N \\
& & = c \log N + c \log c - c + \alphaonec - \log c - (c-1) \log \alphaonec + \ve_N \\
& & = c \log N + (c-1) \log[c/\alphaonec] + \alphaonec - c + \ve_N, 
\eeas
where $\ve_N \goto 0$ as $N \goto \infty$. The proof of the theorem is complete. \ink

\skp
According to Theorem \ref{thm:transfer}, for $b=0$ and $b=1$ the large deviation lower bound, proved in part (b) of Theorem \ref{thm:ldpthetankm} for $\Pnbm$,
is also valid for $\Pnb$. Thus for any open subset $G$ of $\pnbc$ 
\be 
\label{eqn:lowerlower}
\liminf_{N \goto \infty} \frac{1}{N} \log \Pnb(\omega \in \omeganb : \Thetankm(\omega) \in G) \geq -R(G | \rhoalphabc).
\ee

For $b \in \N$ satisfying $b \geq 2$ the quantity $\mbox{card}(\omeganb)$ is related to the $b$-associated Stirling number $S_b(K,N)$ 
of the second kind by the formula
$\mbox{card}(\omeganb) = N! \cdot S_b(K,N)$. The quantity $S_b(K,N)$ is the number of ways to partition a set of $K$ elements into $N$ subsets, each of which contains at least $b$ elements \cite[pp.\ 221--222]{Comtet}. One could verify condition (\ref{eqn:wantliminf}) for these values of $b$ if there were an asymptotic formula for $S_b(K,N)$ analogous to the formula derived in Example 5.4 in \cite{Bender}. However, we are unable to locate such a formula. Nevertheless, based on our calculation for $b=0$ and $b=1$ it is reasonable to conjecture that condition (\ref{eqn:wantliminf}) holds for any $b \in \N$ satisfying $b \geq 2$, which would imply the large deviation lower bound (\ref{eqn:lowerlower}) for these values. 

We now explore whether we can extend to $\Pnb$ the large deviation upper bound proved in parts (c) and (d) of Theorem \ref{thm:ldpthetankm} for $\Pnbm$. 
If we could do this, then we could transfer to $\Pnb$ the fact,
proved in Theorem \ref{thm:equilibrium} and Corollary \ref{cor:equilibrium}, that with respect to $\Pnbm$, 
$\rhoalphabc$ is the equilibrium distribution of $\Thetanb$ and of $K_\ell$. Unfortunately, we are unable to prove the large
deviation upper bound for $\Pnb$ using either of two possible approaches
explained briefly below. 

Concerning the statement about the equilibrium distribution, the best that we can do is to use the large deviation lower bound
for $b=0$ and $b=1$ to prove that with respect to $\Pnb$ for these values of $b$, $\rhoalphabc$ is the equilibrium
distribution of $\Thetanb$ in the following weak form: for any $\ve >0$
\[
\lim_{N \goto \infty} \frac{1}{N} \log \Pnb(\omega \in \omeganb : \Thetanb(\omega) \in B_\pi(\rhoalphabc,\ve)) = 0,
\]
where $B_\pi(\rhoalphabc,\ve)$ is the open ball in $\pnbc$ with center $\rhoalphabc$ and radius $\ve$ with respect to the Prohorov metric $\pi$. 
This follows from (\ref{eqn:lowerlower}) with $G = B_\pi(\rhoalphabc,\ve)$ and from the facts that $R(B_\pi(\rhoalphabc,\ve) | \rhoalphabc) = 0$ and
\[
\limsup_{N \goto \infty} \frac{1}{N} \log \Pnb(\omega \in \omeganb : \Thetanb(\omega) \in B_\pi(\rhoalphabc,\ve)) 
\leq \limsup_{N \goto \infty} \frac{1}{N} \log 1 = 0.
\]

We end this section by discussing two possible approaches to transferring to $\Pnb$
the large deviation upper bound proved in parts (c) and (d) of Theorem \ref{thm:ldpthetankm} for $\Pnbm$. 
The first approach is based on the following upper bound valid for any subset $A$ of $\omeganb$:
\[
\Pnb(A) \leq \Pnbm(A) + \Pnb(A \cap (\omeganb \setminus \omeganbm)).
\]
This formula is a consequence of (\ref{eqn:giveitaname}) and the fact that $\mbox{card}(\omeganbm)/
\mbox{card}(\omeganb) \leq 1$. 
Now let $F$ be a compact subset of $\pnbc$, and define $A = \{\omega \in \omeganb : \Thetanb \in F\}$. The case where $F$ is a closed subset
of $\pnbc$ can be handled analogously. By part (b) of Theorem \ref{thm:ldpthetankm} 
\beas 
\lefteqn{
\limsup_{N \goto \infty} \frac{1}{N} \log \Pnb(A)} \\
&& \leq \max\!\left(\limsup_{N \goto \infty} \frac{1}{N} \log \Pnbm(A), 
\limsup_{N \goto \infty} \frac{1}{N} \Pnb(A \cap (\omeganb \setminus \omeganbm) \right) \\
&& 
\leq \max\!\left(-R(F | \rhoalphabc), 
\limsup_{N \goto \infty} \frac{1}{N} \Pnb(A \cap (\omeganb \setminus \omeganbm) \right).
\eeas
If we could prove that $-R(F | \rhoalphabc)$ is greater than or equal to the second expression on the right side of the last line, then 
we would be able to transfer the large deviation upper bound to $\Pnb$. Unfortunately, however, we are unable prove that 
$-R(F | \rhoalphabc)$ is greater than or equal to the second expression on the right side of the last line. 

The second approach to transferring to $\Pnb$ the large deviation upper bound in parts (c) and (d) of Theorem \ref{thm:ldpthetankm} rests on a careful
analysis of how these upper bounds follow from the local estimate in part (b) of Theorem \ref{thm:mainestimate} and from Theorem \ref{thm:balls} as 
applied to Theorem \ref{thm:ldlimitballs}, for which we need only the large deviation upper bound for the sets appearing in Theorem \ref{thm:ldlimitballs}. 
Omitting the details, we claim that the crucial step is to show that
\[
\lim_{N \goto \infty} \min_{\nu \in \anbm} R(\thetanbnu | \rhoalphab) = \min_{\theta \in \pnbc} R(\theta | \rhoalphab).
\]
At the end of the proof of part (b) of Lemma \ref{lem:omegankm} we prove this limit by applying the approximation procedure in appendix B, which 
requires the condition that $m^2/N \goto 0$ as $N \goto \infty$. If we could prove this limit without invoking the approximation procedure and under a condition that is compatible with $m/\sqrt{N} \goto \infty$ as $n \goto \infty$, then the large deviation upper bound in parts (c) and (d) of Theorem \ref{thm:ldpthetankm} would hold with $\Pnb$ replacing $\Pnbm$. Unfortunately, we have not been able to carry this out.

We end this section by proposing an interesting test case for gaining insight into whether the conditioned measure $\Pnbm$ could be replaced by $\Pnb$ in the LDP for $\Thetanb$ in Theorem \ref{thm:ldpthetankm}. This test case would be to use the methods of this paper to prove 
Sanov's Theorem for the empirical measures of i.i.d.\ random variables taking values 
in $\Nb$. This theorem, of course, can be proved directly without the methods of this chapter  \cite[Thm.\ 6.2.10]{DemboZeitouni}, \cite[Thm.\ 4.5]{DonVar3}. If one uses the methods of this paper, then one would first have to prove it for the analogue of the measure $\Pnbm$ restricted to the analogue of the restricted configuration space $\omeganbm$, where the number of positive components of $N_j$ is restricted by $m = m(N)$. The quantity $m(N) \goto \infty$ at an appropriate rate. It would be instructive to see if this restriction can be eliminated using one of the approaches proposed in this appendix.

\end{document}